\documentclass[11pt,a4paper]{article}
\usepackage[utf8]{inputenc}
\usepackage[T1]{fontenc}
\usepackage{verbatim} 
\usepackage[section]{placeins}
\usepackage{url}
\usepackage[sc,osf]{mathpazo}
\usepackage{enumitem,textcomp}
\usepackage{linguex}
\usepackage[english]{babel}
\usepackage{geometry}

\newcommand{\amrand}[2]
  {\leavevmode
   \marginpar
     [\raggedleft\scriptsize \leavevmode\normalcolor #1: #2]
     {\raggedright\scriptsize \leavevmode\normalcolor #1: #2}
  }

\usepackage{amsmath,enumitem}
\usepackage{amssymb}
\usepackage{amsthm}
			\newtheorem{thm}{Theorem}[section]

			\newtheorem{prop}[thm]{Proposition}

			\newtheorem{defn}[thm]{Definition}
			\newtheorem{fact}[thm]{Fact}

\usepackage{tikz}

\usepackage{units}
\usetikzlibrary{arrows,decorations.pathmorphing,
backgrounds,positioning,fit,matrix}
\usepackage{tikz-qtree}
\usepackage{tikz-qtree-compat,textcomp}
\usetikzlibrary{calc}
\usetikzlibrary{backgrounds}
\usetikzlibrary{shapes.geometric}
\usetikzlibrary{shapes.arrows}
\usepackage{units}
\usepackage{booktabs} 
\usepackage{graphicx}
\usepackage{bussproofs} 
\usepackage{proof} 
\usepackage{eso-pic} 

\usepackage{color, colortbl}
\definecolor{Gray}{gray}{.85}
\definecolor{Blue}{RGB}{60,20,255}
\definecolor{Red}{RGB}{255,20,60}

\usepackage{titling}
\pretitle{\begin{center}\huge}

\usepackage{hyperref}
\hypersetup{colorlinks=true,allcolors=blue}

\usepackage[sectionbib,longnamesfirst]{natbib}
\setcitestyle{aysep={}} 
\shortcites{cobreros2015vagueness,Cobreros2012tolerant,OverEtAl2007,EvansEtAl2007}

\newcommand{\cctt}[0]{{\sf CC/TT}}
\newcommand{\dftt}[0]{{\sf DF/TT}}
\newcommand{\half}{\nicefrac{1}{2}}  

\newcommand{\Lc}{\mathcal{L}_{\to}}  

\newcommand*{\prob}{\mathsf{Pr}}

\vspace{1cm}
\author{}%
\date{}
\title{De Finettian Logics of Indicative Conditionals
}

\begin{document}

\maketitle
\thispagestyle{empty}
\vspace{-2.5cm}
\centerline{\large Part I: Trivalent Semantics and Validity}\vspace{1cm}

\centerline{\small Paul \'Egr\'e, Lorenzo Rossi, Jan Sprenger}

%

\begin{abstract}

\noindent This paper explores \textit{trivalent truth conditions} for indicative conditionals, examining the ``defective'' table put forward by \citet{definetti1936logique}, as well as \citet{reichenbach1944philosophic}, first sketched in \cite{reichenbach1935wahr}. On their approach, a conditional takes the value of its consequent whenever its antecedent is True, and the value Indeterminate otherwise. Here we deal with the problem of choosing an adequate notion of validity for this conditional. We show that all standard trivalent schemes are problematic, and highlight two ways out of the predicament: one pairs de Finetti's conditional ({\sf DF}) with validity as the preservation of non-False values ({\sf TT}-validity), but at the expense of Modus Ponens; the other modifies de Finetti's table to restore Modus Ponens.
In Part I of this paper, we present both alternatives, with specific attention to a variant of de Finetti's table ({\sf CC}) proposed by \citet{cooper1968propositional} and \citet{cantwell2008logic}. In Part II, we give an in-depth treatment of the proof theory of the resulting logics, {\sf DF/TT} and {\sf CC/TT}: both are connexive logics, but with significantly different algebraic properties.

\end{abstract}


\section{Introduction}

Choosing a semantics for the indicative conditional of natural language ``if $A$, then $C$'' (henceforth, $A \to C$) usually involves substantial tradeoffs. The most venerable account identifies the indicative conditional with the material conditional $\neg A \vee C$ and has several attractive features: it is truth-functional, allows for a straightforward treatment of nested conditionals, and satisfies various intuitive principles such as Conditional Proof and Import-Export. However, the material conditional account severs the link between antecedent and consequent: Suppose John was not in Paris yesterday. Then ``if John was in Paris yesterday, then he will be in Milan tomorrow'' is true regardless of John's travels plans. The inferential dimension of conditionals, and in particular the link between truth and justified assertion, are completely lost in this picture. 

Seeking a way out of this predicament, \citet{Stalnaker1968,stalnaker1975indicative} proposed to give up truth-functionality and to strengthen the truth conditions of the indicative conditional as follows: $A \to C$ is true if and only if $C$ is true in the closest possible $A$-world---i.e. the closest world in which the antecedent is true. This proposal has many virtues but also some limitations, on which we say more in the next section. 

A second strategy admits that the truth conditions of the indicative conditional may not be truth-functional, or perhaps agree with those of the material conditional \citep[e.g.,][]{Jackson1987}, but in any case they are a matter of secondary importance. What matters, ultimately, is the assertability or ``reasonableness'' of a conditional $A \to C$. This notion is often explicated in probabilistic terms, by analyzing the conditional as expressing a \textit{supposition} that precedes the evaluation of the consequent, and by focusing on the probability of $C$ given $A$, in symbols $\prob(C|A)$. This strategy is popular among cognitive scientists \citep[e.g.,][]{EvansEtAl2007,OverEtAl2007}, and among philosophers who focus on the evidential and inferential dimension of a conditional \citep[e.g.,][]{Adams1965,Adams1975,Edgington1995,Krzyzanowska2015,Douven2016}. To our mind, however, it would be preferable to have a theory that \textit{explains} how assertability conditions are related to, and can be motivated from, the truth conditions of a conditional. 

This paper is an attempt to connect the dimensions of truth and assertability in a principled way, and to construct a semantics that preserves the most attractive features of both propositional and non-propositional accounts. Due to well-known impossibility results \citep[e.g.,][]{gibbard1980two,Lewis1986}, this means that we have to leave the familiar framework of bivalent logic. Our starting point is the intuition voiced by \citet{definetti1936logique}, \citet[168]{reichenbach1944philosophic} and \citet{quine1950methods} (crediting Ph.~Rhinelander for the idea), that uttering a conditional amounts to making a conditional assertion: the speaker is committed to the truth of the consequent when the antecedent is true, but committed to neither truth nor falsity of the consequent when the antecedent is false.

\begin{figure}[t]
\[
\begin{tabular}[t]{c|cc}
 & 1 & 0\\
 \hline
 1 & 1 & 0\\
 0 & \# & \#\\
 \end{tabular}
 \qquad
\begin{tabular}[t]{c|ccc}
 & 1 & \# & 0\\
 \hline
  1 & 1 & $\cdot$ & 0\\
 \# & $\cdot$ & $\cdot$ & $\cdot$\\
  0 & \# & $\cdot$ & \#\\
  \end{tabular}
 \]
\caption{``Defective'' bivalent table (left) and trivalent incomplete expansion (right)}\label{fig:defective}
 \end{figure}

The idea that a conditional with a false antecedent has no classical truth value is sometimes summarized in what \citet{kneale1962development} have named the ``defective'' truth table, where the symbol `\#' marks a truth value gap (Figure \ref{fig:defective}), and whose first appearance may be found in \citet[381]{reichenbach1935wahr}.\footnote{De Finetti presented his paper in Paris in the same year 1935, with explicit reference to \citealt{reichenbach1935wahr}, but criticizing the latter's objective interpretation of probability. To the best of our knowledge, Reichenbach's 1935 book does not quite present de Finetti's three-valued table, but some variants instead. However \citet[168, fn.2]{reichenbach1944philosophic} traces quasi-implication back to his previous opus. In our view, the de Finetti conditional may therefore be called the de Finetti-Reichenbach conditional, but for simplicity and partly for established usage, we stick to calling it the {\sf DF} conditional. See also \citealt{milne1997bruno}, \citealt{baratgin2013uncertainty}, and \citealt{over2017defective} on the history of the defective table.}   When the gap is handled as a value of its own (we represent it by $\nicefrac{1}{2}$, for ``indeterminate''), and so as a possible input for semantic evaluation, then the ``defective'' two-valued conditional naturally leads to truth conditions within a \emph{trivalent} (=three-valued) logic. For de Finetti, asserting a conditional of the form ``if $A$ then $C$'' is a \textit{conditional assertion}: an assertion that is retracted, or \emph{void}, if the antecedent turns out to be false. In this respect, it is akin to making a bet that if $A$ then $C$. When $A$ is realized and $C$ is false, the bet is lost; when $A$ is realized and $C$ is true, the bet is won; when $A$ is not realized, however, the bet is simply called off (more on this in Section 2). The trivalent table proposed by de Finetti for the conditional (henceforth $\rightarrow_{\sf DF}$) is given in Figure \ref{fig:DF}. The same table is put forward by \citet{reichenbach1944philosophic}, who calls it \emph{quasi-implication}. Like de Finetti, Reichenbach considers that some conditionals are void when the antecedent is false, though Reichenbach's interpretation of the third truth value differs, being driven by measurement-theoretic considerations in quantum physics.\footnote{Closer to the interpretation of the third truth value that features in \citet{bochvar1937three}, Reichenbach considers that some conditionals are \emph{meaningless} when the antecedent concerns an event whose precise measurement is impossible (for instance, we cannot in general simultaneously measure position and momentum of a particle with arbitrary degree of precision). 
Reichenbach treats the third truth value as objectively indeterminate rather than as expressing a notion of subjective ignorance, as de Finetti does. In motivating this interpretation, Reichenbach refers explicitly to the Bohr-Heisenberg interpretation of quantum mechanics. }



\begin{figure}[h!]\label{fig:finetti}
\[
\begin{tabular}[t]{l|ccc}
$f_{\to_{\sf DF}}$ & 1 & $\nicefrac{1}{2}$ & 0\\
 \hline
 1 & 1 &  $\nicefrac{1}{2}$ & 0\\
$\nicefrac{1}{2}$ & $\nicefrac{1}{2}$& $\nicefrac{1}{2}$& $\nicefrac{1}{2}$\\
 0 & $\nicefrac{1}{2}$ &  $\nicefrac{1}{2}$& $\nicefrac{1}{2}$\\
 \end{tabular}
 \]\caption{The truth table for de Finetti's trivalent conditional.}\label{fig:DF}
 \end{figure}

The truth table given by de Finetti and Reichenbach mirrors an interpretation on which the conditional is indeterminate when its antecedent is \emph{not true} ($\neq 1$).\footnote{Related proposals include \citealt{jeffrey1963indeterminate}, \citealt{belnap1970conditional}, \citealt{manor1975propositional}, \citealt{farrell1986implication}, \citealt{dubois1994conditional}, \citealt{mcdermott1996truth}, \citealt{rothschild2014capturing}, and \citealt{kapsner2018humble}. We note that the {\sf DF} table was reintroduced several times in the past decades, very often without prior notice of either de Finetti or Reichenbach, and sometimes with separate motivations in mind, viz. \citet{blamey1986partial}, who calls it \emph{transplication}, to highlight its hybrid character between a conjunction and an implication, 
or recently \citet{kapsner2018humble}, who came up with the scheme specifically to deal with connexiveness. More on this will be said below.}
 However, understanding the conditional as a conditional assertion is compatible with more choices of a trivalent truth table, especially in the second line, namely when the antecedent is indeterminate (viz.~the antecedent might be a conditional with false antecedent). Two notable proposals, the first due to \citet{cooper1968propositional} and independently to \citet{cantwell2008logic}, the second to \citet{farrell1979material}, are given in Figure \ref{fig:cooper}:
\begin{figure}[htb]
\[
\begin{tabular}[t]{l|ccc}
$f_{\to_{\sf CC}}$ & 1 & $\nicefrac{1}{2}$ & 0\\
 \hline
 1 & 1 &  $\nicefrac{1}{2}$ & 0\\
$\nicefrac{1}{2}$ & 1& $\nicefrac{1}{2}$& 0\\
 0 & $\nicefrac{1}{2}$ &  $\nicefrac{1}{2}$& $\nicefrac{1}{2}$\\
 \end{tabular}
 \qquad
 \begin{tabular}[t]{l|ccc}
$f_{\to_{\sf F}}$ & 1 & $\nicefrac{1}{2}$ & 0\\
 \hline
 1 & 1 &  $\nicefrac{1}{2}$ & 0\\
$\nicefrac{1}{2}$ & $\nicefrac{1}{2}$& $\nicefrac{1}{2}$& 0\\
 0 & $\nicefrac{1}{2}$ &  $\nicefrac{1}{2}$& $\nicefrac{1}{2}$\\
 \end{tabular}
 \]\caption{Truth tables for the Cooper-Cantwell conditional (left) and the Farrell conditional (right).}\label{fig:cooper}
 \end{figure}

Which trivalent table is the most adequate? \citet{baratgin2013uncertainty} approached this question experimentally. They asked participants to evaluate various indicative conditional sentences as ``true'', ``false'' and ``neither'', by manipulating the truth value of the antecedent and consequent (making them clearly true, false, or uncertain). They conclude that the original de Finetti table is better-supported than its competitors and that participants's judgments are well-correlated with the de Finettian bet interpretation of conditionals. However, the focus on (intuitions about) truth tables neglects the inferential properties of conditionals, that is, how we should reason with them. For that, we need an analysis of the notion of logical validity. Indeed, the same truth tables can support radically distinct entailments, depending on how validity is defined. 

 
 In trivalent logic, several notions of validity can be considered, and they yield significantly distinct predictions (\citealt{egre2016conditionals}). 
 Consider validity as preservation of truth (i.e., the value 1) from premises to conclusion in an argument. Following the terminology of \citet{Cobreros2012tolerant}, we call this \emph{strict-to-strict} validity, or {\sf SS}-validity. %
 An alternative is to define validity as the preservation of non-falsity ($\{1, \nicefrac{1}{2}\}$), also known as  \emph{tolerant-to-tolerant} or {\sf TT}-validity. Other schemes considered in the literature are the intersection of {\sf SS} and {\sf TT} (see \citealt{mcdermott1996truth}), as well as  so-called \emph{mixed} (strict-to-tolerant, tolerant-to-strict) consequence relations ({\sf ST, TS}). All schemes have advantages and drawbacks, but some combinations of a conditional operator with a validity scheme appear better than others.

In this paper, we bring together the research strands on validity in trivalent logic and trivalent semantics for indicative conditionals. More precisely, we conduct a systematic investigation of the main trivalent semantics for defective conditionals, and isolate the most promising combinations of truth tables and validity relations. To the best of our knowledge, no such systematic comparison has been conducted so far. In particular, apart from \citet{cooper1968propositional}, we are not aware of an axiomatization of the logics based on a trivalent semantics for the indicative conditional. 

We fill this gap in our paper and proceed in two main parts. Part I of this paper focuses on semantics: it reviews the main motivations for the de Finetti conditional (Section \ref{S:2}) and expounds the problems it faces when selecting an adequate trivalent consequence relation. This is what we call the ``validity trilemma'' for the de Finetti conditional (Section \ref{S:3}): the de Finetti conditional must either fail to support any sentential validity, support unacceptable arguments, or fail Modus Ponens. We present two ways out of this predicament: the first bites the bullet and associates de Finetti's conditional with a notion of tolerant-to-tolerant validity that fails Modus Ponens (Section \ref{sec:DF+MP}). The other consists in modifying de Finetti's table so as to restore Modus Ponens for the same notion of validity. We specify the class of trivalent conditionals that support Modus Ponens and are adequate for {\sf TT}-validity (``Jeffrey conditionals"), and we distinguish, among those, the conditional introduced independently by Cooper and Cantwell (Section \ref{S:5}). We end part I of this paper with a comparison between the two logics that ensue from those considerations, {\sf DF/TT} (de Finetti-{\sf TT}) and {\sf CC/TT} (Cooper-Cantwell-{\sf TT}), with an indication of their commonalities (in particular both are connexive logics, Section \ref{S:5}) and limitations (in particular both retain the Linearity principle of two-valued logic, see Section \ref{S:6}).  
In part II, we further this comparison with an in-depth investigation of the proof theory and algebraic properties of those two logics.

\section{The de Finetti Conditional}\label{sec:definetti}
\label{S:2}

\subsection{Philosophical Motivation}
\label{S:2.1}

\citet{Ramsey1926} was likely the first philosopher to connect an assertion of a proposition $A$ with an implicit disposition to bet on $A$, and to interpret an indicative conditional $A \to C$ as a conditional assertion where we suppose the antecedent, and reason on that basis about the consequent. His views strongly influenced de Finetti, who combined both ideas of Ramsey's by postulating an isomorphism between the conditions that settle the truth of a (conditional) proposition, and the conditions that settle the winner of a (conditional) bet. Evaluating the truth or falsity of a conditional proposition, assertion or event requires supposing the antecedent in the same way that a conditional bet on $C$ given $A$ can only be won or lost if $A$ is true; if $A$ is false, the bet will be called off. 

Hence, while the truth value of an ordinary, non-conditional proposition $A$ is settled by either $A$ or $\neg A$, the truth value of a \textit{conditional proposition or assertion}---de Finetti uses the notation $C/A$---is settled by the corresponding pair $A \wedge C$ and $A \wedge \neg C$ \citep[][568, emphasis in original]{definetti1936logique}:
\footnote{
In the English translation of R.~Angell, the quote goes: ``It is here that introduction of a special logic of three values seems indicated, as we have already announced: $C$ and $A$ being any two events (propositions) whatever, we will speak of the \emph{tri-event} ${C}/
{A}$ ($C$ given $A$), the logical entity which is considered:
\begin{enumerate}
\item \emph{true} if $C$ and $A$ are true;
\item \emph{false} if $C$ is false and $A$ true; 
\item \emph{null} if $A$ is false
\end{enumerate}
(one
does not distinguish between ``not $A$ and $C$" and ``not $A$ and not $C$", the
tri-event being only a function of $A$ and $A\wedge C$).''
}
\begin{quote}
``C'est ici qu'il para\^it indiqu\'e d'introduire une logique sp\'eciale \`a trois valeurs, comme nous l'avions d\'ej\`a annonc\'e : $C$ et $A$ \'etant deux \'ev\'enements (propositions) quelconques, nous dirons \textit{tri\'ev\'enement} $C$/$A$ ({$C$} subordonn\'e \`a {$A$}), l'entit\'e logique qui est consid\'er\'ee 

\begin{enumerate}
  \item \emph{vraie} si $C$ et {$A$} sont vrais;
  \item \emph{fausse} si $C$ est faux et $A$ est vrai; 
  \item \emph{nulle} si $A$ est faux 
\end{enumerate}
(on n'a pas de distinction entre ``non {$A$} et {$C$}'' et ``non {$A$} et non {$C$}'', le tri\'ev\'enement ne devant \^etre fonction que de {$A$} et ${C \wedge A}$).''

\end{quote}

This approach explains the intuition that upon observing $A \wedge C$, we feel compelled to say that the (previously made) conditional assertion $C/A$ was \emph{right}, that it has been \emph{verified}.\footnote{See also \citealt{cantwell2008logic}, and the ``hindsight problem'' in \citealt{khoo2015indicative}.} Similarly, the conditional assertion $C/A$ is \emph{falsified} by the observation of $A \wedge \neg C$: we have been proven \emph{wrong} by the facts. The indicative conditional $A \to C$ shall, in the rest of this paper, be understood as a conditional assertion $C/A$ whose truth conditions correspond to the conditions that determine the result of a conditional bet. We now define a corresponding class of conditional operators:

\begin{defn}[de Finettian operators]
A trivalent binary operator is called \emph{de Finettian} if it agrees with de Finetti's truth conditions when the antecedent is determinate, that is, when the antecedent takes the value 1 or the value 0.
\label{D:definettian}
\end{defn}

Equivalently, an operator is de Finettian if it agrees on the first and third row of the table in Figure \ref{fig:DF}: it takes the value indeterminate when its antecedent is false, and the value of its consequent when its antecedent is true. From the class of de Finettian operators, de Finetti selects the truth conditions that assign value $\nicefrac{1}{2}$ to the conditional whenever the antecedent is itself indeterminate. Note that this grouping of indeterminate with false antecedents is not covered by the above epistemological motivation; in fact, this choice is a classical point of contention between trivalent logics of conditionals. De Finetti's choice resembles Bochvar's scheme for trivalent operators (a.k.a. the Weak Kleene scheme), where the value $\half$ is carried over from any part of a sentence to the whole sentence (\citealt{bochvar1937three}). Similarly, he assumes that a conditional is undefined as soon as antecedent or consequent are undefined. As we know from the theory of presupposition projection (\citealt{beaver2001partial}), however, Bochvar's choice is not the most adequate to account for the transmission of indeterminate values from smaller to larger constituents, and therefore it should not be viewed as mandated by the rest of de Finetti's motivations for the conditional. In fact, de Finetti himself does not handle conjunction and disjunction \`a la Bochvar/Weak Kleene, but in line with the Strong Kleene scheme (see below).

\subsection{Main benefits of the approach}\label{sec:benefits}
\label{S:2.2}

De Finetti's trivalent approach has the potential to avoid the paradoxes of material implication and yields a variety of benefits.\footnote{In particular, paired with {\sf SS}-validity, the de Finetti conditional supports neither the entailment from $\neg A$ to $ (A \to C)$, nor the entailment from $C$ to$ (A \to C)$. For {\sf TT}-validity, only the former scheme is blocked.} First of all, it is very simple and has a clear motivation: asserting a conditional amounts to making a conditional assertion; conditionals express dispositions to bet just as ordinary assertions do. The trivalent approach treats conditionals as expressing propositions, in agreement with their linguistic form and assertive usage; only their truth conditions cannot be expressed in bivalent logic. This is a substantial advantage over non-propositional views that have to explain the gap between linguistic form and philosophical theorizing.  

Second, de Finettian conditionals keep the epistemic notion of assertability and the semantic notion of truth separate, while allowing for a fruitful interaction: degrees of assertability can be defined directly in terms of the truth conditions. For a probability function $\prob$ on a propositional language, and assuming $X$ is a Boolean sentence or a simple conditional,\footnote{We refer to section \ref{S:6} for a discussion of compounds of conditionals.} we define the degree of assertability to be:

\begin{equation}\tag{A}\label{eq:ast}
 {\sf Ast}(X) = 
                \prob(X \; \text{is true}|X \; \text{has a classical truth value})
                \end{equation}
(see also \citealt{mcdermott1996truth, cantwell2006laws, rothschild2014capturing}). Trivalent semantics replaces the familiar norm of asserting what is probably true by the equally plausible norm of asserting what is (much) more likely to be true than to be false. This collapses to the classical picture ${\sf Ast}(X) = \prob(X \; \text{is true})$ for bivalent propositions. For $X = A \to C$, and assuming that $C$ is bivalent, we obtain 
\begin{eqnarray*}
{\sf Ast}(A \to C) &=& \prob(A \to C \; \text{is true}|A \to C \; \text{has a classical truth value}) \\
&=& \prob(A \wedge C\ \text{is true} |A\ \text{is true}) \\
&=&\prob(C|A)
\end{eqnarray*}

We thus obtain \emph{Adams' Thesis} (sometimes also called ``The Equation'', and read as a thesis about the \emph{probability} of $A \to C$), a plausible principle for the assertability of conditionals supported by patterns observed in natural language \citep[][]{Stalnaker1968,Adams1975,dubois1994conditional,EvansEtAl2007,OverEtAl2007,egre2011if,Over2016}.\footnote{For recent criticisms of Adams' Thesis, see \citealt[][]{DouvenVerbrugge2010} and \citealt{SkovgaardOlsenEtAl2016a}.} Similarly, the suppositional reading of conditionals as expressing conditional degrees of belief \citep[e.g.,][]{Ramsey1926,Edgington1995} can be naturally grounded in trivalent semantics. Other theories, notably \cite{spohn2013ranking}'s, accept Adams thesis, but view is as dependent on a more fundamental notion of conditional belief, captured by rank-ordering instead of probability. That theory too is compatible with de Finetti's trivalent approach.\footnote{Compare \citet[1083]{spohn2013ranking}'s definition of conditional rank with the probabilistic derivation of {\sf Ast} above. The conditional rank of $C$ given $A$ is $\kappa(C|A)=\kappa(A\cap C)-\kappa(A)$, provided $\kappa(A)< \infty$. As highlighted by Spohn, a bridge is given by a logarithmic transformation.}

The close relationship between truth and assertability allows us to explain intuitions which conflict at first with the trivalent view. For example, a sentence such as:

\ex. If Paul is in Paris, then Paul is in France.

would typically be judged as true, whereas trivalent semantics regard this as an empirical question: when Paul is in Berlin, the sentence has indeterminate truth value. However,  the trivalent view can offer an error theory since \Last is \textit{maximally assertable} regardless of Paul's whereabouts ($\prob(C|A)=1$). When we call sentences such as \Last ``true'', what we really mean is that they command consent, that they are ``maximally assertable'' \citep[see also][]{Adams1975}. Since assertability conditions are fully defined in terms of truth conditions, this defense is arguably not \textit{ad hoc}. In sum, on this view, indicative conditionals are \textit{factual}---their truth and falsity is a matter of correspondence with the world---, like for predictions about future events, while their assertability is epistemic and is represented probabilistically.

Thirdly, the de Finetti conditional satisfies the following identity:
\begin{equation}\label{eqn:IE}
A \to (B \to C) \equiv (A  \wedge B) \to C \tag{Import-Export}
\end{equation}
Here, ``$\equiv$'' means that the truth values of $A \to (B \to C)$ and $(A  \wedge B) \to C$ coincide according to the de Finetti tables. \ref{eqn:IE} expresses the idea that right-nesting a conditional is just the same as adding a further supposition. \citet{gibbard1980two} proved that there is no truth-conditional operator $\to$ 
that (i) satisfies \ref{eqn:IE}; (ii) validates $A \to C$ whenever $A$ classically entails $C$; (iii) is strictly stronger than the material conditional. In Stalnaker's and Lewis's possible world semantics, \ref{eqn:IE} thus fails. \citet{mcgee1989conditional} proposed a modification of Stalnaker's semantics that restores \ref{eqn:IE} and is stronger than the material conditional, giving up (ii).\footnote{See \citet{McGee1985, mcgee1989conditional} on the failure of Modus Ponens in that logic. \citet{Mandelkern2019a} observes a certain tension between Import-Export and classical conjunction, suggesting to restrict Import-Export accordingly. However, our findings show that the canonical extension of classical conjunction to trivalent logics (i.e., Strong Kleene truth tables) is perfectly compatible with Import-Export. The observed tensions may therefore be a peculiar feature of bivalent logic.} However, it involves syntactic restrictions on the sentences appearing as antecedents. 
{The advantage of de Finetti's conditional is that it can satisfy Import-Export without any syntactic restrictions, and within a truth-conditional framework. Depending on which notion of validity it is paired with, it may or may not obviate the conditions of Gibbard's theorem. As we will see, however, even when it falls prey to Gibbard's result, it need not have all properties of its material counterpart.}



\section{Comparing Schemes for Validity}
\label{S:3}

We now introduce and compare the main notions of validity that can be used in relation to de Finetti's conditional. By so doing, we expose a problem for the de Finetti conditional: all of the basic schemes available for validity in trivalent logic appear to overgenerate or to undergenerate relative to general principles of conditional reasoning.

\subsection{Evaluations and Validity}
\label{S:3.1}

Throughout the paper, we let $\mathcal{L}$ be a propositional language featuring denumerably many propositional variables (indicated as $p_0$, $p_1$, $\ldots$), whose logical connectives include $\neg$ and $\wedge$ (the others, $\vee$ and $\supset$, are defined as usual). We call $\Lc$ the language obtained from $\mathcal{L}$ by adding a new conditional connective, in symbols $\to$, to the primitive stock of logical constants of $\mathcal{L}$. We use uppercase Latin letters ($A$, $B$, $C$, $\ldots$) as meta-variables for $\mathcal{L}$- and $\Lc$-sentences, and {\sf For} to denote the set of formulae of the language $\Lc$. With a slight notational abuse, we will write $\Gamma, A$ rather than $\Gamma \cup \{A\}$ (for $\Gamma$ a set of $\Lc$-formulae and $A$ a $\Lc$-formula), in order to improve readability. 

For all trivalent semantics of the conditional that we consider, negation and conjunction are interpreted via the familiar Strong Kleene truth tables (introduced by \citealt{lukasiewicz1920three}, also featuring in \citealt{definetti1936logique}): 

\begin{figure}[h!]
\[
\begin{tabular}{|c|c|}
\hline
& $f_{\neg}$\\
\hline
$1$ & $0$\\
$\nicefrac{1}{2}$ & $\nicefrac{1}{2}$\\
$0$ & $1$\\
\hline
\end{tabular}
\hspace{10pt}
\begin{tabular}{|c|ccc|}
\hline
$f_{\wedge}$ & $1$ & $\nicefrac{1}{2}$ & $0$\\
\hline
$1$ & $1$ & $\nicefrac{1}{2}$ & $0$\\
$\nicefrac{1}{2}$ & $\nicefrac{1}{2}$ & $\nicefrac{1}{2}$ & $0$\\
$0$ & $0$ & $0$ & $0$\\
\hline
\end{tabular}
\]
\caption{Truth tables for negation and conjunction}
\end{figure}

We can now proceed to define evaluations and consequence relations for the de Finetti conditional. %

\begin{defn}[Classical, {\sf SK}-, and {\sf DF}-evaluation]
\item[-] A \textnormal{classical} evaluation is a function from $\mathcal{L}$-sentences to $\{1, 0\}$ that interprets $\neg$ and $\wedge$ by the functors $f_{\neg}$ and $f_{\wedge}$ restricted to the values $1$ and $0$.
\item[-] A \textnormal{Strong Kleene evaluation} (or {\sf SK}-evaluation) is a function from $\mathcal{L}$-sentences to $\{1, \nicefrac{1}{2}, 0\}$ that interprets $\neg$ and $\wedge$ by the functors $f_{\neg}$ and $f_{\wedge}$.
\item[-] A \textnormal{de Finetti} evaluation (or {\sf DF}-evaluation) is a function from {\sf For} to $\{1, \nicefrac{1}{2}, 0\}$ interpreting $\neg$, $\wedge$, and $\to$ by the functors $f_{\neg}$, $f_{\wedge}$ and $f_{\to_{\sf DF}}$.
\end{defn}

Given an evaluation, we can distinguish {two levels of \emph{truth}} for a sentence, namely {\sf T}-truth (for tolerant truth) and {\sf S}-truth (for strict truth), following \citealt{Cobreros2012tolerant} and \citealt{cobreros2015vagueness}.\footnote{\citet{zardini2008model} talks of levels of \emph{goodness} for a sentence, and \citet{cobreros2015vagueness} talk of levels of \emph{assertability}, rather than truth. Given our separation of truth and assertability in the previous section, we avoid this terminology.} Identifying the value 1 with the True, the value $\half$ with the Indeterminate, and the value 0 with the False, then {\sf S}-truth is for a sentence to be True, whereas {\sf T}-truth is for a sentence to be non-False. The two notions obviously coincide relative to classical evaluations, but they come apart relative to trivalent evaluations.

\begin{defn}[{\sf T}-truth and {\sf S}-truth]$ $
\begin{itemize}
\item[-] An evaluation $v: {\sf For} \longmapsto \{1, \half, 0\}$ makes a sentence $A$ \textnormal{strictly true (or {\sf S}-true)} provided $v(A)=1$. 
\item[-] An evaluation $v: {\sf For} \longmapsto \{1, \half, 0\}$ makes a sentence $A $ \textnormal{tolerantly true (or {\sf T}-true)} provided $v(A)>0$. 
\end{itemize}
\end{defn}

Following \citet{chemla2017characterizing} and \citet{chemla2018suszko}, we single out five notions of validity in a trivalent setting, depending on whether validity is defined as the preservation of truth, non-falsity, or as some combination of those. Those five notions of validity are not the only conceivable ones in trivalent logic, but there is a sense in which they form a natural class.\footnote{See \citet{chemla2017characterizing} for general arguments regarding the oddness of $\mathsf{SS\cup TT}$ in particular. In the present case, taking the union of $\mathsf{SS}$ and $\mathsf{TT}$ would obviously not solve the overgeneration problem raised in the next section, in particular regarding the entailment to the converse conditional. \citet{cooper1968propositional} restricts $\mathsf{TT}$ to bivalent atomic valuations -- what \citet[][\S 7.19, 1044 and following]{humberstone2011conn} calls `atom-classical' valuations: we set aside that restriction, which makes no essential difference to our discussion here. \citet{farrell1979material} sketches another variant, which we can set aside on the same grounds (see next footnote).} In particular, the five schemata under discussion are all monotonic, and they are all the monotonic trivalent schemata (see \citealt{chemla2018suszko} for a proof), meaning that an inference remains valid by the inclusion of more premises. We leave open whether a nonmonotonic scheme for validity might offer a good fit for the original de Finetti table.\footnote{\citet{farrell1979material} introduces a notion of sentential validity that may be generalized into a nonmonotonic notion of argument-validity. On his definition, $A$ is valid provided it is {\sf TT}-valid, and there is a valuation that gives $A$ the value 1. We may generalize this to: $\Gamma\models A$ provided $\Gamma$ {\sf TT}-entails $A$ and there is at least one valuation that gives the formulae in $\Gamma$ and $A$ the value 1. On that definition, $A\models A$, but $A, \neg A\not \models A$ (we are indebted to a remark by T. Ferguson in relation to that fact). We note that like standard {\sf TT}-validity, this nonmonotonic restriction still fails Modus Ponens. As such, it would not add a separate route from the one described with standard {\sf TT}-validity.}
 
\begin{defn}[{\sf SS}-, {\sf TT}-, {\sf SS$\cap$TT}-, {\sf ST}- and {\sf TS}-validity]

For every $\{\Gamma, A\} \subseteq {\sf For}$, for every $\sf X$-evaluation, we say that: 
\begin{itemize}
\item[-] {$\Gamma \models_{{\sf X/SS}} A$, provided every $\sf X$-evaluation that makes all sentences of $\Gamma$ $\mathsf{S}$-true also makes $A$ $\mathsf{S}$-true.}

\item[-] {$\Gamma \models_{{\sf X/TT}} A$, provided every $\sf X$-evaluation that makes all sentences of $\Gamma$ $\mathsf{T}$-true also makes $A$ $\mathsf{T}$-true.}

\item[-] {$\Gamma \models_{{\sf X/(SS\cap TT)}} A$, provided every $\sf X$-evaluation that makes all sentences of $\Gamma$ $\mathsf{S}$-true also makes $A$ $\mathsf{S}$-true, and every $\sf X$-evaluation that makes all sentences of $\Gamma$ $\mathsf{T}$-true also makes $A$ $\mathsf{T}$-true.}

\item[-] {$\Gamma \models_{{\sf X/ST}} A$, provided every $\sf X$-evaluation that makes all sentences of $\Gamma$ $\mathsf{S}$-true also makes $A$ $\mathsf{T}$-true.}

\item[-] {$\Gamma \models_{{\sf X/TS}} A$, provided every $\sf X$-evaluation that makes all sentences of $\Gamma$ $\mathsf{T}$-true also makes $A$ $\mathsf{S}$-true.}
\end{itemize}
\end{defn}
Relative to $\mathcal{L}$ and to {\sf SK}-evaluations, ${\sf SS}$-validity determines the so-called Strong Kleene logic, whereas ${\sf TT}$-validity determines Priest's logic \textsf{LP}. ${\sf SS\cap TT}$ corresponds to the so-called Symmetric Kleene logic, whereas ${\sf TS}$ and ${\sf ST}$ correspond to the so-called Tolerant-Strict and Strict-Tolerant Logics \citep[also called the logics of $q$-consequence and $p$-consequence:][]{malinowski1990q,Cobreros2012tolerant,frankowski2004formalization}. In general, our definitions of validity are relative to the choice of a type of evaluation function (e.g., classical, ${\sf SK}$, ${\sf DF}$); however, in the rest of this section, we always refer to ${\sf DF}$-evaluations, in line with our focus on the de Finetti conditional.


An interesting feature of the {\sf DF/TT}-logic is that it implies mutual entailment between its conditional and the material conditional. The following inferences are {\sf DF/TT}-valid:
\[\neg A \vee B \models_{\sf DF/TT} A \to B \hspace{2cm} A \to B \models_{\sf DF/TT} \neg A \vee B\]
Moreover, we also have: 
\[\models_{\sf DF/TT} (\neg A \vee B) \leftrightarrow (A \to B)\]
where $\leftrightarrow$ is de Finetti's biconditional, that is, $A \leftrightarrow B$ is defined as $(A \to B) \wedge (B \to A)$. {In fact, all of Gibbard's conditions (i) to (iii) are met by de Finetti's conditional in {\sf DF/TT}, and so the mutual entailment between de Finetti's conditional and its material counterpart may be seen as an instance of Gibbard's collapse result.}

Remarkably, however, although $\supset$ and $\to$ are equivalent in {\sf DF/TT}-logic, they don't obey the same principles. For instance:

\[A \to B \models_{\sf DF/TT} \neg (A \to \neg B) \hspace{.5cm} \mbox{ but } \hspace{.5cm} \neg A \vee B \not\models_{\sf DF/TT} \neg (\neg A \vee \neg B).\]


\subsection{A trilemma for de Finetti's conditional}
\label{S:3.2}

Among the previous schemes, which one is the most adequate relative to de Finetti's conditional? We begin with applying the {\sf SS}-validity scheme over {\sf DF}-evaluations, and similarly, mutatis mutandis, for the other schemes. It is easy to see that:
\[
A \to B \models_{\sf DF/SS} A \wedge B
\]
That is, the conditional entails conjunction. This property is not intuitive, but perhaps less bad than it seems since the trivalent approach is based on de Finetti's idea of identifying the truth conditions for conditionals with the conditions for winning a conditional bet. Worse is that the de Finetti conditional entails its converse on a {\sf SS}-validity scheme:\footnote{\citet[152]{reichenbach1944philosophic} claims the contrary, but because he seems to focus on the fact that $A\to B$ and $B\to A$ have different tables. He does not appear to see that they take the value 1 exactly in the same place, despite electing {\sf SS}-validity as his default notion of validity, in particular to guarantee Modus Ponens.}

\[
A \to B \models_{\sf DF/SS} B \to A
\]

The {$\sf SS$}-scheme is thus very distant from an intuitive notion of reasonable inference with conditionals since supposing $A$ and asserting $B$ is very different from supposing $B$ and asserting $A$. The $\sf TT$-scheme avoids this problem since  
\begin{align*}
A \to B &\not\models_{\sf DF/TT} A \wedge B & A \to B &\not\models_{\sf DF/TT} B \to A
\end{align*}
\citet{mcdermott1996truth} therefore proposes the {$\mathsf{SS\cap TT}$}-scheme to preserve the idea that validity is preservation of the value 1, but to weed out the implication from a conditional to the conjunction and to its converse. The {$\mathsf{SS\cap TT}$} consequence relation suffers, however, from the drawbacks of both of its constituents, as evidenced by the following observations: 
\begin{align*}
&\not\models_{\sf DF/SS} A \to A & A, A\to B &\models_{\sf DF/SS} B \\
&\models_{\sf DF/TT} A \to A & A, A\to B &\not\models_{\sf DF/TT} B
\end{align*}
{$\mathsf{DF}/{\sf(SS\cap TT)}$}  fails both the Identity Law ($A \to A$) and Modus Ponens: the first because {\sf DF/SS} has no sentential validities (as is the case in the Strong Kleene logic {\sf SK/SS}), the second because Modus Ponens is not valid in {\sf DF}/{\sf TT} (as is the case for the material conditional in Priest's {\sf LP = SK/TT}). As a result, the logic {$\mathsf{DF}/{\sf(SS\cap TT)}$} ends up being very weak. 

Consider now the so-called ``mixed consequence'' schemes, namely {\sf TS} and {\sf ST}, in which the level of truth varies from premises to conclusion \citep{Cobreros2012tolerant}. 
{\sf DF/TS} squares well with the degrees of assertability defined in Section \ref{S:2} since ${\it Ast}(A) \le {\it Ast}(B)$ for all underlying probability functions if and only if either $A$ and $B$ are logically equivalent, or $A \models_{\sf TS} B$ \citep[][166]{cantwell2006laws}. Hence, the logic connects well to epistemology, and it also eschews the  conjunction- and converse-conditional fallacies. Unfortunately, Modus Ponens and the Identity Law fail (like other sentential validities), not to mention other oddities of the logic, in which $A\not \models_{\sf DF/TS}A$. In {\sf DF/ST}, on the other hand, Modus Ponens and the Identity Law are retained, but also the entailment of the conditional to conjunction and to its converse remain.

We may summarize these observations in the form of a trilemma:

\begin{fact}\label{prop:trilemma}
Irrespective of whether {$\sf SS, TT, ST, TS, SS\cap TT$} is chosen for validity, a logic on $(\Lc, f_{\to_{\sf DF}})$ must either (1) fail Modus Ponens; or  (2) fail the Identity Law (and other sentential validities); or (3) validate the inference from a conditional to its converse. 
\end{fact}

\noindent The trilemma at a glance:
\[
\begin{tabular}{c|ccc}
 {\sf DF}/$\cdot$& MP & Identity & $\to\ \models\ \leftarrow$ \\
 
 \hline
 
 $\mathsf{SS}$ & $\checkmark$ & $\times$ & $\checkmark$\\
 
 $\mathsf{TT}$ & $\times$ & $\checkmark$ & $\times$ \\
 
 $\mathsf{ST}$ & $\checkmark$ & $\checkmark$ & $\checkmark$\\
 
 $\mathsf{TS}$ & $\times$ & $\times$ & $\times$ \\
 
  $\mathsf{SS\cap TT}$ & $\times$ & $\times$ & $\times$ \\
  \hline
  Ideal case & $\checkmark$ & $\checkmark$ & $\times$ \\

\end{tabular}
\]\medskip

The interest of this trilemma is that it involves schemata that depend on no other connective than the conditional. In what follows, we explore two main ways out of the trilemma: both select {\sf TT} validity as comparatively the best choice for validity, but the second moreover involves a modification of the de Finetti table so as to restore Modus Ponens.

\section{Giving up Modus Ponens: {\sf DF/TT} }\label{sec:DF+MP}

{Given that no validity scheme satisfies the three desiderata of making the {\sf DF} conditional validate Modus Ponens, avoid the entailment to its converse, and validity the Identity Law, one way out of the trilemma is to follow \citealt{quine1970philosophy}'s maxim of ``minimum mutilation'', and to elect as optimal the scheme or schemes that violate the fewer of those constraints}.\footnote{As in Optimality Theory (see \citealt{prince2008optimality}), we also assume that constraints can be rank-ordered in terms of how their comparative importance. We don't state the ordering explicitly here, the discussion makes it clear enough.}

{Three of the schemes violate two constraints, but $\mathsf{DF/TT}$ and $\mathsf{DF/ST}$ violate only one. However, $\mathsf{DF/ST}$ badly overgenerates (by validating the entailment to the converse), whereas $\mathsf{DF/TT}$ mildly undergenerates (by failing Modus Ponens, but still satisfying Conditional Introduction, see below). Arguably therefore, $\mathsf{DF/TT}$ appears to be the less inadequate of all options: it retains the Identity Law and avoids the entailment to the converse conditional, only at the expense of losing Modus Ponens---a principle that is given up in other logics such as Priest's LP (i.e., $\mathsf{SK/TT}$) for the material conditional.}\footnote{Note that unlike McGee's logic \citep{mcgee1989conditional}, which fails Modus Ponens for complex conditionals, $\mathsf{DF/TT}$ can fail Modus Ponens for simple conditionals, composed of atomic sentences.}

{Two more facts are worth highlighting about $\mathsf{DF/TT}$. Firstly, despite the failure of Modus Ponens, the conditional supports Conditional Introduction, namely $\Gamma, A\models B$ implies $\Gamma \models A \to B$. In $\mathsf{DF/SS}$, the situation is reversed, since Conditional Introduction fails despite Modus Ponens holding. Secondly, $\mathsf{DF/TT}$ supports full commutation of the conditional with negation, a schema widely regarded as plausible in natural language (see \citealt{cooper1968propositional, cantwell2008logic}, and Section 4.1 below).

\begin{fact}
For every $\{\Gamma, A, B\} \subseteq {\sf For}$:
\begin{description}
  \item[Conditional Introduction] If $\Gamma, A \models_{\mathsf{DF/TT}} B$, then $\Gamma \models_{\mathsf{DF/TT}} A \to B$.
  \item[Commutation with Negation] $\neg (A \to B) \equiv_{\mathsf{DF/TT}} A \to \neg B$. 
\end{description}
\end{fact}
\begin{proof}$\,$
\begin{itemize} 
\item[-] Suppose $\Gamma \not\models_{\mathsf{DF/TT}} A \to B$. Then there exists a $\mathsf{DF}$-evaluation $v$ such that for all $C \in \Gamma$, $v(C)>0$, but $v(A \to B)=0$. Hence $v(A)=1$, and $v(B)=0$, and $\Gamma, A \not\models_{\mathsf{DF/TT}} B$.
\item[-] Consider any $\mathsf{DF}$-evaluation $v$ such that $v(A\to \neg B)=0$. Then $v(A)=1$, $v(\neg B)=0$, so $v(B)=1$, and $v(A\to B)=1$, hence $v(\neg (A\to B))=0$, and the converse entailments hold.
\end{itemize}
\end{proof}
}
Despite blocking the entailment to the converse conditional, $\mathsf{DF/TT}$ validates several sentential schemata that are intuitively controversial. \citet{farrell1979material} for example points out that it validates the problematic schema $(B \wedge (A\to B)) \to A$, a sentential version of the fallacy of affirming the antecedent. More generally, we have:

\begin{fact}\label{fact:supp}
For every $A, B \in {\sf For}$:
\begin{align*}
&\models_{\mathsf{DF/TT}} (A \to B) \to A
\end{align*}
\end{fact}

\begin{proof}
For the principle to fail, there must be a $\mathsf{DF}$-evaluation $v$ such that $v(A\to B)=1$ and $v(A)=0$. But then $v(A\to B)=\half$, contradiction.
\end{proof}

Given the conditions the de Finetti conditional puts on {\sf TT}-validity, however, this schema does not necessarily constitute an unwelcome prediction. Firstly, it does not hold in argument form (that is, $A\to B\not \models_{\mathsf{TT}} A$), consistently with the fact that $\mathsf{TT}$-validity does not satisfy Modus Ponens. Secondly, consider the left-nested conditional sentence:

\ex. If Peter visits if Mary visits, then Mary will visit [indeed].

This seems intuitively acceptable, in line with the suppositional reading of the conditional.  

The upshot is that $\mathsf{DF/TT}$ loses some classical inferences based on the conditional (like Modus Ponens), and introduces some conditional sentences as validities that are not classical (viz. Fact \ref{fact:supp}), though not necessarily problematic under a suppositional reading. 

If, on the other hand, we wish to retain Modus Ponens as a central property of the conditional along with the Identity Law, then the trilemma presented in Fact \ref{prop:trilemma} implies that \emph{either some further notion of validity must be sought for the de Finetti conditional, or the de Finetti conditional itself is not adequate}. However, we have already argued that the notions of validity considered in this section exhaust the most natural and well-motivated class of monotonic notions of consequence defined over trivalent evaluations. For this reason, in the next section we explore that second option and explore alternatives to the de Finetti conditional.

\section{Retaining Modus Ponens: {\sf CC/TT}}
\label{S:5}


In this section, we show that under a {\sf TT}-definition of validity, de Finetti's table can be modified, and his motivations preserved, so as to preserve Modus Ponens and to avoid the previous trilemma. We first isolate the class of what we call Jeffrey conditionals. Within that class, we discuss some reasons to favor the Cooper-Cantwell conditional.

\subsection{Jeffrey conditionals}
\label{S:5.1}

In a short and underappreciated note, \citet{jeffrey1963indeterminate} highlighted the following condition for a trivalent operator to satisfy Modus Ponens when {$\sf TT$} is used for validity:

\begin{fact}
Under a $\mathsf{TT}$-notion of validity, a trivalent conditional operator $f_{\to}$ validates Modus Ponens only if $f_{\to}(1,0)=f_{\to}(\half,0)=0$.
\end{fact}

\begin{proof}
Assume $f_{\to}(1,0)\neq 0$ or $f_{\to}(\half, 0)\neq 0$. Then it is possible to have $v(A)> 0, v(A\to B)> 0$ and $v(B)=0$, which invalidates Modus Ponens.
\end{proof}


\noindent We may therefore call a conditional operator \emph{Jeffrey} if it extends the bivalent ``gappy'' conditional as follows (\citealt{jeffrey1963indeterminate}):

 \begin{defn}\label{def:jeff}

 A \emph{Jeffrey conditional} is any binary trivalent operator of the form:
 \[
 \begin{tabular}[t]{c|ccc}
 $f_{\to}$& $1$ & $\nicefrac{1}{2}$ & $0$\\
  \hline
  $1$ & $1$ &  $d_1$ & $0$\\
 $\nicefrac{1}{2}$ & $d_2$& $d_3$& $0$\\
  $0$ & $\half$ &  $d_4$& $\half$\\
  \end{tabular}\]
 \noindent where $d_i \in \{\half, 1\}$ for $1 \le i \le 4$.\footnote{Jeffrey contends that any completion of the gappy truth table must satisfy this schema; to prove this claim he demands that any acceptable logic satisfy Modus Ponens, Syllogism, the Deduction Theorem and Contraposition. His argument depends on choosing a negation operator mapping designated values under the $\mathsf{TT}$-scheme to a nondesignated value, and conversely \citep[for more details on this connection, see][]{chemla2018mvcons}.} 
 \end{defn}

An operator can therefore satisfy Jeffrey's constraint and be de Finettian at the same time, namely comply with the truth conditions of de Finetti's conditional when the antecedent has a classical truth value (see Definition \ref{D:definettian}). 
We thus say that:
 
 \begin{fact}
 A Jeffrey conditional is \emph{de Finettian} provided it is of the form:
  \[
 \begin{tabular}[t]{c|ccc}
 $f_{\to}$& $1$ & $\nicefrac{1}{2}$ & $0$\\
  \hline
  $1$ & $1$ &  $\half$ & $0$\\
 $\nicefrac{1}{2}$ & $d_2$& $d_3$& $0$\\
  $0$ & $\half$ &  $\half$& $\half$\\
  \end{tabular}
  \]
 \noindent where $d_2, d_3 \in \{\half, 1\}$.
 \end{fact}
 

\begin{figure}[h]
\[
\begin{tabular}[t]{l|ccc}
$f_{\to_{\sf CC}}$ & 1 & $\nicefrac{1}{2}$ & 0\\
 \hline
 1 & 1 &  $\nicefrac{1}{2}$ & 0\\
$\nicefrac{1}{2}$ & 1& $\nicefrac{1}{2}$& 0\\
 0 & $\nicefrac{1}{2}$ &  $\nicefrac{1}{2}$& $\nicefrac{1}{2}$\\
 \end{tabular}
 \qquad
 \begin{tabular}[t]{l|ccc}
$f_{\to_{\sf F}}$ & 1 & $\nicefrac{1}{2}$ & 0\\
 \hline
 1 & 1 &  $\nicefrac{1}{2}$ & 0\\
$\nicefrac{1}{2}$ & $\nicefrac{1}{2}$& $\nicefrac{1}{2}$& 0\\
 0 & $\nicefrac{1}{2}$ &  $\nicefrac{1}{2}$& $\nicefrac{1}{2}$\\
 \end{tabular}
 \]
 \[
 \begin{tabular}[t]{l|ccc}
$f_{\to_{\sf J2}}$ & 1 & $\nicefrac{1}{2}$ & 0\\
 \hline
 1 & 1 &  $\nicefrac{1}{2}$ & 0\\
$\nicefrac{1}{2}$ & 1& 1 & 0\\
 0 & $\nicefrac{1}{2}$ &  $\nicefrac{1}{2}$& $\nicefrac{1}{2}$\\
 \end{tabular}\qquad
 \begin{tabular}[t]{l|ccc}
$f_{\to_{\sf J1}}$ & 1 & $\nicefrac{1}{2}$ & 0\\
 \hline
 1 & 1 &  $\nicefrac{1}{2}$ & 0\\
$\nicefrac{1}{2}$ & $\nicefrac{1}{2}$& 1& 0\\
 0 & $\nicefrac{1}{2}$ &  $\nicefrac{1}{2}$& $\nicefrac{1}{2}$\\
 \end{tabular}\]

 \caption{The de Finettian Jeffrey conditionals}\label{fig:jeffrey}
 \end{figure}

\noindent Clearly, there exist four de Finettian Jeffrey conditionals (see Figure \ref{fig:jeffrey}). Two of them are the Cooper-Cantwell ({\sf CC}) and the Farrell conditional ({\sf F}). We call the other two {\sf J1} and {\sf J2}. For each such table, we modify the notion of {\sf DF}-evaluation accordingly (call it a {\sf CC}-, {\sf F}-, {\sf J1}-, and {\sf J2}-evaluation respectively).

It is straightforward to see that Jeffrey conditionals (whether de Finettian or not) eschew the trilemma faced by de Finetti's:
 
 \begin{prop}[Trilemma Resolution]
Under a $\mathsf{TT}$-notion of validity, any Jeffrey conditional 
\begin{itemize}
  \item[-] satisfies Modus Ponens and the Identity Law;
  \item[-] invalidates the entailment of the conditional to its converse.
 \end{itemize}
\end{prop}

\begin{proof}$\,$
\begin{itemize}
\item[-] Modus Ponens: Assume $v(A)>0$ and $v(A\to B)> 0$, then clearly $v(B)>0$.
\item[-] Identity: All values on the diagonal of any Jeffrey conditional differ from $0$.
\item[-] Avoiding the entailment to the converse: When $v(A)=0$ and $v(B)=\half$, $v(A\to B)>0$ but $v(B\to A)=0$: this invalidates the entailment from $A\to B$ to $B\to A$. 
\end{itemize}
\end{proof}

Like de Finetti's conditional, all Jeffrey conditionals {\sf TT}-validate Conditional Introduction, but unlike the de Finetti conditional they satisfy the converse, namely the full Deduction Theorem. In fact, there is a precise sense in which {\sf TT}-validity and Jeffrey conditionals fit each other:\footnote{Compare with \citet{chemla2018mvcons}, who examine which conditionals of a specific form are admitted by a given consequence relation in trivalent and higher-valued logic. Here we partly reverse this problem, by looking at which validity scheme, if any, is most appropriate to a given conditional operator.}

\begin{prop}[Deduction Theorem]
\item Any Jeffrey conditional $\mathsf{TT}$-validates both directions of the Deduction Theorem, that is for every ${\sf J}$-evaluation,
\begin{equation}
\Gamma, A \models_{\mathsf{J/TT}} B \qquad \text{if and only if} \qquad \Gamma\models_{\mathsf{J/TT}} A \to B \tag{Deduction Theorem}\label{eqn:DT}
\end{equation}    
\item No Jeffrey conditional validates the full Deduction Theorem for $\mathsf{SS}$-, $\mathsf{TT \cap SS}$, $\mathsf{ST}$ and $\mathsf{TS}$-validity. 
\end{prop}
 
\begin{proof}$\,$
\begin{description}
\item[Deduction Theorem for ${\sf TT}$-validity:] $\,$\\
($\Rightarrow$) Suppose $\Gamma \not \models_{\mathsf{J/TT}} A\to B$. Then there is some ${\sf J}$-evaluation $v$ such that for all $C \in \Gamma$,  $v(C)>0$ but $v(A)>0$ and $v(B)=0$. This implies that $\Gamma, A\not\models_{\mathsf{J/TT}} B$. \\
($\Leftarrow$) Suppose $\Gamma, A\not \models_{\mathsf{J/TT}} B$: there there is some ${\sf J}$-evaluation $v$ such that for all $C \in \Gamma$,  $v(C)>0$, $v(A)>0$, but $v(B)=0$. Hence, $v(A\to B)=0$, and $\Gamma \not \models_{\mathsf{J/TT}} A\to B$.

\item[Failure of Deduction Theorem for ${\sf SS}$-, ${\sf ST}$-, ${\sf SS} \cap {\sf TT}$-, and ${\sf TS}$-validity:] $\,$\\
($\Leftarrow$) For $\mathsf{SS}$-validity and failure of the Deduction Theorem, consider $v(A) = \nicefrac{1}{2}$ and $v(B) = 0$. Then the Jeffrey conditional $A \to B$ is false, but the entailment $A \models_{\sf SS} B$ holds. The same case shows failure of the Deduction Theorem for $\mathsf{ST}$-validity. For $\mathsf{SS \cap TT}$- and $\mathsf{TS}$-validity, consider the case $v(A) = 0$ and $v(B) = 1$.
\end{description}
\end{proof}
 
This result is important since our consequence relation is meant to capture a suitable \emph{logic of suppositional reasoning}, in line with de Finetti's original motivation. Just as the truth table for the trivalent conditional is motivated by the idea of evaluating the consequent under the supposition of the antecedent, the consequence relation should describe the inferences that are licensed by supposing the antecedent. Therefore, a deduction theorem is an important adequacy condition for a logic of trivalent conditionals, making a strong case for $\mathsf{TT}$-validity in combination with Jeffrey conditionals. Relatedly, it can be seen that no Jeffrey conditional supports $(A\to B) \to A$ as a valid schema relative to {\sf TT}-validity (to see this, let $v(A)=0, v(B)=\half$), unlike de Finetti's conditional (see Fact \ref{fact:supp} and compare \citealt{farrell1979material}, whose motivation for $\to_{\sf F}$ lies precisely here). {Finally, we saw that relative to {\sf TT}-validity de Finetti's conditional is logically equivalent to the material conditional. By contrast, every Jeffrey conditional relative to that same scheme is strictly stronger than the material conditional. Relative to {\sf TT}-validity, Jeffrey conditionals do not fall prey to Gibbard's collapse result, basically because they do not support Gibbard's condition (ii): when $A$ classically entails $C$, $A\to C$ need not be valid.}

\subsection{Negation and {\sf CC/TT}}
\label{S:5.2}

To choose between the various Jeffrey conditionals, we suggest to look at the interplay of the conditional with the other logical connectives. The interplay between conditional and negation is especially relevant, since several of the most debated principles involving indicative conditionals concern negation as well. One common fact about Jeffrey conditionals is that they fail contraposition relative to Strong Kleene negation:

\begin{prop} 
For any Jeffrey conditional, $A \to B \not \models_{\mathsf{J/TT}} \neg B \to \neg A$.
\label{P:contra}
\end{prop}
\begin{proof}
Suppose $v(A)=1, v(B)=\half$. Then $v(A\to B)=\half$, but $v(\neg B \to \neg A)=0$. Hence, $A\to B\not \models_{\mathsf{J/TT}} \neg B\to \neg A$.
\end{proof}

\noindent The failure of Contraposition may be seen as a welcome prediction. First of all, supposing $A$ and supposing $\neg B$ are just two different things. For example, when $v(A)=v(B)=1$, then $A \to B$ is obviously true, whereas $\neg B \to \neg A$ is now ``void''---the conditions for evaluating its truth or falsity are not satisfied. Therefore $v(\neg B \to \neg A) = \half$. Second, contraposition does not always preserve meaning. The contrapositive of a sentence like ``if Sappho did not die in 570 BC, then she is dead by now'' would be ``if Sappho is not dead by now, then she died in 570 BC''. The latter obviously conveys a different thought. Hence the inference to the contrapositive is not warranted in all situations.\footnote{Of course, we are assuming double negation elimination inside conditionals---this seems entirely unproblematic. Accounts where Contraposition holds, such as the refined material conditional view of \citet{Jackson1979,Jackson1987}, have to go to some length to explain away the counterintuitive feel of such examples.} 
Since all Jeffrey conditionals satisfy the deduction theorem relative to ${\sf TT}$-validity, this means they also fail to validate Modus Tollens. Modus Tollens is not ${\sf DF/TT}$-valid either, though Contraposition is.

On the other hand, as noted by \citet{cooper1968propositional} and \citet{cantwell2008logic}, the Cooper-Cantwell conditional supports the full commutation of Strong Kleene negation with the conditional, namely the logical equivalence between $\neg (A\to B)$ and $(A\to \neg B)$.\footnote{Prof.~Farrell (p.c.) draws our attention to the fact that his table supports full commutation for conditionals involving atomic sentences, when restricted to atom-classical valuations. Because the restriction to atom-classical valuations is defended by Cooper, the two accounts mostly differ on nested conditionals.} In fact, it is the only Jeffrey conditional that does so:

\begin{prop}\label{fact:commut}

Among all Jeffrey conditionals, only the Cooper-Cantwell conditional validates the full commutation schema for negation. For de Finettian Jeffrey conditionals, in particular, $\mathsf{SK}$-negation is a separating connective:

\[
\begin{tabular}[t]{l|c|c}

${\sf J}=\cdot$ & $\neg (A \to B)\models_{\mathsf{J/TT}} A\to \neg B$ & $A\to \neg B\models_{\mathsf{J/TT}}  \neg (A \to B)$\\
 \hline
$\mathsf{CC}$ & $\checkmark$ & $\checkmark$\\

$\mathsf{F}$ & $\mathsf{\times}$ & $\checkmark$\\

$\mathsf{J1}$ & $\mathsf{\times}$ & $\times$\\

$\mathsf{J2}$ & $\checkmark$ & $\times$\\

\end{tabular}\]

\end{prop}

\begin{proof}
From the definition of a Jeffrey conditional in Definition \ref{def:jeff}, the truth tables for $\neg(A \to B)$ and $A \to \neg B$ look like this:
\[
\begin{tabular}[t]{c|ccc}
$\neg (A \to B)$ & $1$ & $\nicefrac{1}{2}$ & $0$\\
\hline
  $1$ & $0$ & $\neg d_1$ & $1$\\
 $\nicefrac{1}{2}$ & $\neg d_2$ & $\neg d_3$ & $1$\\
  $0$ & $\half$ &  $\neg d_4$ & $\half$
 \end{tabular} \qquad
\begin{tabular}[t]{c|ccc}
$A \to \neg B$ & $1$ & $\nicefrac{1}{2}$ & $0$\\
 \hline
  $1$ & $0$ & $d_1$ & $1$\\
 $\nicefrac{1}{2}$ & $0$ & $d_3$& $d_2$\\
  $0$ & $\half$ &  $d_4$ & $\half$\\
 \end{tabular}
\]

For $\mathsf{TT}$-entailment to go in both directions, necessarily, $\neg d_{2}=0$, hence $d_{2}=1$, and $d_{1}, d_{3}, d_4$ must all equal $\half$, which yields the table for the Cooper-Cantwell conditional.

For the other de Finettian Jeffrey cases: let $v$ be an $\mathsf{F}$-evaluation, or a $\mathsf{J1}$-evaluation: assume $v(A)=\half$ and $v(B)=1$, then $v(\neg (A\to B))=\half)$, but $v(A \to \neg B)=0$. 
Let $v$ be a $\mathsf{J1}$-evaluation, or a $\mathsf{J2}$-evaluation: assume $v(A)=\half$ and $v(B)=\half$, then $v(A\to \neg B)=1$, but $v(\neg (A \to B))=0$.
Consider any $\mathsf{J2}$-evaluation. To show that $\neg (A\to B) \models A \to \neg B$, assume that there is a $v$ such that $v(A \to \neg B)=0$, but $v(\neg (A \to B))>0$. Necessarily, $v(A)>0$, but $v(\neg B)=0$, so $v(B)=1$. But then $v(A\to B)=1$, and $v(\neg (A\to B))=0$, contradiction.
\end{proof}

In classical logic, only the commutation from outer to inner negation is valid. On the other hand, inferences in natural language appear to support both directions in many contexts.  \citet{Ramsey1929}, \citet{Adams1965}, \citet{cooper1968propositional}, \citet{cantwell2008logic} and \citet{francez2016natural} give a theoretically motivated defense of the commutation scheme, while the studies by \citet{HandleyEvansThompson2006} and \cite{politzer2009could} provide some empirical support. See, however, \citet{egre2013negation}, \citet{olivier2018enonces} and \citet{skovgaard2019cancellation} for a more complex picture.

\subsection{Connexivity}

We conclude this section by briefly relating our discussion of the {\sf TT}-logics of de Finettian and Jeffrey conditionals to a slightly wider logical context. A conditional logic is called \emph{connexive} if it validates the two following schemata:
\begin{equation}\label{eqn:AT}
\neg (\neg A \to A) \tag{Aristotle's Thesis} \\
\end{equation}
and
\begin{equation}\label{eqn:BT}
(A \to C) \to \neg (A \to \neg C) \tag{Boethius' Thesis} 
\end{equation}
(see \citealt{pizzi1977boethius, sep-logic-connexive}). Both de Finetti's conditional and the Cooper-Cantwell conditional are connexive when paired with {\sf TT}-validity (and Strong Kleene negation).\footnote{\citealt{kapsner2018humble} is a recent paper that restates de Finetti's table specifically in relation to connexivity.} Neither Aristotle's Thesis nor Boethius' Thesis are classical tautologies: indeed, connexive logics are not subsystems of classical logic.\footnote{Classical counterexamples are easily obtained by assigning value $1$ to $A$ in Aristotle's Thesis, and value $0$ in Boethius' Thesis.} On the other hand, systems of connexive logic lack some classical principle, lest they are trivial (of course, {\sf DF/TT} and {\sf CC/TT} are no exception). Informally construed, Aristotle's Thesis requires that it is never the case that a formula is implied by its own negation, while Boethius' Thesis requires that if a conditional $A \rightarrow C$ holds, then it is not the case that the conditional that results from the former by negating the consequent, i.e. $A \rightarrow \neg C$ (which is equivalent to the negated conditional $\neg (A \rightarrow C)$ in both {\sf DF/TT} and {\sf CC/TT}) hold. Now, since both {\sf DF/TT} and {\sf CC/TT} employ a tolerant-tolerant notion of validity, the fact that they satisfy  Boethius' Thesis can hardly be interpreted as saying that they show that a conditional is `incompatible' with its negation (and similarly for Aristotle's Thesis). Nevertheless, in requiring such a strict, extra-classical connection between antecedent and consequent of a conditional, connexive logics---including {\sf DF/TT} and {\sf CC/TT}---arguably ensure that the conditional interacts reasonably well with negation.

Nevertheless, the interaction of conditional and negation displayed by connexive logics of De Finettian and Jeffrey conditionals, {\sf DF/TT} and {\sf CC/TT} in particular, is not entirely free from worries. For one thing, connexivity comes at a price when it comes to reductio proofs (see \citealt{cooper1968propositional} for discussion). For another, like de Finetti's conditional, the Cooper-Cantwell conditional also validates the following equivalence, where $\equiv_{\it m}$ is the material biconditional: 
\[\neg (A \to B) \equiv_{\it m} (A\to\neg B)\]
As a consequence, both conditionals validate 
\begin{description}
  \item[Conditional Excluded Middle] $\models_{\mathsf{\cdot/TT}} \hspace{5pt} (A \to B) \vee (A\to \neg B)$
\end{description}
Conditional Excluded Middle is a moot principle, but it is a natural one to have if negation is to commute with the conditional.\footnote{For a criticism of Conditional Excluded Middle, see \citet{Lewis1973}, for a defense see \citet{Stalnaker1980defense}.} Moreover, since every de Finettian Jeffrey conditional validates Conditional Excluded Middle, this does not tell against the Cooper-Cantwell variant. Thanks to the fact that it is the only one, within the de Finettian Jeffrey conditionals, to support the full commutation with negation, the Cooper-Cantwell conditional stands out as the closest to de Finetti's original connective.

\section{Comparisons and Limits}
\label{S:6}

We have distinguished two trivalents logics of indicative conditionals, namely {\sf DF/TT} and {\sf CC/TT}, whose proof theory and algebraic semantics we will explore in Part II of this paper. Before doing so, let us summarize the commonalities between the two logics, their principal differences, and draw comparisons with other logics of conditionals.

Four main features are common to \dftt\ and \cctt: they are truth-functional logics, they share the same de Finettian semantic core, they are connexive, and both support the law of Import-Export without restriction. The main difference between \dftt\ and \cctt\ is that the former fails Modus Ponens, whereas the latter preserves it, so that only \cctt\ supports the full Deduction Theorem. This property is in line with the fact that for {\sf TT}-validity, the designated values are 1 and $\half$, and the Cooper-Cantwell conditional is only evaluated as false when the antecedent is designated and the consequent undesignated. Conversely, relative to Strong Kleene negation the Cooper-Cantwell conditional fails Contraposition, whereas de Finetti's conditional supports Contraposition, but both fails Modus Tollens.

The preservation of Modus Ponens may be seen as virtue of \cctt\ compared to \dftt. However, one common fact about both logics, given our assumption that they share the same Strong Kleene disjunction, is that they fail the rule of Disjunctive Syllogism ($\neg A, A\vee B\models B$). Clearly, this concerns the table for disjunction for a {\sf TT}-consequence relation (see \citealt{priest1979logic, cantwell2008logic}), independently of the particular truth conditions for the conditional.

Because the Law of Import-Export is validated, in both \cctt\ and \dftt\ only one of the paradoxes of material implication is blocked, namely the schema $A \to (\neg A \to B)$. On the other hand, $A \to (B \to A)$ holds in both logics, consistent with the fact that $A\wedge B \to A$ is valid. As discussed in Section \ref{sec:DF+MP}, this property squares well with the proposed suppositional interpretation of the conditional. Given the way conjunction and disjunction are handled in \dftt\ and \cctt, we can therefore conclude that whereas both logics are connexive, neither is relevantist, except in a weak sense (by failing one of the paradoxes of material implication). 

We now discuss some limitations of our logics. First, both \cctt\ and \dftt\ validate the so-called Linearity principle $(A\to B) \vee (B\to A)$. This schema was famously criticized by \citet{maccoll1908if}, who pointed out that neither of ``if John is red-haired, then John is a doctor'' and ``if John is a doctor, then he is red-haired'' seems acceptable in ordinary reasoning. 


Second, there is a certain tension between our extensional semantics of conditionals and the intensional use to which they are often put. Suppose Mary believes the following conditional: 

\ex. If the Church is East of the City Hall, then the City Hall is West of the Church \label{ex:church}

Intuitively the proposition that Mary believes appears analytically true. Nonetheless, on the de Finettian analysis its truth value depends on the position of the City Hall with respect to the Church: the conditional may be evaluated either as true or as indeterminate. The apparent analyticity of \ref{ex:church} has to be explained by reference to it being maximally assertable, regardless of its actual truth value. In fact, also \citet[315]{Lewis1986} observes that ``there is a discrepancy between truth- and assertability-preserving inference involving indicative conditionals; and that our intuitions about valid reasoning with conditionals are apt to concern the latter, and so to be poor evidence about the former.'' In other words, while \dftt\ and \cctt\ aim at describing a logic of suppositional reasoning and  their analysis of \ref{ex:church} should be evaluated by these criteria, reasonable inferences with conditionals, including ``apparent analytic truths'', may need to be analyzed in terms of a (probabilistic) theory of assertability. This theory can again be anchored in, and motivated by, trivalent truth conditions for conditionals---see Section \ref{sec:definetti}. Detailing the division of labor between semantics (truth conditions, validity) and epistemology (degrees of assertability) is, however, a project for future work.

Third, some conjunctive sentences can never be true on \dftt\ or \cctt, because one of the conjuncts will always be indeterminate. An ``obvious truth'' such as $(A \to A) \wedge (\neg A \to \neg A)$ is always classified as indeterminate (we are indebted to Paolo Santorio for this example). Likewise, a ``partitioning sentence'' of the form $(A \to B) \wedge (\neg A \to C)$ will always be indeterminate or false \citep[][368--370]{Bradley2002}. However, a sentence such as: 

\ex. If the sun shines tomorrow, John goes to the beach; and if it rains, he goes to the museum. \label{ex:beach}

seems to be true (with hindsight) if the sun shines tomorrow and John goes indeed to the beach. 


Both cases raise a challenge for an account of the assertability of compound conditionals. An attempt to deal with them may be to extend the assertability principle \ref{eq:ast} to sentences of arbitrary logical complexity, and to stipulate logical validities to have degree of assertability $1$. This would make  $(A \to A) \wedge (\neg A \to \neg A)$ maximally assertable. However, it would remain silent on the assertability of Bradley-type sentences. For the latter, an unrestricted version of \ref{eq:ast} would predict the degree of assertability to be $0$ (since Bradley's sentence must be false when its truth value is classical), obviously a bad prediction. Because of that, principle \ref{eq:ast} may indeed have to be restricted to non-compound conditionals, and a full account of assertability for compound sentences may have to be a recursive account, based on the assertability of simpler sentences.\footnote{For example, both for Bradley's example, and for Santorio's example, the assertability of each conjunction is better predicted if thought of as the minimum of the assertability of each conjunct as predicted by principle \ref{eq:ast}. We leave an elaboration of that idea for further work.} 




A possibly more elegant way of overcoming these limitations, proposed by \citet{cooper1968propositional}, is to introduce different truth tables for trivalent conjunction and disjunction---see Table \ref{tab:quasi} \citep[see also][1044--1053]{humberstone2011conn}. These truth tables, where the conjunction of the true and the indeterminate is the true (and vice versa for disjunction), can be motivated by the isomorphism between bets and truth values introduced in Section \ref{S:2}: a system of bets should be classified as winning if it consists of a winning and a called-off bet. Adopting this quasi-conjunction and quasi-disjunction, invalidates Linearity and resolves the problem with the truth and assertability conditions of partitioning sentences. In particular, $(A \to A) \wedge (\neg A \to \neg A)$ is always true, and so is $(A \to B) \wedge (\neg A \to C)$ when its first conjunct is true. However, when paired with \dftt, quasi-conjunction leads to a violation of Import-Export; so it should be considered only as a possible modification of \cctt.

\begin{figure}[h]
\[
\begin{tabular}{|c|ccc|}
\hline
$f'_{\wedge}$ & $1$ & $\nicefrac{1}{2}$ & $0$\\
\hline
$1$ & $1$ & $1$ & $0$\\
$\nicefrac{1}{2}$ & $1$ & $\nicefrac{1}{2}$ & $0$\\
$0$ & $0$ & $0$ & $0$\\
\hline
\end{tabular}
\hspace{10pt}
\begin{tabular}{|c|ccc|}
\hline
$f'_{\vee}$ & $1$ & $\nicefrac{1}{2}$ & $0$\\
\hline
$1$ & $1$ & $1$ & $1$\\
$\nicefrac{1}{2}$ & $1$ & $\nicefrac{1}{2}$ & $0$\\
$0$ & $1$ & $0$ & $0$\\
\hline
\end{tabular}
\]
\caption{Truth tables for trivalent quasi-conjunction and quasi-disjunction, as advocated by \citet{cooper1968propositional}.}\label{tab:quasi}
\end{figure}

On a general level, quasi-disjunction violates Disjunction Introduction ($A \models A \vee B$, and quasi-conjunction the dual inference from $\neg A$ to $\neg (A \wedge B)$), but this feature is in line with a relevantist solution to the paradoxes of material implication. More surprising is perhaps that the material conditional $\neg A \vee C$ is now logically \textit{stronger} than the indicative conditional $A \to C$---a feature that needs closer analysis. On the positive side, the two connectives in Table \ref{tab:quasi} are dual to each other and thus satisfy the de Morgan rules. Conjunction Elimination ($A \wedge B \models A$) still holds for \cctt, and so do most other desirable principles (e.g., Import-Export, Distributivity, Conditional Excluded Middle). The results of Section \ref{S:5}---the classification of Jeffrey conditionals, the Deduction Theorem, commutation with negation, the connexive principles---also stay intact since they do not depend on the choice of the connective for conjunction and disjunction. An additional benefit of operating with quasi-disjunction instead of Strong Kleene disjunction is the validation of the rule of Disjunctive Syllogism in \cctt\ ($\neg A, A \vee B \models B$). Ultimately, the choice between these two versions of \cctt\ is up to the reader, dependent on how they weigh the distinctive features of the two resulting logics. In any case, quasi-conjunction and -disjunction offer a principled way of responding to philosophically minded objections that have long plagued advocates of trivalent semantics for indicative conditionals.

\section{Summary and Perspectives}
\label{S:7}

De Finetti's trivalent conditional was put forward by de Finetti to qualitatively model the way in which conditional statements are probabilistically represented. Since its discovery, the {\sf DF} table has received a fair amount of attention from mathematicians as well as psychologists, but there have been surprisingly few investigations of the trivalent logics supported by the conditional as well as the variants in its vicinity. Our main motivation for this paper has been to fill this gap.

We started with the observation that de Finetti's truth table faces a trilemma when confronted with the choice of a trivalent validity relation: give up the Identity Law and other sentential validities, support the entailment from a conditional to its converse, or give up Modus Ponens. We have argued that the latter option is the less costly in relation to its alternatives, if the {\sf DF} conditional is paired with a notion of {\sf TT}-validity. On the other hand, trivalent Jeffrey conditionals, which have the property $f_{\to}(\half, 0) =  0$, avoid this trilemma when endowed with the same $\mathsf{TT}$-consequence relation: they block the entailment to the converse conditionals, they support the Identity Law, and moreover they support the full Deduction Theorem (Modus Ponens and Conditional Introduction), in line with the fact that the values $1$ and $\half$ are designated for consequence, and pattern in the same way for those conditionals.

Zooming in on Jeffrey conditionals, we see that the Cooper-Cantwell conditional stands out in that it satisfies the full commutation schema for negation, a schema widely regarded as plausible in natural language, also supported by the de Finetti conditional. Prima facie therefore, the Cooper-Cantwell conditional appears to strike the best balance between logical and epistemological properties: like Farrell's conditional, but unlike de Finetti's, it satisfies Modus Ponens. Its motivation for the middle line of its truth table---to treat an indeterminate antecedent like a true one---is more stringent than Farrell's, and well-aligned with the $\mathsf{TT}$-consequence relation.

As pointed out in the previous section, both $\mathsf{CC/TT}$ and $\mathsf{DF/TT}$ share features which may be seen as problematic, such as the Linearity principle and the treatment of partitioning sentences. 
A principled way out of these problems that merits further attention is to modify \cctt\ by changing the connectives for conjunction and disjunction along the lines of \citealt{cooper1968propositional}. 
From a methodological point of view, however, we think it matters to any further work on conditionals to locate exactly the (actual and alleged) limits of the trivalent approach, in particular because they should be carefully compared to some of the benefits we highlighted.  In Part II of this paper, we therefore propose a more elaborate treatment of the proof theory and algebraic semantics of both \cctt\ and \dftt, in order to give a more informed assessment of both logics.

\pagebreak
\bibliographystyle{../BibFolder/mychicago-ff}

\renewcommand{\bibname}{References}
\bibliography{../BibFolder/v-c(14)}

\end{document}


\maketitle
\thispagestyle{empty}
\vspace{-2.5cm}
\centerline{\large Part II: Proof Theory and Algebraic Semantics}
\vspace{1cm}

\centerline{\small Paul \'Egr\'e, Lorenzo Rossi, Jan Sprenger}

\begin{abstract}
\noindent In Part I of this paper, we identified and compared various schemes for trivalent truth conditions for indicative conditionals, most notably the proposals by \citet{definetti1936logique} and \cite{reichenbach1944philosophic} on the one hand, and by \citet{cooper1968propositional} and \citet{cantwell2008logic} on the other. Here we provide the proof theory for the resulting logics {\sf DF/TT} and {\sf CC/TT}, using tableau calculi and sequent calculi, and proving soundness and completeness results. Then we turn to the algebraic semantics, where both logics have substantive limitations: \dftt\ allows for algebraic completeness, but not for the construction of a canonical model, while \cctt\ fails the construction of a Lindenbaum-Tarski algebra. With these results in mind, we draw up the balance and sketch future research projects.  
\end{abstract}

\medskip


\setcounter{section}{0}

In Part I of this paper, we have reviewed the motivations for a trivalent semantic treatment of indicative conditionals, centered on the  proposal made by \citet{definetti1936logique}, and \citet{reichenbach1944philosophic}, to treat indicative conditionals as conditional assertions akin to conditional bets. We have singled out two de Finettian logics of the indicative conditional, the first based on de Finetti's table, paired with a notion of logical consequence as preservation of non-Falsity ({\sf TT}-validity), the other based on a close kin of de Finetti's table, the Cooper-Cantwell table, paired with the same notion of validity (\citealt{cooper1968propositional, cantwell2008logic}). 
These logics are called {\sf DF/TT} and {\sf CC/TT}, respectively. We repeat the truth tables of the conditional operator in Figure \ref{fig:DFCC} and the definition of {\sf TT-}validity below. 

\begin{figure}[htb]
\[
 \begin{tabular}[t]{l|ccc}
$f_{\to _{\sf DF}}$ & $1$ & $\nicefrac{1}{2}$ & $0$\\
 \hline
$1$ & $1$ &  $\nicefrac{1}{2}$ & $0$\\
$\nicefrac{1}{2}$ & $\nicefrac{1}{2}$& $\nicefrac{1}{2}$& $\nicefrac{1}{2}$\\
$0$ & $\nicefrac{1}{2}$ &  $\nicefrac{1}{2}$& $\nicefrac{1}{2}$\\
 \end{tabular}
 \qquad
\begin{tabular}[t]{l|ccc}
$f_{\to _{\sf CC}}$ & $1$ & $\nicefrac{1}{2}$ & $0$\\
 \hline
$1$ & $1$ &  $\nicefrac{1}{2}$ & $0$\\
$\nicefrac{1}{2}$ & 1& $\nicefrac{1}{2}$& $0$\\
$0$ & $\nicefrac{1}{2}$ &  $\nicefrac{1}{2}$& $\nicefrac{1}{2}$\\
 \end{tabular}
 \]\caption{truth tables for the de Finetti conditional (left) and the Cooper-Cantwell conditional (right).}\label{fig:DFCC}
 \end{figure}
 
\begin{description}
  \item[{\sf TT-}validity] $\Gamma \models_{\sf X/TT} A$ provided every {\sf X}-evaluation that makes all sentences of $\Gamma$ $\mathsf{T}$-true also makes $A$ $\mathsf{T}$-true. In other words, if for every {\sf X}-evaluation function $v: {\sf For} \longmapsto \{0, \half, 1\}$ (where {\sf X} specifies the interpretation of $\to$ and other connectives), for every sentence $B \in \Gamma$, $v(B) \in \{\half, 1\}$, then also $v(A) \in \{\half, 1\}$.
  \end{description}
  

As easily seen from the tables, both conditionals are \emph{de Finettian} in the following sense: they take the value of the consequent when the antecedent is true ($=1$), and the value indeterminate when the antecedent is false ($=0$). They differ when the antecedent itself is indeterminate ($=\half$): whereas the {\sf DF} conditional groups $\half$ with $0$ in antecedent position,  the {\sf CC} conditional groups $\half$ with $1$ instead. In both logics (whether {\sf X=DF} or {\sf CC}), the non-conditional connectives `$\neg$', `$\wedge$' and `$\vee$' are interpreted according to the \L ukasiewicz/de Finetti/Strong Kleene truth tables, where negation swaps 1 and 0, mapping $\half$ to itself, and where conjunction and disjunction are interpretable as $min$ and $max$ respectively. Although alternative tables are given in \citet{cooper1968propositional} for disjunction and conjunction, in this part as well as in the previous one, our focus remains on the standard interpretation of the connectives.

In Part I, it was pointed out that the resulting logics {\sf DF/TT} and {\sf CC/TT} share some distinctive features, in particular both satisfy Conditional Introduction and the law of Import-Export, and both are connexive logics (supporting $\neg (A \to \neg A)$ and the inference from $(A\to \neg B)$ to $\neg (A \to B)$). In fact, both support unrestricted commutation of the conditional with negation. This feature sets them apart from other de Finettian logics, in particular \citet{farrell1979material}'s and further variants (called de Finettian-Jeffrey in part I, after \citealt{jeffrey1963indeterminate}).\footnote{As explained in part I, \citet{farrell1979material}'s {\sf F}-variant differs from the Cooper-Cantwell table just on the entry $(\half, 1)$, where it returns the value $\half$: the resulting logic, {\sf F/TT}, is also connexive, but does not support full commutation.} 
However, they differ foremost on Modus Ponens, which is preserved in {\sf CC/TT} but given up in {\sf DF/TT}.

In this second part of our inquiry, we turn to an investigation of the proof theory of {\sf DF/TT} and {\sf CC/TT}. We proceed in three main steps: in Section \ref{S:tableaux}, we give sound and complete tableaux calculi for either logic; in Section \ref{S:sequents}, we present sound and complete sequent calculi; in Section \ref{S:algebra}, finally, we examine the prospect for an algebraic semantics for both {\sf DF/TT} and {\sf CC/TT}. As we shall see, neither logic admits a `nice' algebraic semantics, but there is a sense in which {\sf CC/TT}, despite satisfying Modus Ponens, falls even shorter than {\sf DF/TT} in that regard. We give a discussion of that result and further the comparison between both logics in Section \ref{S:Discussion2}. 

\section{Tableau Calculi}
\label{S:tableaux}

\noindent In this section, we introduce sound and complete tableau calculi for {\sf CC/TT} and {\sf DF/TT}. Tableau calculi are a proof-theoretical formalism that is very close to the semantics. To prove a sentence, tableaux employ \emph{trees} that can be conceptualized as reverse truth tables. In building a tableau, one starts from the assumption that certain sentences $A_0, \ldots, A_n$ have certain semantic values, and iteratively works out all the value assignments to the sub-sentences of $A_0, \ldots, A_n$ that result from the initial assignment. In the propositional case, this process always terminates after a finite number of steps, resulting in either an \emph{open} or a \emph{closed} tableau: in the former case, the initial assignment is possible according to the chosen semantics, whereas in the latter it is not. Therefore, in order to prove that $A$ follows from a (finite) set of sentences $\Gamma$ in a tableaux system, one shows that \emph{all} the tableaux resulting from the initial assignments in which all the sentences in $\Gamma$ have a designated value but $A$ does not are closed. 

\subsection{Tableau calculus for {\sf CC/TT}}

\noindent The {\sf CC/TT} tableau calculus, in symbols {\sf CC/TTt}, is given by the following tableau construction rules: 

\begin{center}
  \begin{minipage}{0.95\linewidth}
    \centering
\begin{itemize}
\item[] 
\begin{tikzpicture}[baseline=0pt]
\Tree
[.\node(A){$\neg A:1$};
[.\node(B){$A:0$};
]
]
\end{tikzpicture}
\hspace{5pt}
\begin{tikzpicture}[baseline=0pt]
\Tree
[.\node(A){$\neg A : 0$};
[.\node(B){$A:1$};
]
]
\end{tikzpicture}
\hspace{5pt}
\begin{tikzpicture}[baseline=0pt]
\Tree
[.\node(A){$\neg A : \nicefrac{1}{2}$};
[.\node(B){$A : \nicefrac{1}{2}$};
]
]
\end{tikzpicture}
\item[] 
\item[] 
\begin{tikzpicture}[baseline=0pt]
\Tree
[.\node(A){$A \wedge B : 1$};
[.\node(B){$A : 1$ ; $B : 1$};
]
] 
\end{tikzpicture} 
\hspace{5pt}
\begin{tikzpicture}[baseline=0pt]
\Tree
[.\node(A){$A \wedge B: 0$};
[.\node(B){$A : 0$};
]
[.\node(C){$B : 0$};
]
]\end{tikzpicture}
\hspace{5pt}
\begin{tikzpicture}[baseline=0pt]
\Tree
[.\node(A){$A \wedge B : \nicefrac{1}{2}$};
[.\node(B){$A : 1$ ; $B : \nicefrac{1}{2}$};
]
[.\node(C){$A : \nicefrac{1}{2}$ ; $B : \nicefrac{1}{2}$};
]
[.\node(D){$A : \nicefrac{1}{2}$ ; $B : 1$};
]
]
\end{tikzpicture}
\item[] 
\item[]
\begin{tikzpicture}[baseline=0pt]
\Tree
[.\node(A){$A \rightarrow B : 1$};
[.\node(B){$A : 1$ ; $B : 1$};
]
[.\node(C){$A : \nicefrac{1}{2}$ ; $B : 1$};
]
]
\end{tikzpicture} 
\hspace{5pt}
\begin{tikzpicture}[baseline=0pt]
\Tree
[.\node(A){$A \rightarrow B : 0$};
[.\node(B){$A : 1$ ; $B : 0$};
]
[.\node(C){$A : \nicefrac{1}{2}$ ; $B : 0$};
]
]
\end{tikzpicture}
\hspace{5pt}
\begin{tikzpicture}[baseline=0pt]
\Tree
[.\node(A){$A \rightarrow B : \nicefrac{1}{2}$};
[.\node(B){$A : 0$};
]
[.\node(C){$B : \nicefrac{1}{2}$};
]
]
\end{tikzpicture}
\end{itemize}
\end{minipage}
\end{center}
Tableau rules are essentially versions of the truth-table semantics for the target logic, in our case {\sf CC/TT}. To see this, consider the first rule for the conditional, the one having $A \rightarrow B:1$ as premise. From this premise, one derives two nodes as consequences, one labelled with $A:1$ ; $B:1$ and one labelled with $A:\nicefrac{1}{2}$ ; $B:1$. But these two consequences correspond exactly to the conditions for a conditional $A \rightarrow B$ to have value $1$ in the truth tables for the Cooper-Cantwell conditional: $A \rightarrow B$ has value $1$ in these tables if either both $A$ and $B$ have value $1$, or if $A$ has value $\nicefrac{1}{2}$ and $B$ has value $1$ (see Figure \ref{fig:DFCC}). Similar considerations apply to the other rules. 

We now give a precise characterization of the tableaux generated according to the above rules, and of {\sf CC/TTt}-derivability. 

\begin{defn}$\,$
\begin{itemize}
\item[-] For every formula $A$, the \emph{{\sf CC/TTt}-$n$-tableau of $A$} (for $n = 0$, $\nicefrac{1}{2}$, or $1$) is the tree whose root is $A : n$, and that is obtained by applying the rules of {\sf CC/TTt}.
\item[-] For every \emph{finite} set of formulae $\Gamma = \{B_0, \ldots, B_k\}$, the \emph{{\sf CC/TTt}-$\langle n_0; \ldots; n_k\rangle$-tableau of $A$} (for $n_i = 0$, $\nicefrac{1}{2}$, or $1$, and $i \in \{0, \ldots, k\}$) is the tree whose root is $\langle B_0 : n_0 ; \ldots ; B_k : n_k\rangle$, obtained by applying the rules of {\sf CC/TTt}.\footnote{To ensure uniqueness in the definition of tableaux for more than one sentence, one should fix a convention for the order in which the {\sf CC/TTt} rules are applied. Nothing crucial hinges on this, so we don't specify any such convention for the sake of readability.} 
\end{itemize}
\label{CC/TTtdefinition1}
\end{defn}

Since we are only concerned with the tableau calculus for {\sf CC/TT} in this subsection, we suppress the label `{\sf CC/TTt}' whenever possible, to improve readability.

\begin{defn}$\,$
\begin{itemize}
\item[-] A branch $\mathcal{B}$ of an $n$-tableau is \emph{closed} if, for some formula $A$, there are at least two nodes in $\mathcal{B}$ that have $A : m$ and $A : n$ in their labels, and $n \neq m$. A branch is \emph{open} if it is not closed. 
\item[-] An $n$-tableau is \emph{closed} if all its branches are closed, and \emph{open} otherwise. 
\end{itemize}
\label{CC/TTtdefinition2}
\end{defn}

\begin{defn}
For every \emph{finite} set of formulae $\Gamma$ of cardinality $k$ and every formula $A$, $A$ is \emph{{\sf CC/TTt}-deducible} from a $\Gamma$, in symbols $\Gamma \vdash_{\sf CC/TTt} A$, if and only if all its $\langle n_0; \ldots; n_k \rangle$-tableaux are closed, where we use the indices  $n_1, \ldots, n_{k-1}$ to range over elements in $\{\nicefrac{1}{2},1\}$, and we set $n_{k} = 0$.
\label{CC/TTtdefinition3}
\end{defn}

Before proving soundness and completeness for {\sf CC/TTt}, we give a sample of how to reason in this calculus. In particular, we prove (one direction of) the commutation with negation in {\sf CC/TTt}. The following two closed tableaux establish that $A \rightarrow \neg B$ follows from $\neg (A \rightarrow B)$ in {\sf CC/TTt}. The first tableau shows that $\neg (A \rightarrow B)$ cannot have value $1$ while $A \rightarrow \neg B$ has value $0$. 

\begin{center}
\scriptsize
\begin{tikzpicture}[baseline=0pt]
\Tree
[.\node{$\neg (A \rightarrow B) : 1$ ; $A \rightarrow \neg B : 0$};
[.\node{$A \rightarrow B : 0$ ; $A \rightarrow \neg B : 0$};
[.\node{$A: 1$ ; $B : 0$ ; $A \rightarrow \neg B : 0$};
[.\node{$A: 1$ ; $B : 0$ ; $A: 1$ ; $\neg B : 0$};
[.\node{$A: 1$ ; $B : 0$ ; $A: 1$ ; $B : 1$};
]
]
[.\node(B){$A : 1$ ; $B : 0$ ; $A : \nicefrac{1}{2}$ ; $\neg B : 0$};
]
]
[.\node(B){$A : \nicefrac{1}{2}$ ; $B : 0$ ; $A \rightarrow \neg B: 0$};
[.\node(B){$A: \nicefrac{1}{2}$ ; $B : 0$ ; $A : 1$ ; $\neg B: 0$};
]
[.\node(B){$A : \nicefrac{1}{2}$ ; $B : 0$ ; $A : \nicefrac{1}{2}$ ; $\neg B: 0$};
[.\node(B){$A : \nicefrac{1}{2}$ ; $B : 0$ ; $A : \nicefrac{1}{2}$ ; $B: 1$};
]
]
]
]
]
\end{tikzpicture} 
\end{center}

\noindent The second tableau shows that $\neg (A \rightarrow B)$ cannot have value $\nicefrac{1}{2}$ while $A \rightarrow \neg B$ has value $0$.  

\begin{center}
\footnotesize
\begin{tikzpicture}[baseline=0pt]
\Tree
[.\node{$\neg (A \rightarrow B) : \nicefrac{1}{2}$ ; $A \rightarrow \neg B : 0$};
[.\node{$A \rightarrow B : \nicefrac{1}{2}$ ; $A \rightarrow \neg B : 0$};
[.\node{$A: 0$ ; $A \rightarrow \neg B : 0$};
[.\node{$A: 0$ ; $A : 1$ ; $\neg B : 0$};
]
[.\node(B){$A: 0$ ; $A: \nicefrac{1}{2}$ ; $\neg B : 0$};
]
]
[.\node(B){$B : \nicefrac{1}{2}$ ; $A \rightarrow \neg B : 0$};
[.\node(B){$B : \nicefrac{1}{2}$ ; $A : 1$ ; $\neg B : 0$};
[.\node(B){$B : \nicefrac{1}{2}$ ; $A : 1$ ; $B: 1$};
]
]
[.\node(B){$B: \nicefrac{1}{2}$ ; $A : \nicefrac{1}{2}$ ; $\neg B: 0$};
[.\node(B){$B : \nicefrac{1}{2}$ ; $A : \nicefrac{1}{2}$ ; $B: 1$};
]
]
]
]
]
\end{tikzpicture} 
\end{center}

\noindent The claim that $\neg (A \rightarrow B)$ follows from $A \rightarrow \neg B$ in {\sf CC/TTt} is established in a similar fashion. 

We now prove that {\sf CC/TTt} is sound and complete with respect to {\sf CC/TT}-validity. 

\begin{defn}
A \emph{quasi-{\sf CC}-evaluation} is a non-total function from the formulae of $\mathcal{L}_{\rightarrow}$ to $\{0, \nicefrac{1}{2},1\}$ that is compatible with the {\sf CC} truth tables.
\label{D:quasi}
\end{defn}

\noindent  More compactly, a quasi-{\sf CC}-evaluation is a proper subset of a {\sf CC}-evaluation. For example, the function that sends $p$ and $(p \wedge q)$ to 1 is a quasi-{\sf CC}-evaluation.

\begin{lem}
\label{presound}
For every finite set $\{B_1, \ldots, B_k\}$ of formulae and every {\sf CC}-evaluation $v$, the completed {\sf CC/TTt}-tableau whose root is 
\begin{center}
$B_1 : v(B_1)$ ; $\ldots$ ; $B_k : v(B_k)$
\end{center}
is open, and all partial functions from sentences to $\{1, \nicefrac{1}{2}, 0\}$ induced by its open branches are quasi-{\sf CC}-evaluations.
\end{lem}

\begin{proof}
By induction on the height of the tree. 
\begin{itemize}
    \item[-] The tableau consisting only of the root $B_1 : v(B_1)$ ; $\ldots$ ; $B_k : v(B_k)$ is open. For suppose it is closed. Then, there are at least two sentences $B_i$ and $B_j$ s.t. $B_i = B_j$ but $v(B_i) \neq v(B_j)$, against the hypothesis that $v$ is a {\sf CC}-evaluation: no {\sf CC}-evaluation assigns two different values to the same sentence, because {\sf CC}-evaluations are functions.
    \item[-] Assume by the inductive hypothesis (IH) that the (incomplete) tableau $\mathcal{T}_n$ whose root is $B_1 : v(B_1)$ ; $\ldots$ ; $B_k : v(B_k)$ and that has height $n$ is open, and that its open branches induce quasi-{\sf CC}-evaluations. Suppose also (in contradiction with the lemma to be shown) that the tableau $\mathcal{T}_{n+1}$ of height $n+1$ resulting by applying one tableau rule to the terminal nodes of $\mathcal{T}_n$ is closed. We reason by cases, according to the last rule applied to the nodes in a branch of $\mathcal{T}_n$ (we only do two cases): 
    \begin{itemize}
        \item[$(\wedge)$] Suppose a conjunction rule is applied to a node $\mathtt{v}$ occurring in an open branch $\mathcal{B}_n$ of height $n$ in $\mathcal{T}_n$, and all the branches of height $n+1$ resulting from this application are closed. There are three possibilities: $\mathtt{v}$ has in its label $A \wedge B : 1$, or $A \wedge B : 0$, or $A \wedge B : \nicefrac{1}{2}$.
        \begin{itemize}
            \item[-] If $\mathtt{v}$ has $A \wedge B : 1$ in its label, then there is exactly one successor node $\mathtt{v_1}$ in the resulting branch $\mathcal{B}_{n+1}$ of height $n$, and $\mathtt{v_1}$ has $A : 1$ ; $B : 1$ in its label. If $\mathcal{B}_{n+1}$ is closed as a result of the addition of $\mathtt{v_1}$, this means that there is at least one node $\mathtt{w}$, a predecessor of $\mathtt{v}$, such that: 
            \begin{align*}
                & \mathtt{w} \mbox{ has } A : 0 \mbox{ in its label}, \mbox{ or }\\
                & \mathtt{w} \mbox{ has } A : \nicefrac{1}{2} \mbox{ in its label}, \mbox{ or }\\
                & \mathtt{w} \mbox{ has } B : 0 \mbox{ in its label}, \mbox{ or }\\       
                & \mathtt{w} \mbox{ has } B : \nicefrac{1}{2} \mbox{ in its label}
            \end{align*} 
            Since we assumed that $\mathcal{B}_{n+1}$ is closed, $\mathcal{B}_n$ has a node (namely $\mathtt{v}$) that has $A \wedge B : 1$ in its label, and a node (namely $\mathtt{w}$) whose label is as in one of the cases just listed. By IH, $\mathcal{B}_n$ induces a quasi-{\sf CC}-evaluation. But no quasi-{\sf CC}-evaluation assigns value $1$ to a conjunction and a value different from $1$ to both conjuncts. Contradiction. 
            \item[-] If $\mathtt{v}$ has $A \wedge B : 0$ or $A \wedge B : \nicefrac{1}{2}$ in its label, the reasoning is exactly analogous to the previous case. 
        \end{itemize}
        \item[$(\rightarrow)$] Suppose a conditional rule is applied to a node $\mathtt{v}$ occurring in an open branch $\mathcal{B}_n$ of height $n$ in $\mathcal{T}_n$, and all the branches of height $n+1$ resulting from this application are closed. There are three possibilities: $\mathtt{v}$ has in its label $A \rightarrow B : 1$, or $A \rightarrow B : 0$, or $A \rightarrow B : \nicefrac{1}{2}$.
        \begin{itemize}
            \item[-] If $\mathtt{v}$ has $A \rightarrow B : \nicefrac{1}{2}$ in its label, then there are exactly two branches $\mathcal{B}^1_{n+1}$ and $\mathcal{B}^2_{n+1}$ of height $n+1$ extending $\mathcal{B}_{n}$ with three successor nodes of $\mathtt{v}$, call them $\mathtt{v_1}$ and $\mathtt{v_2}$ such that:
            \begin{align*}
                & \mathtt{v_1} \mbox{ has } A : 0 \mbox{ in its label}\\
                & \mathtt{v_2} \mbox{ has } B : \nicefrac{1}{2} \mbox{ in its label}
            \end{align*}
Since we assumed that $\mathcal{B}^1_{n+1}$ and $\mathcal{B}^2_{n+1}$ are both closed, then $\mathcal{B}_{n}$ has two nodes $\mathtt{w_1}$ and $\mathtt{w_2}$, predecessors of $\mathtt{v}$, such that: 
            \begin{align*}
                & \mathtt{w_1} \mbox { has } A : 1\mbox{ or } A : \nicefrac{1}{2} \mbox{ in its label}\\
                & \mathtt{w_2} \mbox{ has } B : 0 \mbox{ or } B :1\mbox{ in its label}
            \end{align*}
By IH, $\mathcal{B}_n$ induces a quasi-{\sf CC}-evaluation. But no quasi-{\sf CC}-evaluation assigns value $\nicefrac{1}{2}$ to a conditional while assigning any of the following pairs of values to its antecedent and consequent respectively: $\langle 1, 0 \rangle$, $\langle 1, 1 \rangle$, $\langle \nicefrac{1}{2}, 0 \rangle$, and $\langle \nicefrac{1}{2}, 1 \rangle$. Contradiction. 
            \item[-] If $\mathtt{v_0}$ has $A \rightarrow B : 1$ or $A \rightarrow B : 0$ in its label, the reasoning is exactly analogous to the previous case. 
        \end{itemize}
    \end{itemize} 
\end{itemize}
\end{proof}

\begin{prop}[Soundness] For every \emph{finite} set $\Gamma$ of formulae and every formula $A$:
\[\mbox{if } \Gamma \vdash_{\sf CC/TTt} A, \mbox{ then } \Gamma \models_{\sf CC/TT} A\]
\label{CC/TTtsoundness}
\end{prop}

\begin{proof}
We prove the contrapositive. Suppose that $\Gamma \not\models A$, for $\Gamma = \{B_1, \ldots, B_k\}$. Then there is at least one {\sf CC}-evaluation $v$ such that $v(B_1) \in \{1, \nicefrac{1}{2}\}, \ldots, v(B_k) \in \{1, \nicefrac{1}{2}\}$ but $v(A)=0$. Then, by Lemma \ref{presound}, the tree whose root is labeled as
\begin{center}
$B_1 : v(B_1)$ ; $\ldots$ ; $B_k : v(B_k)$ ; $A : 0$
\end{center}
is open. Therefore, not all the trees whose root is labeled as 
\begin{center}
$B_1 : i$ ; $\ldots$ ; $B_{k} : j$ ; $A : 0$
\end{center}
where $i, j\in \{1,\nicefrac{1}{2}\}$, are closed. But this means that $\Gamma \nvdash_{\sf CC/TTt} A$.
\end{proof}

We finally show that {\sf CC/TTt} is complete with respect to {\sf CC/TT}-validity (for inferences with \emph{finite} sets of premises). 

\begin{lem}\label{quasi}
Every open branch of a completed {\sf CC/TTt}-tableau induces a quasi-{\sf CC}-evaluation that has all the formulae appearing in the branch as its domain and assigns to such formulae the values assigned in the labels appearing in the branch.
\end{lem}

\begin{proof}[Proof (Sketch)]
Let $\mathcal{T}$ be a completed {\sf CC/TTt}-tableau with $\mathcal{B}$ an open branch. The branch is finite and it has a unique terminal node $\mathtt{v}$ of the form $p_i : k$. Consider now the partial function that only sends $p_i$ to $k$ (i.e., that is constituted by the single pair $\langle p_i, k \rangle$). This is clearly a quasi-{\sf CC}-evaluation. Call this function $v_0^\mathcal{B}$. Then construct a new function $v_1^\mathcal{B}$ that simply adds to $v_0^\mathcal{B}$ every pair $\langle A, k\rangle$, where $A : k$ is in the label of the predecessor of $\mathtt{v}$ in $\mathcal{B}$. More generally, let $v_{n+1}^\mathcal{B}$ be the function that results from adding to $v_{n+1}^\mathcal{B}$ all the pairs of sentences and values (recall that functions are extensionally construed as pairs) $\langle A, k\rangle$ such that $A : k$ is in the label of the predecessor of $\mathtt{v}$ in $\mathcal{B}$. More generally, let $v_{i+1}^\mathcal{B}$ be the function that results from adding to $v_{i}^\mathcal{B}$ all the pairs of sentences and values (recall that functions are extensionally construed as pairs) $\langle A, k\rangle$ such that $A : k$ is in the label of the predecessor of $\mathtt{v}$ in $\mathcal{B}$. Proceed in this fashion until the root of $\mathcal{T}$ is reached. It is easy to show that the resulting function $v_n^\mathcal{B}$ (for $n+1$ the length of $\mathcal{B}$) is a quasi-{\sf CC}-evaluation.\footnote{More precisely, the construction of $v_n^\mathcal{B}$ is a positive elementary definition that closes at ordinal stage $\omega$. That $v_n^\mathcal{B}$ is a quasi-{\sf CC}-evaluation follows by a straightforward induction, similar to the one used in the proof of Lemma \ref{presound}.}
\end{proof}

\begin{prop}[Completeness]
For every \emph{finite} set $\Gamma$ of formulae and every formula $A$:
\[\mbox{if } \Gamma \models_{\sf CC/TT} A, \mbox{ then } \Gamma \vdash_{\sf CC/TTt} A\]
\label{CC/TTtcompleteness}
\end{prop}

\begin{proof}
We prove the contrapositive. Assume $\Gamma \nvdash_{\sf CC/TTt}A$. By definition this means that not all the {\sf CC/TTt}-tableaux whose root is labeled as 
\begin{center}
$B_1 : i$ ; $\ldots$ ; $B_{k} : j$ ; $A : 0$
\end{center}
are closed, where $i, j\in \{1,\nicefrac{1}{2}\}$. Therefore, at least one such tableau is open. Let $\mathcal{B}$ be an open branch in that tree. By Lemma \ref{quasi}, $\mathcal{B}$ induces a quasi-{\sf CC}-evaluation $v^\mathcal{B}$ such that 
\[v^\mathcal{B}(B_1) \in \{1, \nicefrac{1}{2}\}, \ldots, v^\mathcal{B}(B_{k}) \in \{1, \nicefrac{1}{2}\}, \mbox{ and } v^\mathcal{B}(A) = 0.\] 
Then $v^\mathcal{B}$ can be extended to at least one {\sf CC}-evaluation using Zorn's Lemma. Call one such evaluation $v$. $v$ and $v^\mathcal{B}$ agree on $\Gamma$ and $A$, and therefore \[v(B_1) \in \{1, \nicefrac{1}{2}\}, \ldots, v(B_{k}) \in \{1, \nicefrac{1}{2}\}, \mbox{ and } v(A) = 0.\] But this means that $\Gamma \not\models_{\sf CC/TT} A$.
\end{proof}

\subsection{Tableau calculus for {\sf DF/TT}}

\noindent The tableau calculus for {\sf DF/TT}, in symbols {\sf DF/TTt}, is given by the rules of {\sf CC/TTt}, with the conditional rule replaced by the following one:
\begin{center}
  \begin{minipage}{0.95\linewidth}
    \centering
\begin{itemize}
\item[] 
\begin{tikzpicture}[baseline=0pt]
\Tree
[.\node(A){$A \rightarrow B : 1$};
[.\node(B){$A : 1$ ; $B : 1$};
]
]
\end{tikzpicture} 
\hspace{5pt}
\begin{tikzpicture}[baseline=0pt]
\Tree
[.\node(A){$A \rightarrow B : 0$};
[.\node(B){$A : 1$ ; $B : 0$};
]
]
\end{tikzpicture} 
\hspace{5pt}
\begin{tikzpicture}[baseline=0pt]
\Tree
[.\node(A){$A \rightarrow B : \nicefrac{1}{2}$};
[.\node(B){$A : 0$};
]
[.\node(C){$A : \nicefrac{1}{2}$};
]
[.\node(D){$B : \nicefrac{1}{2}$};
]
]
\end{tikzpicture}
\end{itemize}
\end{minipage}
\end{center}
  
The notions of {\sf DF/TTt}-$n$-tableau, open and closed branch and open and closed $n$-tableau, and {\sf DF/TTt}-deducibility (in symbols $\vdash_{\sf DF/TTt}$) are easily adapted from the corresponding definitions for {\sf CC/TTt} (Definitions \ref{CC/TTtdefinition1}-\ref{CC/TTtdefinition3}).

\begin{prop}[Soundness and completeness]
For every \emph{finite} set $\Gamma$ of formulae and every formula $A$: 
\[\Gamma \vdash_{\sf DF/TTt} A \mbox{ if and only if } \Gamma \models_{\sf DF/TT} A\]
\end{prop}
\noindent The proof is entirely similar to the proof of Propositions \ref{CC/TTtsoundness} and \ref{CC/TTtcompleteness}.

As mentioned at the beginning of \S \ref{S:tableaux}, tableau calculi are very close to truth table semantics. They are also quite informative: their construction determines all the possible truth value assignments that follow from the hypothesis that a given inference is valid. However, tableau calculi are not a particularly convenient formalism to work with. In particular, since tableau calculi are \emph{refutation} calculi, in order to show that $A$ follows from $\Gamma$ in a tableau system, one has to show that the hypothesis that $\Gamma$ holds while $A$ doesn't cannot be maintained. In a classical setting, this amounts to showing that it is not the case that all the sentences in $\Gamma$ can be assigned value $1$ while $A$ is assigned value $0$ by the corresponding tableau. However, in {\sf CC/TT} and {\sf DF/TT} we have three values, two of which are designated, so this is not enough: we have to exclude that all the sentences in $\Gamma$ can be assigned a designated value, that is either $1$ or $\nicefrac{1}{2}$, while $A$ is assigned value $0$. And this requires to consider all the possible combinations of assignments of values $1$ and $\nicefrac{1}{2}$ to sentences in $\Gamma$ (keeping the assignment of value $0$ to $A$ fixed). Of course, as soon as $\Gamma$ contains more than $1$ sentence, showing that $A$ follows from $\Gamma$ requires more than one tableau---more precisely, it requires $2^k$ tableaux, for $k$ the cardinality of $\Gamma$. 
One might avoid this specific problem by considering the conjunction of all the sentences in $\Gamma$, in symbols $\bigwedge \Gamma$, but this move would not really make the calculus more convenient to work with. In fact, rather than having to consider a number of tableaux (depending on the cardinality of $\Gamma$), considering $\bigwedge \Gamma$ would yield a single tableau which, however, is in general much larger than each of the tableaux generated by considering all the possible assignments of designated values to sentences in $\Gamma$. This is because, essentially, one has to iteratively apply the tableau rules to $\bigwedge \Gamma$ until this formula has been decomposed into all its conjuncts -- and there are as many ways to assign value $1$ or $\nicefrac{1}{2}$ to a conjunction as there are ways to assign values $1$ or $\nicefrac{1}{2}$ to the conjuncts. Informally put, considering $\bigwedge \Gamma$ would in general amount to constructing a large tableaux where the `smaller' tableaux generated by the single assignments to the sentences in $\Gamma$ are `sewn together'. 

In addition, tableau calculi are given for inferences with \emph{finite} sets of premises. They can be generalized to the case of infinite sets of premises, but this results in an infinitary formalism, namely a formalism in which one either constructs infinitely many tableaux, or infinitary tableaux, that is, well-founded trees of transfinite ordinal lengths.\footnote{One might argue that inferences with infinite sets of premises do not pose a real problem insofar as the model-theoretical notion of validity $\models$ one aims at capturing is \emph{compact}: for every (possibly infinite) set of formulae $\Gamma$ and every formula $A$, $\Gamma \models A$ if and only if for some finite $\Gamma_0 \subseteq \Gamma$, $\Gamma_0 \models A$. However, if one knows that $A$ follows from an infinite set premises $\Gamma$, compactness only guarantees that \emph{there is} at least one finite $\Gamma_0 \subseteq \Gamma$ from which $A$ follows, but gives no information as to \emph{which} finite subsets of $\Gamma$ entail $A$. Therefore, having a calculus that handles easily infinite sets of premises is crucial when one knows that $A$ follows from an infinite set $\Gamma$ but does not know which finite set of sentences to employ as premises. We also note that, if one's calculus does not allow infinite sets of premises, one cannot derive compactness from completeness, as it is standardly done, and compactness must then be proven with other tools (e.g., via the Ultrafilter Lemma).}

For these reasons, we now present another formalisms to capture {\sf CC/TT}- and {\sf DF/TT}-validity: \emph{many-sided sequent calculi}, in particular three-sided sequent calculi. Three-sided sequent calculi are a generalization of standard sequent calculi: instead of building derivation trees labeled with sequents, the rules of the calculus generate derivation trees labeled with triples of sets of sentences, called three-sided sequents. Unlike tableaux, sequent calculi are not refutation calculi, and therefore any derivation of $A$ from $\Gamma$ establishes that $A$ is provable from $\Gamma$. In addition, sequent calculi handle arbitrary sets of premises, including infinite ones. They can also handle (possibly infinite) sets of conclusions, and therefore generalize {\sf CC/TT}- and {\sf DF/TT}-validity to multiple conclusions. All these advantages have little costs for the intuitiveness of the calculus. Even though one cannot represent in a sequent calculus all the possible outcomes of assigning a given value to a set of sentences, the sequent rules that we are going to use are very close to the tableau rules, and mirror closely the evaluations of their target sentences according to the {\sf CC} and {\sf DF} truth tables.

\section{Three-sided sequent calculi}
\label{S:sequents}

\noindent In this section, we introduce sound and complete three-sided sequent calculi for {\sf CC/TT} and {\sf DF/TT}. Since both {\sf CC/TT} and {\sf DF/TT} are super-logics of {\sf LP} (they extend the latter with a new conditional), we can obtain a sequent calculus by extending an existing calculus for {\sf LP}, in particular the three-sided sequent axiomatization of {\sf LP} provided by \citealt{Ripley2012RSL}.\footnote{{\sf LP}, for the `Logic of Paradox', is a paraconsistent logic adopted in some approaches to the semantic paradoxes. It is the sub-logic of both {\sf CC/TT} and {\sf DF/TT} that results from removing the conditional ($\rightarrow$) from the latter, and it is therefore quite natural to axiomatize {\sf CC/TT} and {\sf DF/TT} over axiomatizations of {\sf LP}. For more on {\sf LP} and some of its developments, see \citet{asen:acal66, priest1979logic,  Priest2006, Goodship1996, Beall2009}.} A three-sided sequent, or a sequent for short, is an object of the form 
\[\Gamma \;|\; \Delta \;|\; \Sigma\]
where $\Gamma$, $\Delta$, and $\Sigma$ are sets of formulae. As above, we focus on the calculus for {\sf CC/TT}, and then indicate how to adapt it to the case of {\sf DF/TT}.

\subsection{Three-sided sequent calculus for {\sf CC/TT}}
\label{S:2.2}

\noindent Let {\sf CC/TTm} be the calculus given by the following principles: \\ 

\noindent Axiom: 

\begin{center}
\AxiomC{$\,$}
\RightLabel{\scriptsize{$\mathsf{SRef}$}}
\UnaryInfC{$\Gamma, A \; | \; \Delta, A \; | \; \Sigma, A$}
\DisplayProof
\end{center}

\noindent Rules: 

\begin{itemize}
\item[] 

\AxiomC{$\Gamma \; | \; \Delta \; | \; \Sigma, A$}
\RightLabel{\scriptsize{$\neg$-0}}
\UnaryInfC{$\Gamma, \neg A \; | \; \Delta \; | \; \Sigma$}
\DisplayProof
\hspace{5pt}
\AxiomC{$\Gamma \; |  \; \Delta, A \; | \; \Sigma$}
\RightLabel{\scriptsize{$\neg$-$\nicefrac{1}{2}$}}
\UnaryInfC{$\Gamma \; | \; \Delta, \neg A \; | \; \Sigma$}
\DisplayProof
\hspace{5pt}
\AxiomC{$\Gamma, A \; |  \; \Delta \; | \; \Sigma$}
\RightLabel{\scriptsize{$\neg$-1}}
\UnaryInfC{$\Gamma \; | \; \Delta \; | \; \Sigma, \neg A$}
\DisplayProof\\
\item[]
\item[]
\AxiomC{$\Gamma, A, B \; |  \; \Delta \; | \; \Sigma $}
\RightLabel{\scriptsize{$\wedge$-0}}
\UnaryInfC{$\Gamma , A \wedge B \; | \; \Delta \; | \; \Sigma$}
\DisplayProof
\hspace{5pt}
\AxiomC{$\Gamma \; | \; \Delta, A \; | \; \Sigma, B$}
\AxiomC{$\Gamma \; | \; \Delta, B \; | \; \Sigma, A$}
\AxiomC{$\Gamma \; | \; \Delta, A, B \; | \; \Sigma$}
\RightLabel{\scriptsize{$\wedge$-$\nicefrac{1}{2}$}}
\TrinaryInfC{$\Gamma \; | \; \Delta, A \wedge B \; | \; \Sigma$}
\DisplayProof\\
\item[]
\begin{center}
\AxiomC{$\Gamma \; | \; \Delta \; | \; \Sigma, A$}
\AxiomC{$\Gamma \; | \; \Delta \; | \; \Sigma, B$}
\RightLabel{\scriptsize{$\wedge$-1}}
\BinaryInfC{$\Gamma \; | \; \Delta \; | \; \Sigma, A \wedge B$}
\DisplayProof\\
\end{center}
\item[]
\item[]
\begin{center}
\AxiomC{$\Gamma \; |  \; \Delta, A \; | \; \Sigma, A$}
\AxiomC{$\Gamma, B \; |  \; \Delta \; | \; \Sigma$}
\RightLabel{\scriptsize{$\rightarrow$-0}}
\BinaryInfC{$\Gamma , A \rightarrow B \; | \; \Delta \; | \; \Sigma$}
\DisplayProof
\hspace{5pt}
\AxiomC{$\Gamma, A \; | \; \Delta, B \; | \; \Sigma$}
\RightLabel{\scriptsize{$\rightarrow$-$\nicefrac{1}{2}$}}
\UnaryInfC{$\Gamma \; | \; \Delta, A \rightarrow B \; | \; \Sigma$}
\DisplayProof\\
\end{center}
\item[]
\item[]
\begin{center}
\AxiomC{$\Gamma \; | \; \Delta, A \; | \; \Sigma, A$}
\AxiomC{$\Gamma \; | \; \Delta \; | \; \Sigma, B$}
\RightLabel{\scriptsize{$\rightarrow$-1}}
\BinaryInfC{$\Gamma \; | \; \Delta \; | \; \Sigma, A \rightarrow B$}
\DisplayProof
\end{center}
\end{itemize}

\noindent Three-sided sequents have an immediate semantic reading. Consider a sequent $\Gamma \; | \; \Delta \; | \; \Sigma$. Intuitively, sentences in $\Gamma$ should be thought of as having value $0$, sentences in $\Delta$ should be thought of as having value $\nicefrac{1}{2}$, and sentences in $\Sigma$ should be thought of as having value $1$. This makes it easy to understand the rationale behind the sequent rules. For example, consider the rule $\neg$-0: if $A$ is thought of as having value $1$ (i.e. it appears in the rightmost position in a sequent), then $\neg A$ is thought of as having value $0$ (and is therefore placed in the leftmost position in the sequent). Similar considerations apply to all the other sequent rules. 

A derivation of a sequent $\Gamma \; | \; \Delta \; | \; \Sigma$ in {\sf CC/TTm} is a tree labeled with sequents, whose leaves are axioms of {\sf CC/TTm} and whose remaining nodes are obtained from their predecessors by applying the {\sf CC/TTm}-rules. Let $\Gamma \vdash_{\sf CC/TTm} \Delta$ be a shorthand for `there is a derivation of $\Gamma \; | \; \Delta \; | \; \Delta$ in {\sf CC/TTm}'. 

\begin{defn}[Satisfaction and Validity]
A {\sf C}-evaluation $v$ \emph{satisfies} a sequent $\Gamma \; | \; \Delta \; | \; \Sigma$ if:
\begin{itemize}
\item[-] there is an $A \in \Gamma$ s.t. $v(A) = 0$, or 
\item[-] there is a $B \in \Delta$ s.t. $v(B) = \nicefrac{1}{2}$, or 
\item[-] there is a $C \in \Sigma$ s.t. $v(C) = 1$.
\end{itemize}
A {\sf C}-evaluation $v$ is a \emph{countermodel} for a sequent $\Gamma \; | \; \Delta \; | \; \Sigma$ if $v$ does not satisfy it. A sequent $\Gamma \; | \; \Delta \; | \; \Sigma$ is {\sf CC/TT}-\emph{valid} if it is satisfied by every {\sf CC/TT}-evaluation. 
\label{CC/TTmdefinition1}
\end{defn}
\noindent Let's expand the notion of {\sf CC/TT}-validity to allow for multiple conclusion, and say that $\Gamma \models_{\sf CC/TT} \Delta$ if every {\sc CC}-evaluation that makes all sentences of $\Gamma$ {\sf T}-true makes at least one sentence in $\Delta$ {\sf T}-true. The following lemma, adapted from \citealt{Ripley2012RSL}, is immediate from the definition of satisfaction and validity.

\begin{lem}
For every sets of formulae $\Gamma$ and $\Delta$:
\[\Gamma \models_{\sf CC/TT} \Delta\ \mbox{ if and only if } \Gamma \; | \; \Delta \; | \; \Delta \mbox{ is {\sf CC/TT}-valid}\]
\end{lem}

Before establishing soundness and completeness for {\sf CC/TTm}, we provide an example of how one can reason with this calculus. More precisely, we show the equivalence of $A \rightarrow \neg B$ and $\neg (A \rightarrow B)$ within it. By the above lemma, this amounts to deriving the sequents $\neg (A \rightarrow B) \; | \; A \rightarrow \neg B \; | \; A \rightarrow \neg B$ and $A \rightarrow \neg B \; | \; \neg (A \rightarrow B) \; | \; \neg (A \rightarrow B)$. In the following examples, we use the empty set symbol $\varnothing$ only in order to make the derivations more readable. 
The following derivation establishes the first sequent:

\begin{center}
\AxiomC{$\,$}
\RightLabel{\scriptsize{$\mathsf{SRef}$}}
\UnaryInfC{$A \; | \; A, B \; | \; A$}
\AxiomC{$\,$}
\RightLabel{\scriptsize{$\mathsf{SRef}$}}
\UnaryInfC{$A \; | \; B, A \; | \; B, A$}
\AxiomC{$\,$}
\RightLabel{\scriptsize{$\mathsf{SRef}$}}
\UnaryInfC{$A, B \; | \; B \; | \; B$}
\RightLabel{\scriptsize{$\neg$-1}}
\UnaryInfC{$A \; | \; B \; | \; B, \neg B$}
\RightLabel{\scriptsize{$\rightarrow$-1}}
\BinaryInfC{$A \; | \; B \; | \; B, A \rightarrow \neg B$}
\RightLabel{\scriptsize{$\rightarrow$-1}}
\BinaryInfC{$A \; | \; B \; | \; A \rightarrow B, A \rightarrow \neg B$}
\RightLabel{\scriptsize{$\neg$-$\nicefrac{1}{2}$}}
\UnaryInfC{$A \; | \; \neg B \; | \; A \rightarrow B, A \rightarrow \neg B$}
\RightLabel{\scriptsize{$\rightarrow$-$\nicefrac{1}{2}$}}
\UnaryInfC{$\varnothing \; | \; A \rightarrow \neg B \; | \; A \rightarrow B, A \rightarrow \neg B$}
\RightLabel{\scriptsize{$\neg$-$0$}}
\UnaryInfC{$\neg(A \rightarrow B) \; | \; A \rightarrow \neg B \; | \; A \rightarrow \neg B$}
\DisplayProof
\end{center}

\noindent The following derivation establishes the second sequent: 

\begin{center}
\AxiomC{$\,$}
\RightLabel{\scriptsize{$\mathsf{SRef}$}}
\UnaryInfC{$A \rightarrow \neg B, A \; | \; A, B \; | \; A$}
\AxiomC{$\,$}
\RightLabel{\scriptsize{$\mathsf{SRef}$}}
\UnaryInfC{$A, B \; | \; A, B \; | \; A$}
\AxiomC{$\,$}
\RightLabel{\scriptsize{$\mathsf{SRef}$}}
\UnaryInfC{$A, B \; | \; B \; | \; B$}
\RightLabel{\scriptsize{$\neg$-0}}
\UnaryInfC{$A, B, \neg B \; | \; B \; | \; \varnothing$}
\RightLabel{\scriptsize{$\rightarrow$-0}}
\BinaryInfC{$A \rightarrow \neg B, A, B \; | \; B \; | \; \varnothing$}
\RightLabel{\scriptsize{$\rightarrow$-0}}
\BinaryInfC{$A \rightarrow \neg B, A \rightarrow B, A \; | \; B \; | \; \varnothing$}
\RightLabel{\scriptsize{$\rightarrow$-$\nicefrac{1}{2}$}}
\UnaryInfC{$A \rightarrow \neg B, A \rightarrow B \; | \; A \rightarrow B \; | \; \varnothing$}
\RightLabel{\scriptsize{$\neg$-$\nicefrac{1}{2}$}}
\UnaryInfC{$A \rightarrow \neg B, A \rightarrow B \; | \; \neg(A \rightarrow B) \; | \; \varnothing$}
\RightLabel{\scriptsize{$\neg$-$1$}}
\UnaryInfC{$A \rightarrow \neg B \; | \; \neg(A \rightarrow B) \; | \; \neg (A \rightarrow B)$}
\DisplayProof
\end{center}

We now proceed to establish soundness and completeness for {\sf CC/TTm}.

\begin{prop}[Soundness]
If $\Gamma \vdash_{\sf CC/TTm} \Delta$, then $\Gamma \models_{\sf DF/TT} \Delta$. 
\end{prop}

\begin{proof}
By induction on the length of the derivation of $\Gamma \; | \; \Delta \; | \; \Delta$. 
\end{proof}

To prove completeness, we prove the following more general result. 

\begin{prop}
\label{P:derivation-countermodel}
For every triple of sets of formulae $\Gamma$, $\Delta$, and $\Sigma$, exactly one of the two following cases is given: 
\begin{itemize}
\item[-] there is a derivation of $\Gamma \; | \; \Delta \; | \; \Sigma$ in {\sf CC/TTm}
\item[-] $\Gamma \; | \; \Delta \; | \; \Sigma$ has a countermodel.
\end{itemize}
\end{prop}

\begin{proof}
We employ the method of Sch\"{u}tte's search trees, adapted to {\sf CC/TTm}.\footnote{Search trees were originally introduced in \citealt{Schutte1956}. See \citet[][Chapter 4]{Pohlers2009} for an application of this method to one-sided sequents for classical logic, \citet{BaazFermullerZach1992, BaazFermullerZach1993} for generalizations and applications to $n$-sided sequents, and \citet{Ripley2012RSL} for an application to {\sf ST}.} For every sequent $\Gamma \; | \; \Delta \; | \; \Sigma$, such method provides the means to construct a tree labeled with sequents which either constitutes a derivation of $\Gamma \; | \; \Delta \; | \; \Sigma$ in {\sf CC/TTm} or can be used to extract a countermodel to $\Gamma \; | \; \Delta \; | \; \Sigma$.

We begin by defining three inductive jumps, that extend a given directed tree labeled with  sequents by applying all the rules of {\sf CC/TTm}. Formally, such a tree is constituted by a pair $\langle \mathtt{N}, \mathtt{S} \rangle$, where $\mathtt{N}$ is the set of nodes and $\mathtt{S}$ is the set of edges, together with a labeling function, that is, a function from $\mathtt{N}$ to their labels (that is, sequents). To simplify our presentation, we identify nodes with their labels, and pairs of nodes with pairs of labels. For every labeled directed tree $\langle \mathtt{N}, \mathtt{S} \rangle$, define the following sets by positive elementary induction:
\begin{itemize}
\item[-] $\mathtt{v_0} \in \mathtt{N}^+$ and $\langle \mathtt{v}, \mathtt{v_0} \rangle \in \mathtt{S}^+$ if: \\

$\mathtt{v} \in \mathtt{N}$ and
\begin{enumerate}
\item[] $\mathtt{v}$ is labeled with $\Gamma, \neg A \; | \; \Delta \; | \; \Sigma$, and $\mathtt{v_0}$ with $\Gamma, \; | \; \Delta \; | \; \Sigma, A$, or 
\item[] $\mathtt{v}$ is labeled with $\Gamma \; | \; \Delta, \neg A \; | \; \Sigma$, and $\mathtt{v_0}$  with $\Gamma, \; | \; \Delta, A \; | \; \Sigma$, or 
\item[] $\mathtt{v}$ is labeled with $\Gamma \; | \; \Delta \; | \; \Sigma, \neg A$, and $\mathtt{v_0}$ with $\Gamma, A \; | \; \Delta \; | \; \Sigma$, or
\item[] $\mathtt{v}$ is labeled with $\Gamma, A \wedge B \; | \; \Delta \; | \; \Sigma$, and $\mathtt{v_0}$ with $\Gamma, A, B \; | \; \Delta \; | \; \Sigma$, or 
\item[] $\mathtt{v}$ is labeled with $\Gamma \; | \; \Delta, A \rightarrow B \; | \; \Sigma$, and $\mathtt{v_0}$ with $\Gamma, A \; | \; \Delta, B \; | \; \Sigma$.
\end{enumerate} 
\item[-] $\mathtt{v_0}, \mathtt{v_1} \in \mathtt{N}^\dagger$ and $\langle \mathtt{v}, \mathtt{v_0} \rangle, \langle \mathtt{v}, \mathtt{v_1} \rangle \in \mathtt{S}^\dagger$ if:\\

$\mathtt{v} \in \mathtt{N}$ and 
\begin{enumerate}
\item[] $\mathtt{v}$ is labeled with $\Gamma \; | \; \Delta \; | \; \Sigma, A \wedge B$, $\mathtt{v_0}$ is labeled with $\Gamma \; | \; \Delta \; | \; \Sigma, A$, and $\mathtt{v_1}$ with $\Gamma \; | \; \Delta \; | \; \Sigma, B$, or
\item[] $\mathtt{v}$ is labeled with $\Gamma \; | \; \Delta \; | \; \Sigma, A \rightarrow B$, $\mathtt{v_0}$ is labeled with $\Gamma \; | \; \Delta, A \; | \; \Sigma, A$, and $\mathtt{v_1}$ with $\Gamma \; | \; \Delta \; | \; \Sigma, B$, or
\item[] $\mathtt{v}$ is labeled with $\Gamma, A \rightarrow B \; | \; \Delta \; | \; \Sigma$, $\mathtt{v_0}$ is labeled with $\Gamma \; | \; \Delta, A \; | \; \Sigma, A$, and $\mathtt{v_1}$ with $\Gamma, B \; | \; \Delta \; | \; \Sigma$.
\end{enumerate}
\item[-] $\mathtt{v_0}, \mathtt{v_1}, \mathtt{v_2} \in \mathtt{N}^\ddagger$ and $\langle \mathtt{v}, \mathtt{v_0} \rangle, \langle \mathtt{v}, \mathtt{v_1} \rangle, \langle \mathtt{v}, \mathtt{v_2} \rangle \in \mathtt{S}^\ddagger$ if:\\

$\mathtt{v} \in \mathtt{N}$ and $\mathtt{v}$ is labeled with $\Gamma \; | \; \Delta, A \wedge B \; | \; \Sigma$, $\mathtt{v_0}$ is labeled with $\Gamma \; | \; \Delta, A \; | \; \Sigma, B$, $\mathtt{v_1}$ is labeled with $\Gamma \; | \; \Delta, B \; | \; \Sigma, A$, and $\mathtt{v_2}$ with $\Gamma \; | \; \Delta, A, B \; | \; \Sigma$.
\end{itemize}
Informally, one can see the jumps $^+$, $^\dagger$, and $^\ddagger$ as corresponding to the operations of extending a given labeled tree to another labeled tree, where the sequents that are added result from applying the rules of {\sf CC/TTm} `upside down', that is, going from a sequent to all its possible premises according to the {\sf CC/TTm} rules. 

Now we construct a \emph{search tree} for a every sequent, that is, a labeled tree where the above jumps are systematically applied as many times as possible. For every sequent $\Gamma \;|\; \Delta \;|\; \Sigma$ define (for a limit ordinal $\delta$):
\begin{align*}
\mathtt{N}_{0} &\; := \Gamma \;|\; \Delta \;|\; \Sigma, &\; \mathtt{S}_{0} &\; := \varnothing\\ 
\mathtt{N}_{\alpha+1} &\; := (\mathtt{N}_{\alpha})^+ \cup (\mathtt{N}_{\alpha})^\dagger \cup (\mathtt{N}_{\alpha})^\ddagger &\; \mathtt{S}_{\alpha+1} &\; := (\mathtt{S}_{\alpha})^+ \cup (\mathtt{S}_{\alpha})^\dagger \cup (\mathtt{S}_{\alpha})^\ddagger \\
\mathtt{N}_{\delta} &\; := \bigcup_{\alpha < \delta} \mathtt{N}_{\alpha} &\; \mathtt{S}_{\delta} &\; := \bigcup_{\alpha < \delta} \mathtt{S}_{\alpha}
\end{align*}
Finally, define (where $\mathsf{Ord}$ is the class of all ordinals):\footnote{Since we construct one search tree per sequent, a more perspicuous notation would indicate the dependence of $\mathtt{N}_{\infty}$ on the starting sequent $\Gamma \; | \; \Delta \; | \; \Sigma$, for example by writing $\mathtt{N}_{\infty}(\Gamma \; | \; \Delta \; | \; \Sigma)$ (similarly for $\mathtt{S}_{\infty}$ and the ordinal stages). We stick to the simpler notation for readability and because the sequent in question is clear from the context.}
\begin{align*}
\mathtt{N}_{\infty} := \bigcup_{\alpha \,\in\, {\sf Ord}} \mathtt{N}_{\alpha} && \mathtt{S}_{\infty} := \bigcup_{\alpha \,\in\, {\sf Ord}} \mathtt{S}_{\alpha}
\end{align*}
The tree $\langle \mathtt{N}_{\infty}, \mathtt{S}_{\infty} \rangle$ is the search tree for $\Gamma \; | \; \Delta \; | \; \Sigma$. We say that $\langle \mathtt{N}_{\infty}, \mathtt{S}_{\infty} \rangle$ is \emph{closed} if all its branches have finite length and have an axiom of {\sf CC/TTm} as their topmost nodes, and that it is \emph{open} otherwise. Clearly, if $\langle \mathtt{N}_{\infty}, \mathtt{S}_{\infty} \rangle$ is closed, this very tree provides a proof of $\Gamma \; | \; \Delta \; | \; \Sigma$, since its topmost nodes are axioms and all the other nodes are obtained from their predecessors by applying {\sf CC/TTm}-rules. Now we show that if $\langle \mathtt{N}_{\infty}, \mathtt{S}_{\infty} \rangle$ is open, one can use it to construct a countermodel for $\Gamma \; | \; \Delta \; | \; \Sigma$. 

Suppose $\langle \mathtt{N}_{\infty}, \mathtt{S}_{\infty} \rangle$ is open, and let $\mathcal{B}$ be an open branch in it. Let $\Gamma_{\infty} \; | \; \Delta_{\infty} \; | \; \Sigma_{\infty}$ be the sequent defined as the union of all the sequents in $\mathcal{B}$. More formally: 
\begin{align*}
\Gamma_{\infty} &\; := \bigcup_{\Gamma_0 \; | \; \Delta_0 \; | \; \Sigma_0 \,\in\, \mathcal{B}} \Gamma_0 & \Delta_{\infty} &\; := \bigcup_{\Gamma_1 \; | \; \Delta_1 \; | \; \Sigma_1 \,\in\, \mathcal{B}} \Delta_1 &\Sigma_{\infty} &\; := \bigcup_{\Gamma_2 \; | \; \Delta_2 \; | \; \Sigma_2 \,\in\, \mathcal{B}} \Sigma_2
\end{align*}

We now have to show that no formula is in $\Gamma_{\infty} \cap \Delta_{\infty} \cap \Sigma_{\infty}$ (otherwise the branch we are constructing might not provide a countermodel). Suppose that there is a formula $A$ and there are sequents $\Gamma_0 \; | \; \Delta_0 \; | \; \Sigma_0$, $\Gamma_1 \; | \; \Delta_1 \; | \; \Sigma_1$, and $\Gamma_2 \; | \; \Delta_2 \; | \; \Sigma_2$ such that $A \in \Gamma_0 \cap \Delta_1 \cap \Sigma_2$. We reason by cases, in order to reach contradictions: 
\begin{itemize}
\item[-] Suppose $A$ is a propositional variable $p_i$. Since $\Gamma_0 \; | \; \Delta_0 \; | \; \Sigma_0$, $\Gamma_1 \; | \; \Delta_1 \; | \; \Sigma_1$, and $\Gamma_2 \; | \; \Delta_2 \; | \; \Sigma_2$ all belong to the same open branch $\mathcal{B}$, then they occur at different heights within $\mathcal{B}$. Suppose without loss of generality that $\Gamma_0 \; | \; \Delta_0 \; | \; \Sigma_0$ occurs at height $n$ (counting upwards the nodes appearing in $\mathcal{B}$ starting from the lowest node, labeled with $\Gamma \; | \; \Delta \; | \; \Sigma$), that $\Gamma_1 \; | \; \Delta_1 \; | \; \Sigma_1$ occurs at height $n + j$, and that $\Gamma_2 \; | \; \Delta_2 \; | \; \Sigma_2$ occurs at height $n + j + k$ (considering different orders would not make a difference). Since $p_i \in \Gamma_0$ and all the rules of ${\sf CC/TTm}$ are context-sharing,\footnote{That is, sequent rules with more than one premises are applied to sequents with identical side-formulae, i.e. sharing the sets of premises $\Gamma$, $\Delta$, and $\Sigma$. See \citet[][64 and following]{TroelstraSchwichtenberg2000} for more details.} $p_i$ is `carried upwards' during the construction of successive stages of $\mathcal{B}$. Therefore, at height $n + j$ we have that $p_i \in \Gamma_1$ and $p_i \in \Delta_1$, and at height $n + j +k$ we have that $p_i \in \Gamma_2$, $p_i \in \Delta_2$, and $p_i \in \Sigma_2$. But this means that $\Gamma_2 \; | \; \Delta_2 \; | \; \Sigma_2$ is an axiom of {\sf CC/TTm}, and that $\mathcal{B}$ is closed. Contradiction. 
\item[-] Suppose $A$ is a complex formula of complexity $n+1$, and assume the claim as IH for formulae of complexity up to $n$. Suppose $A$ is $B \rightarrow C$, and that $\Gamma_0 \; | \; \Delta_0 \; | \; \Sigma_0$ occurs at height $n$, that $\Gamma_1 \; | \; \Delta_1 \; | \; \Sigma_1$ occurs at height $n + j$, and that $\Gamma_2 \; | \; \Delta_2 \; | \; \Sigma_2$ occurs at height $n + j + k$. Then:
\begin{itemize}
\item[$\ast$] $B \in \Delta_0^\prime \cap \Sigma_0^\prime$ and $C \in \Gamma_0^{\prime\prime}$ where $\Gamma_0^{\prime} \; | \; \Delta_1^{\prime} \; | \; \Sigma_0^{\prime}$ and $\Gamma_0^{\prime\prime} \; | \; \Delta_0^{\prime\prime} \; | \; \Sigma_0^{\prime\prime}$ are predecessors of $\Gamma_0 \; | \; \Delta_0 \; | \; \Sigma_0$ in $\langle \mathtt{N}_{\infty}, \mathtt{S}_{\infty} \rangle$, and one of them is in $\mathcal{B}$. 
\item[$\ast$] $B \in \Gamma_1^\prime$ and $C \in \Delta_1^{\prime}$ where $\Gamma_1^{\prime} \; | \; \Delta_1^{\prime} \; | \; \Sigma_1^{\prime}$ is a predecessor of $\Gamma_1 \; | \; \Delta_1 \; | \; \Sigma_1$ in $\langle \mathtt{N}_{\infty}, \mathtt{S}_{\infty} \rangle$, and is in $\mathcal{B}$. 
\item[$\ast$] $B \in \Delta_2^\prime \cap \Sigma_2^\prime$ and $C \in \Sigma_2^{\prime\prime}$ where $\Gamma_2^{\prime} \; | \; \Delta_2^{\prime} \; | \; \Sigma_2^{\prime}$ and $\Gamma_2^{\prime\prime} \; | \; \Delta_2^{\prime\prime} \; | \; \Sigma_2^{\prime\prime}$ are predecessors of $\Gamma_2 \; | \; \Delta_2 \; | \; \Sigma_2$ in $\langle \mathtt{N}_{\infty}, \mathtt{S}_{\infty} \rangle$, and one of them is in $\mathcal{B}$. 
\end{itemize}
Therefore, one of the following is the case: 
\begin{itemize}
\item[(i)] $B \in \Gamma_1^\prime \cap \Delta_0^\prime \cap \Sigma_0^\prime$, where $\Gamma_1^\prime \subseteq \Gamma_{\infty}$, $\Delta_0^\prime \subseteq \Delta_{\infty}$, and $\Sigma_0^\prime \subseteq \Sigma_{\infty}$; or 
\item[(ii)] $B \in \Gamma_1^\prime \cap \Delta_2^\prime \cap \Sigma_2^\prime$, where $\Gamma_1^\prime \subseteq \Gamma_{\infty}$, $\Delta_2^\prime \subseteq \Delta_{\infty}$, and $\Sigma_2^\prime \subseteq \Sigma_{\infty}$; or 
\item[(iii)] $C \in \Gamma_0^{\prime\prime} \cap \Delta_1^\prime \cap \Sigma_2^{\prime\prime}$, where $\Gamma_0^{\prime\prime} \subseteq \Gamma_{\infty}$, $\Delta_1^\prime \subseteq \Delta_{\infty}$, and $\Sigma_2^{\prime\prime} \subseteq \Sigma_{\infty}$;  
\end{itemize} 
But all of (i)-(iii) contradict our IH. The cases of the other connectives are similar.
\end{itemize}

Now that we have shown that the open branch $\mathcal{B}$ is such that the unionsets of the left-, middle-, and right-items in the sequents in $\mathcal{B}$ have an empty intersection, that is that  no formula is in $\Gamma_{\infty} \cap \Delta_{\infty} \cap \Sigma_{\infty}$, we can proceed to extract a proper countermodel from $\Gamma_{\infty}$, $\Delta_{\infty}$, and $\Sigma_{\infty}$. In order to construct our countermodel, we now define the following partial function from formulae to $\{0, \nicefrac{1}{2}, 1\}$ by simultaneous induction:
\begin{align*}
w_{\mathcal{B}}^0(p_i) &:= 
\begin{cases}
0, \mbox{ if } p_i \notin \Gamma_{\infty} \mbox{ and } p \in \Sigma_{\infty} \\
\nicefrac{1}{2}, \mbox{ if } p_i \notin \Delta_{\infty} \mbox{ and } p \in \Gamma_{\infty}\\
1, \mbox{ if } p_i \notin \Sigma_{\infty} \mbox{ and } p \in \Delta_{\infty}
\end{cases}\\
w_{\mathcal{B}}^{\alpha+1}(A) &:= \begin{cases}
0, \mbox{ if } \begin{cases}
A \mbox{ is } \neg B \mbox{ and } w^{\alpha}_{\mathcal{B}}(B) = 1, \mbox{ or}\\
A \mbox{ is } B \wedge C \mbox{ and } w^{\alpha}_{\mathcal{B}}(B) = 0 \mbox{ or } w^{\alpha}_{\mathcal{B}}(C) = 0, \mbox{ or}\\
A \mbox{ is } B \rightarrow C \mbox{ and } w^{\alpha}_{\mathcal{B}}(B) = 1 \mbox{ or } w^{\alpha}_{\mathcal{B}}(B) = \nicefrac{1}{2}, \\
\phantom{A \mbox{ is } B \rightarrow C \mbox{ and }}\mbox{ and } w^{\alpha}_{\mathcal{B}}(C) = 0\\
\end{cases}\\\\
\nicefrac{1}{2}, \mbox{ if } \begin{cases}
A \mbox{ is } \neg B \mbox{ and } w^{\alpha}_{\mathcal{B}}(B) = \nicefrac{1}{2}, \mbox{ or}\\
A \mbox{ is } B \wedge C \mbox{ and } w^{\alpha}_{\mathcal{B}}(B) = 1 \mbox{ and } w^{\alpha}_{\mathcal{B}}(C) = \nicefrac{1}{2}, \\
\phantom{A \mbox{ is } B \wedge C }\mbox{ or } w^{\alpha}_{\mathcal{B}}(B) = \nicefrac{1}{2} \mbox{ and } w^{\alpha}_{\mathcal{B}}(C) = 1,  \\
\phantom{A \mbox{ is } B \wedge C }\mbox{ or } w^{\alpha}_{\mathcal{B}}(B) = \nicefrac{1}{2} \mbox{ and } w^{\alpha}_{\mathcal{B}}(C) = \nicefrac{1}{2} \\
A \mbox{ is } B \rightarrow C \mbox{ and } w^{\alpha}_{\mathcal{B}}(B) = 0 \mbox{ or } w^{\alpha}_{\mathcal{B}}(C) = \nicefrac{1}{2}\\
\end{cases}\\\\
1, \mbox{ if } \begin{cases}
A \mbox{ is } \neg B \mbox{ and } w^{\alpha}_{\mathcal{B}}(B) = 0, \mbox{ or}\\
A \mbox{ is } B \wedge C \mbox{ and } w^{\alpha}_{\mathcal{B}}(B) = 1 \mbox{ and } w^{\alpha}_{\mathcal{B}}(C) = 1, \mbox{ or}\\
A \mbox{ is } B \rightarrow C \mbox{ and } w^{\alpha}_{\mathcal{B}}(B) = 1 \mbox{ or } w^{\alpha}_{\mathcal{B}}(B) = \nicefrac{1}{2},\\
\phantom{A \mbox{ is } B \rightarrow C \mbox{ and }}\mbox{ and } w^{\alpha}_{\mathcal{B}}(C) = 1\\
\end{cases}\\
\end{cases}
\end{align*}
By construction, $w_{\mathcal{B}}$ is a quasi-{\sf C}-evaluation,\footnote{See Definition \ref{D:quasi}. That $w_{\mathcal{B}}$ is a quasi-{\sf C}-evaluation can be rigorously shown by an induction similar to the one employed in the proof of Lemma \ref{presound}.} and it can be extended to a {\sf C}-evaluation $w$ that agrees with $w_{\mathcal{B}}$ on $\Gamma_\infty \cup \Delta_\infty \cup \Sigma_\infty$. Therefore, for every $A \in \Gamma_\infty$, $B \in \Delta_\infty$, and $C \in \Sigma_\infty$:
\begin{align*}
w(A) \neq 0 && w(B) \neq \nicefrac{1}{2} && \vspace{5pt} w(C) \neq 1
\end{align*}
But since $\Gamma \subseteq \Gamma_\infty$, $\Delta \subseteq \Delta_\infty$, and $\Sigma \subseteq \Sigma_\infty$, for every $A \in \Gamma$, $B \in \Delta$, and $C \in \Sigma$: 
\begin{align*}
w(A) \neq 0 && w(B) \neq \nicefrac{1}{2} && \vspace{5pt} w(C) \neq 1
\end{align*}
Therefore, there is a {\sf C}-evaluation that does not satisfy $\Gamma \; | \; \Delta \; | \; \Sigma$, namely $w$. 
\end{proof}

A few observations on the functions $w_{\mathcal{B}}^\alpha$ are in order. First of all, the definition of $w_{\mathcal{B}}^0$ is in part arbitrary, as other choices of truth value assignments to propositional variables would have been possible. In order to get a countermodel, one just needs a function that (i) assigns to the propositional variables in $\Gamma_{\infty}$, $\Delta_{\infty}$, and $\Sigma_{\infty}$ a value that is incompatible with the corresponding position of such variables in the union sequent (and clearly there is more than one choice here) and that (ii) is a quasi-{\sf C}-evaluation. Notice moreover that the construction of every $w_{\mathcal{B}}^\alpha$ is by simultaneous induction, but every $w_{\mathcal{B}}^\alpha$ is \emph{inductive in} $\Gamma_{\infty}$, $\Delta_{\infty}$, and $\Sigma_{\infty}$, since these sets occur also negatively in the definition of $w_{\mathcal{B}}^0$.\footnote{See \citet[][17 and following]{Moschovakis1974}.} 
This seems unavoidable: there seems to be no definition of `having value $1$ if not in $\Gamma_{\infty}$, $\nicefrac{1}{2}$ if not in $\Delta_{\infty}$, and $0$ if not in $\Sigma_{\infty}$' that yields a function and that is positive in $\Gamma_{\infty}$, $\Delta_{\infty}$, and $\Sigma_{\infty}$. However, this causes no problem as far as the existence and uniqueness of $w_{\mathcal{B}}$ is concerned, since the existence and uniqueness of $\Gamma_{\infty}$, $\Delta_{\infty}$, and $\Sigma_{\infty}$ is immediate by their definition. 

Finally, notice that we gave a simplified inductive construction for $w_{\mathcal{B}}$. More specifically, we define $w_{\mathcal{B}}$ directly as a function rather than as a positive elementary set of pairs of sentences and values (then one would have had to show that such set is, indeed, a function). Giving a proper positive elementary definition of $w_{\mathcal{B}}$ would make it clearer that its construction is by simultaneous induction, but would be significantly less readable.

A completeness theorem for {\sf CC/TTm} is now immediate from Proposition \ref{P:derivation-countermodel}. 

\begin{prop}[Completeness]
For every set $\Gamma$ of formulae and every formula $A$:
\[\mbox{ if } \Gamma \models_{\sf CC/TT} \Delta, \mbox{ then } \Gamma \vdash_{\sf CC/TTm} \Delta.\] 
\end{prop}

\subsection{Three-sided sequent calculus for {\sf DF/TT}}

\noindent The three-sided sequent calculus for {\sf DF/TT}, in symbols {\sf DF/TTm}, is given by the rules of {\sf CC/TTm}, with the conditional rules replaced by the following ones: 

\begin{itemize}
\item[]
\begin{center}
\AxiomC{$\Gamma \; |  \; \Delta \; | \; \Sigma, A$}
\AxiomC{$\Gamma, B \; |  \; \Delta \; | \; \Sigma$}
\RightLabel{\scriptsize{$\rightarrow$-0}}
\BinaryInfC{$\Gamma , A \rightarrow B \; | \; \Delta \; | \; \Sigma$}
\DisplayProof
\hspace{5pt}
\AxiomC{$\Gamma, A \; | \; \Delta, A, B \; | \; \Sigma$}
\RightLabel{\scriptsize{$\rightarrow$-$\nicefrac{1}{2}$}}
\UnaryInfC{$\Gamma \; | \; \Delta, A \rightarrow B \; | \; \Sigma$}
\DisplayProof\\
\end{center}
\item[]
\item[]
\begin{center}
\AxiomC{$\Gamma \; | \; \Delta \; | \; \Sigma, A$}
\AxiomC{$\Gamma \; | \; \Delta \; | \; \Sigma, B$}
\RightLabel{\scriptsize{$\rightarrow$-1}}
\BinaryInfC{$\Gamma \; | \; \Delta \; | \; \Sigma, A \rightarrow B$}
\DisplayProof
\end{center}
\end{itemize}
The notions of {\sf DF/TTm}-derivability, as well as of satisfaction and validity of a three-sided sequent are immediate  from the corresponding definitions for {\sf CC/TTm} (Definition \ref{CC/TTmdefinition1}). 

\begin{prop}[Soundness and completeness]
For every set $\Gamma$ of formulae and every formula $A$: 
\[\Gamma \vdash_{\sf DF/TTm} A \mbox{ if and only if } \Gamma \models_{\sf DF/TT} A\]
\label{P:soundnesscompletenessDF}
\end{prop}
\noindent The proof is entirely similar to the proof of soundness and completeness for {\sf CC/TTm}.

\section{Algebraic semantics}\label{S:algebra}

\noindent In this section, we explore the algebraic structures that correspond to {\sf DF/TT} and {\sf CC/TT}, and investigate the prospects for an algebraic semantics of these two logics. We begin by recalling some  structures, and introducing the algebraic counterparts of {\sf DF/TT}. We start with {\sf DF/TT} because, as will be clear in Subsection \ref{S:algebraic-sound-compl-CC}, it is algebraically significantly more tractable than {\sf CC/TT}. We use overlined uppercase Latin letters ($\overline{A}$, $\overline{B}$, $\overline{C}$, $\ldots$) to range over sets (supports of algebraic structures) in order to avoid possible confusions with meta-variables for $\Lc$-formulae, and boldface characters to indicate designated elements of the supports of algebraic structures (${\bf 1}$, ${\bf 0}$, $\nicefrac{\bf 1}{\bf 2}$, $\ldots$), in order to avoid possible confusions with truth values in truth table semantics.

\subsection{De Finetti algebrae}


\noindent We begin by some basic definitions, which will be needed for the algebraic semantics for {\sf DF/TT}.

\begin{defn}
\label{structures}
A structure $\mathcal{A} = \langle \overline{A}, \sqcap, \sqcup, {\bf 0}, {\bf 1} \rangle$, where $\overline{A}$ is a set and ${\bf 0}, {\bf 1} \in \overline{A}$, is a \emph{distributive bounded lattice} if for every $a, b, c \in \overline{A}$:
\begin{itemize}
\item[-] The \emph{lattice} conditions are satisfied: 
\begin{itemize}
\item[] $a \sqcap b = b \sqcap a$ and $a \sqcup b = b \sqcup a$ (commutativity)
\item[] $a \sqcap (b \sqcap c) = (a \sqcap b) \sqcap c$ and $a \sqcup (b \sqcup c) = (a \sqcup b) \sqcup c$ (associativity)
\item[] $a \sqcup (b \sqcap a) = a$ and $a \sqcap (b \sqcup a) = a$ (absorption)
\end{itemize}
\item[-] The lattice is \emph{bounded}:
\begin{itemize}
\item[] $a \sqcup \mathbf{0} = a$ 
\item[] $a \sqcap \mathbf{1} = a$
\end{itemize}
\item[-] The lattice is \emph{distributive}:
\begin{itemize}
\item[] $a \sqcap (b \sqcup c) = (a \sqcap b) \sqcup (a \sqcap c)$
\item[] $a \sqcup (b \sqcap c) = (a \sqcup b) \sqcap (a \sqcup c)$
\end{itemize}
\end{itemize}

\noindent For every lattice $\mathcal{A} = \langle \overline{A}, \sqcap, \sqcup  \rangle$ the \emph{order induced by $\mathcal{A}$} is the binary relation $\sqsubseteq \;\subseteq (\overline{A} \times \overline{A})$ such that:
\[a \sqsubseteq b \mbox{ if and only if } a \sqcap b = a \mbox{ if and only if } a \sqcup b = b.\footnote{Notice that $\sqsubseteq$ is transitive. Suppose that $a \sqsubseteq b$ and $b \sqsubseteq c$. Then, $a \sqcup c = (a \sqcap b) \sqcup c = c \sqcup (a \sqcap b) = (c \sqcup a) \sqcap (c \sqcup b) = (c \sqcup a) \sqcap (b \sqcup c) = (c \sqcup a) \sqcap c = c$ (these identities follow from the assumption that $a \sqsubseteq b$, commutativity, distributivity, the assumption that $b \sqsubseteq c$, distributivity again, and absorption, in this order). The proof for $a \sqcap c = a$ is similar. Therefore, $a \sqsubseteq c$.}\] 
An \emph{involution} on a lattice $\mathcal{A}$ is a unary operation $-$ s.t. for every $a, b \in \overline{A}$:
\begin{itemize}
\item[] If $a \sqsubseteq b$, then $-b \sqsubseteq -a$, and 
\item[] $--a = a$
\end{itemize}
A bounded, distributive, involutive lattice $\mathcal{A} = \langle \overline{A}, \sqcap, \sqcup, -, \mathbf{0}, \mathbf{1} \rangle$ is a \emph{De Morgan algebra} if for every $a, b \in \overline{A}$:
\begin{itemize}
\item[] $-(a \sqcap b) = (-a\sqcup-b)$, and 
\item[] $-(a \sqcup b) = (-a\sqcap-b)$
\end{itemize}
A De Morgan algebra $\mathcal{A} = \langle \overline{A}, \sqcap, \sqcup, -, \mathbf{0}, \mathbf{1} \rangle$ is \emph{Kleene} if, for every $a, b \in \overline{A}$: 
\begin{itemize}
\item[] $a \sqcap -a \sqsubseteq b \sqcup -b$
\end{itemize}
A \emph{relative pseudocomplementation} on a lattice $\mathcal{A}$ is a binary operation $\rightarrowtail$ s.t. for every $a, b \in \overline{A}$: 
\begin{itemize}
\item[] $a \sqcap c \sqsubseteq b$ if and only if $c \sqsubseteq a \rightarrowtail b$ for all $c \in \overline{A}$.
\end{itemize}
A relatively pseudocompletemented Kleene algebra $\mathcal{A} = \langle \overline{A}, \sqcap, \sqcup, -, \rightarrowtail, \mathbf{0}, \mathbf{1} \rangle$ is an \emph{$\L{}3$ algebra} if for every $a \in \overline{A}$:
\begin{itemize}
\item[] $(a \rightarrowtail \mathbf{0}) \sqcup (- a \rightarrowtail a) = \mathbf{1}$
\end{itemize}
An $\L{}3$ algebra $\mathcal{A} = \langle \overline{A}, \sqcap, \sqcup, -, \rightarrowtail, \dfcond, \mathbf{0}, \mathbf{1}, \nicefrac{{\bf 1}}{{\bf 2}} \rangle$ is  \emph{de Finetti} if: 
\begin{itemize}
\item[] There is a distinguished element $\nicefrac{{\bf 1}}{{\bf 2}} \in \overline{A}$ s.t. $- \nicefrac{{\bf 1}}{{\bf 2}} = \nicefrac{{\bf 1}}{{\bf 2}}$, and 
\item[] There is an operation $\dfcond$ defined on $\overline{A} \times \overline{A}$ s.t. $a \dfcond b = (\nicefrac{{\bf 1}}{{\bf 2}} \sqcap -a) \sqcup (a \sqcap b)$.\footnote{The definition of $\L{}3$ algebrae follows  \citet[][517-518]{Milne2004}, and so does the characterization of the algebraic counterpart of the de Finetti conditional over them. We note that Milne considers algebrae of conditional events, while we consider arbitrary supports. Nothing crucial hinges on this.}
\end{itemize}
\end{defn}


\noindent Some remarks on de Finetti algebrae are in order. First, we have defined them over $\L{}3$-algebrae (also known as \L{}ukasiewicz (or Moisil-\L{}ukasiewicz) trivalent algebrae), but other options are possible, including  MV$_3$-algebrae.\footnote{See \citet[][Ch. 9.2]{Gottwald2001} and \citet[][Ch. 5]{Malinowski1993}. Indeed, while MV$_n$ algebrae provide suitable algebraic counterparts for every $n$-valued \L{}ukasiewicz logic (and MV-algebrae algebraically characterize \L{}ukasiewicz continuum-valued logic), \L{}ukasiewicz algebrae only succeed in capturing three- and four-valued cases.} We have adopted $\L{}3$-algebrae both because they are simpler than MV$_n$-algebrae, and in order to better relate our presentation and results to the elegant formalization and the results of \citealt{Milne2004}. Second, de Finetti algebrae have a paraconsistent flavour, suggested by the behaviour of the element $\nicefrac{{\bf 1}}{{\bf 2}}$. Such flavour is more vividly expressed by noticing that they are both special cases of \emph{LP algebrae}. In the characterization offered by \citet{Pynko1995}, a Kleene algebra $\mathcal{A} = \langle \overline{A}, \sqcap, \sqcup, -, \mathbf{0}, \mathbf{1} \rangle$ is \emph{LP} if it has an \emph{inconsistent proper filter} on its support, that is, if there is an $\overline{F} \subset \overline{A}$ s.t. for every $a, b, \in \overline{A}$ and for some $c \in \overline{F}$
\begin{itemize}
\item[(i)] if $a \in \overline{F}$ and $a \sqsubseteq b$, then $b \in \overline{F}$,
\item[(ii)] if $a, b \in \overline{F}$, then $a \sqcap b \in \overline{F}$,
\item[(iii)] both $c \in \overline{F}$ and $-c \in \overline{F}$.
\end{itemize}
It is easily seen that de Finetti algebrae are LP. Let $\mathcal{A}$ be a de Finetti algebra with support $\overline{A}$. The set $\{a \in \overline{A} \,|\, \nicefrac{{\bf 1}}{{\bf 2}} \sqsubseteq a\} \subset \overline{A}$ provides the required inconsistent proper filter. 
\begin{itemize}
\item[] (i) is immediate, because $\sqsubseteq$ is transitive.  
\item[] As for (ii), assume that $\nicefrac{{\bf 1}}{{\bf 2}} \sqsubseteq a$ and that $\nicefrac{{\bf 1}}{{\bf 2}} \sqsubseteq b$. By Definition \ref{structures}, this assumption entails that $(\nicefrac{{\bf 1}}{{\bf 2}} \sqcup a) = a$ and that $(\nicefrac{{\bf 1}}{{\bf 2}} \sqcup b) = b$. But since we have $\nicefrac{{\bf 1}}{{\bf 2}} \sqcup (a \sqcap b) = (\nicefrac{{\bf 1}}{{\bf 2}} \sqcup a) \sqcap (\nicefrac{{\bf 1}}{{\bf 2}} \sqcup b)$ by distributivity, we also have $(\nicefrac{{\bf 1}}{{\bf 2}} \sqcup a) \sqcap (\nicefrac{{\bf 1}}{{\bf 2}} \sqcup b) = a \sqcap b$ by our assumption, that is, $\nicefrac{{\bf 1}}{{\bf 2}} \sqsubseteq a \sqcap b$.
\item[] As for (iii), notice that both $\nicefrac{{\bf 1}}{{\bf 2}}$ and $-\nicefrac{{\bf 1}}{{\bf 2}}$ are in $\{a \in \overline{A} \,|\, \nicefrac{{\bf 1}}{{\bf 2}} \sqsubseteq a\}$.
\end{itemize}

\subsection{Algebraic semantics for {\sf DF/TT}}

In order to prove algebraic soundness and completeness for {\sf DF/TT}, we construct the Lindenbaum-Tarski algebra of a set of formulae, for {\sf DF/TTm}-deducibility. Therefore, we first isolate the relation of {\sf DF/TTm}-provable equivalence (where `equivalence' is formalized via the {\sf DF}-biconditional). 

\begin{defn}
For every $\Gamma \subseteq {\sf For}$, let $\sim^{\sf df}_{\Gamma} \, \subseteq {\sf For} \times {\sf For}$ be the relation defined as follows: 
\[A \sim^{\sf df}_{\Gamma} B \mbox{ if and only if } \Gamma \vdash_{\sf DF/TTm} A \leftrightarrow B\]
\end{defn}

This definition, however, does not partition the set of formulae into equivalence classes, but only into sets that have weaker closure conditions. 

\begin{lem}
\label{noequiv}
In general, $\sim^{\sf df}_{\Gamma}$ is \emph{not} an equivalence relation on ${\sf For} \times {\sf For}$. 
\end{lem}

\begin{proof}
Reflexivity and symmetry hold, since $\Gamma \vdash_{\sf DF/TTm} A \leftrightarrow A$, and if $\Gamma \vdash_{\sf DF/TTm} A \leftrightarrow B$, then also $\Gamma \vdash_{\sf DF/TTm} B \leftrightarrow A$. However, transitivity fails, for otherwise $\vdash_{\sf DF/TTm}$ would be unsound (consider a {\sf DF/TT}-evaluation $v$ in which $v(A) = 1$, $v(B) = \nicefrac{1}{2}$, and $v(C) = 0$). 
\end{proof}

\noindent As the above proof shows, the failure of transitivity for {\sf DF/TTm}-provable equivalence is closely connected to the failure of Modus Ponens for {\sf DF/TT}. However, even though $\sim^{\sf df}_{\Gamma}$ is not an equivalence relation on ${\sf For} \times {\sf For}$, we will see that it is sufficiently well-behaved to support an application of the Lindenbaum-Tarski method. More specifically, due to failure of transitivity, $\sim^{\sf df}_{\Gamma}$ fails to partition ${\sf For} \times {\sf For}$ into equivalence classes. Nevertheless, we can still use $\sim^{\sf df}_{\Gamma}$ to define sets of sets of formulae which still support an application of the Lindenbaum-Tarski method, and therefore a proof of algebraic completeness. This is done in the following definition. 

\begin{defn}
For every $\{\Gamma, A\} \subseteq {\sf For}$, let $[A]^{\sf df}_{\Gamma}$ denote the set of formulae that are provably ${\sf DF/TTm}$-equivalent to $A$. The \emph{quotient} induced by $\sim^{\sf df}_{\Gamma}$ on ${\sf For}$, in symbols ${\sf For} \slash\hspace{-4pt}\sim^{\sf df}_{\Gamma}$, is the set of sets in which $\sim^{\sf df}_{\Gamma}$ subdivides ${\sf For}$.
\end{defn}
\noindent More formally, ${\sf For} \slash\hspace{-4pt}\sim^{\sf df}_{\Gamma}$ is the set of sets $\{A \in {\sf For} \,|\, \Gamma \vdash_{\sf DF/TTm} A \leftrightarrow B \mbox{ for some } B \in {\sf For}\}$.  Since we are only concerned with de Finetti algebrae in this subsection, we drop the superscript $^{\sf df}$ and simply write $[A]_\Gamma$, in order to improve readability. 

\begin{defn}
\label{Lindenbaum-Tarski}
The \emph{de Finetti-Lindenbaum-Tarski algebra of $\Gamma$} is the structure
\[\mathcal{D}(\Gamma) = \langle {\sf For} \slash\hspace{-4pt}\sim_{\Gamma}, \sqcap_{\Gamma}, \sqcup_{\Gamma}, -_{\Gamma}, \dfcond_{\Gamma}, \mathbf{0}_{\Gamma}, \nicefrac{{\bf 1}}{{\bf 2}}_{\Gamma}, \mathbf{1}_{\Gamma} \rangle\]
where: 

\begin{center}
\begin{minipage}{6cm}
\begin{itemize}
\item[] $[A]_{\Gamma} \sqcap_{\Gamma} [B]_{\Gamma} := [A \wedge B]_{\Gamma}$
\item[] $-_{\Gamma} [A]_{\Gamma} := [\neg A]_{\Gamma}$
\item[] $[\bot]_{\Gamma} := \mathbf{0}_{\Gamma}$
\end{itemize}
\end{minipage} 
\begin{minipage}{6cm}
\begin{itemize}
\item[] $[A]_{\Gamma} \sqcup_{\Gamma} [B]_{\Gamma} := [A \vee B]_{\Gamma}$
\item[] $ [A]_{\Gamma} \dfcond_{\Gamma} [B]_{\Gamma} := [A \rightarrow B]_{\Gamma}$
\item[] $[\top]_{\Gamma} := \mathbf{1}_{\Gamma}$
\end{itemize}
\end{minipage} 
\[[\bot \rightarrow \top]_{\Gamma} := \nicefrac{{\bf 1}}{{\bf 2}}_{\Gamma}\]
\end{center}
\end{defn}

\noindent As shown by Lemma \ref{noequiv}, $\sim_{\Gamma}$ is not an equivalence relation, and the sets $[A]_{\Gamma}$ are not equivalence classes. Therefore, there is no guarantee that every formula belongs to exactly one of the elements in ${\sf For} \slash\hspace{-4pt}\sim_{\Gamma}$. So, we have to prove that the operations that characterize de Finetti-Lindenbaum-Tarski algebrae, that is $\sqcap_{\Gamma}$, $\sqcup_{\Gamma}$, $-_{\Gamma}$, and $\dfcond_{\Gamma}$ are actually \emph{well-defined}, and do not depend on the choice of particular formulae: otherwise $\sqcap_{\Gamma}$, $\sqcup_{\Gamma}$, $-_{\Gamma}$, and $\dfcond_{\Gamma}$ might not be \emph{operations} at all. This is done in the following lemma.

\begin{lem}[Independence from representatives]
For every set $\{\Gamma, A, B, C, D\} \subseteq {\sf For}$, the following holds: 
\begin{itemize}
\item[-] If $A \sim_{\Gamma} B$, then $\neg A \sim_{\Gamma} \neg B$.
\item[-] If $A \sim_{\Gamma} B$ and $C \sim_{\Gamma} D$, then $(A \wedge C) \sim_{\Gamma} (B \wedge D)$.
\item[-] If $A \sim_{\Gamma} B$ and $C \sim_{\Gamma} D$, then $(A \vee C) \sim_{\Gamma} (B \vee D)$.
\item[-] If $A \sim_{\Gamma} B$ and $C \sim_{\Gamma} D$, then $(A \rightarrow C) \sim_{\Gamma} (B \rightarrow D)$.
\end{itemize}
\label{L:indipendence}
\end{lem}

\begin{proof}
We only show the cases of negation and conditional (the others are similar). 
For the case of negation, suppose that there is a set $\{\Gamma, A, B\} \subseteq {\sf For}$ such that $A \sim_{\Gamma} B$ but that it is not the case that $\neg A \sim_{\Gamma} \neg B$. This means that $\Gamma \vdash_{\sf DF/TTm} A \leftrightarrow B$ but $\Gamma \not\vdash_{\sf DF/TTm} \neg A \leftrightarrow \neg B$. By the completeness of {\sf DF/TTm} (Proposition \ref{P:soundnesscompletenessDF}), this means that $\Gamma \models_{\sf DF/TT} A \leftrightarrow B$ but $\Gamma \not\models_{\sf DF/TT} \neg A \leftrightarrow \neg B$. Let $v$ be any {\sf DF}-evaluation that assigns value $1$ or $\nicefrac{1}{2}$ to all the sentences in $\Gamma$, value $1$ or $\nicefrac{1}{2}$ to $A \leftrightarrow B$ but value $0$ to $\neg A \leftrightarrow \neg B$ (if there are no {\sf DF}-evaluations that assign values $1$ or $\nicefrac{1}{2}$ to all the sentences in $\Gamma$, the claim is immediate). A biconditional is assigned value $0$ by a {\sf DF}-evaluation just in case that evaluation assigns value $1$ to one side of the biconditional and $0$ to the other. Suppose without loss of generality that $v(\neg A) = 1$ and $v(\neg B) = 0$. Since $v$ is a {\sf DF}-evaluation, $v(A) = 0$ and $v(B) = 1$. But then $v(A \leftrightarrow B) = 0$, against our supposition.\footnote{Notice that the above reasoning breaks down for {\sf CC/TT}, due to the different truth table for the conditional (and hence the biconditional). More specifically, a biconditional $A \leftrightarrow B$ is assigned value $0$ by a {\sf CC}-evaluation if that evaluation assigns value $1$ or $\nicefrac{1}{2}$ to one side of the biconditional and $0$ to the other. In the case of a {\sc CC}-evaluation $v$ that assigns value $\nicefrac{1}{2}$ to $\neg A$ and $0$ to $\neg B$, we have that $v(\neg A \leftrightarrow \neg B) = 0$ but $v(A \leftrightarrow B) = \nicefrac{1}{2}$. See also Lemma \ref{lem:CCsub} below.}

For the case of the conditional, suppose that there is a set $\{\Gamma, A, B, C, D\} \subseteq {\sf For}$ such that $A \sim_{\Gamma} B$ and $C \sim_{\Gamma} D$ but that it is not the case that $(A \rightarrow C) \sim_{\Gamma} (B \rightarrow D)$. This means that $\Gamma \vdash_{\sf DF/TTm} A \leftrightarrow B$ and $\Gamma \vdash_{\sf DF/TTm} C \leftrightarrow D$ but $\Gamma \not\vdash_{\sf DF/TTm} (A \rightarrow C) \leftrightarrow (B \rightarrow D)$. Again by the completeness of {\sf DF/TTm} (Proposition \ref{P:soundnesscompletenessDF}), this means that $\Gamma \models_{\sf DF/TT} A \leftrightarrow B$ and $\Gamma \models_{\sf DF/TT} C \leftrightarrow D$ but $\Gamma \not\models_{\sf DF/TT} (A \rightarrow C) \leftrightarrow (B \rightarrow D)$. Let $v$ be any {\sf DF}-evaluation that assigns value $1$ or $\nicefrac{1}{2}$ to all the sentences in $\Gamma$, value $1$ or $\nicefrac{1}{2}$ to $A \leftrightarrow B$ and to $C \leftrightarrow D$ but value $0$ to $(A \rightarrow C) \leftrightarrow (B \rightarrow D)$ (if there are no {\sf DF}-evaluations that assign values $1$ or $\nicefrac{1}{2}$ to all the sentences in $\Gamma$, the claim is immediate). As noted above, a biconditional is assigned value $0$ by a {\sf DF}-evaluation just in case that evaluation assigns value $1$ to one side of the biconditional and $0$ to the other. Suppose without loss of generality that $v(A \rightarrow C) = 1$ and $v(B \rightarrow D) = 0$. Since $v$ is a {\sf DF}-evaluation, $v(A) = v(C) = v(B) = 1$ and $v(D) = 0$. But then $v(C \leftrightarrow D) = 0$, against our supposition. 
\end{proof}

\begin{lem}
\label{LTDeFinetti}
For every $\Gamma \subseteq {\sf For}$, $\mathcal{D}(\Gamma)$ is a de Finetti algebra. 
\end{lem}

\begin{proof}
It is easy to see that the properties of distributive bounded lattices hold for $\mathcal{D}(\Gamma)$. We do just one case of distributivity in detail.  
\begin{align*}
[A]_{\Gamma} \sqcap_{\Gamma} ([B]_{\Gamma} \sqcup_{\Gamma} [C]_{\Gamma}) \;& = [A]_{\Gamma} \sqcap_{\Gamma} [B \vee C]_{\Gamma} && \mbox{ definition of } \mathcal{D}(\Gamma)\\
\;& = [A \wedge (B \vee C)]_{\Gamma}&&\mbox{ definition of } \mathcal{D}(\Gamma)\\
\;& = [(A \wedge B) \vee (A \wedge C)]_{\Gamma}&& \mbox{ logic } \\
\;& = [A \wedge B]_{\Gamma} \sqcup_{\Gamma} [A \vee C]_{\Gamma}&&\mbox{ definition of } \mathcal{D}(\Gamma)\\
\;& = ([A]_{\Gamma} \sqcap_{\Gamma} [B]_{\Gamma}) \sqcup_{\Gamma} ([A]_{\Gamma} \sqcap_{\Gamma} [C]_{\Gamma})&&\mbox{ definition of } \mathcal{D}(\Gamma)
\end{align*}
The line labeled with `logic' abbreviates the fact that the corresponding identity is proven by the fact that $\Gamma \vdash_{\sf DF/TTm} (A \wedge (B \vee C)) \leftrightarrow ((A \wedge B) \vee (A \wedge C))$. 

As for the involution, suppose that $[A]_{\Gamma} \sqsubseteq_{\Gamma} [B]_{\Gamma}$, where $\sqsubseteq_{\Gamma}$ is the partial order induced on ${\sf For} \slash\hspace{-4pt}\sim^{\sf df}_{\Gamma}$ by $\sqcap_{\Gamma}$ and $\sqcup_{\Gamma}$. By Definition \ref{structures}, this means that $([A]_{\Gamma} \sqcap_{\Gamma} [B]_{\Gamma}) = [A]_{\Gamma}$. However, $(A \wedge B) \leftrightarrow A$ and $A \rightarrow B$ have the same {\sf DF} truth table, and are {\sf DF/TTm}-provably equivalent. Therefore, in particular, $\Gamma \vdash_{\sf DF/TTm} ((A \wedge B) \leftrightarrow A) \leftrightarrow A \rightarrow B$. It follows that $\sqsubseteq_{\Gamma}$ inherits the features of $\dfcond_{\Gamma}$ in $\mathcal{D}(\Gamma)$, and thus the claim is established by contraposition. More precisely:\footnote{Note that the following reasoning does not rely on  uses of Modus Ponens in {\sf DF/TTm}, but can be carried out using the soundness and completeness of the calculus {\sf DF/TTm}.}
\begin{align*}
([A]_{\Gamma} \sqcap_{\Gamma} [B]_{\Gamma}) = [A]_{\Gamma} \; & \mbox{ iff } ([A \wedge B]_{\Gamma}) = [A]_{\Gamma}\\
& \mbox{ iff } \Gamma \vdash_{\sf DF/TTm} (A \wedge B) \leftrightarrow A\\
& \mbox{ iff } \Gamma \vdash_{\sf DF/TTm} A \rightarrow B\\
& \mbox{ iff } \Gamma \vdash_{\sf DF/TTm} \neg B \rightarrow \neg A\\
& \mbox{ iff } \Gamma \vdash_{\sf DF/TTm} (\neg B \wedge \neg A) \leftrightarrow \neg B\\
& \mbox{ iff } ([\neg B \wedge \neg A]_{\Gamma}) = [\neg B]_{\Gamma}\\
& \mbox{ iff } (-_{\Gamma} [B]_{\Gamma} \sqcap_{\Gamma} -_{\Gamma} [A]_{\Gamma}) = -_{\Gamma} [B]_{\Gamma}\\
\end{align*}
The second constraint on involution is also satisfied, since $\Gamma \vdash_{\sf DF/TTm} A \leftrightarrow \neg \neg A$. 

The De Morgan and Kleene properties, as well as the defining equation of $\dfcond_{\Gamma}$, are proven in a similar way. 
\end{proof}


We now provide algebraic counterparts of the notions of evaluation and  {\sf TT}-consequence. More precisely, we provide a local notion of algebraic {\sf TT}-consequence ({\sf TT}-consequence with respect to a single de Finetti algebra) and a global one ({\sf TT}-consequence with respect to a class of de Finetti algebrae). This is done in the next two definitions.

\begin{defn}
\label{evaluations}
Let $\mathcal{D}$ be a de Finetti algebra with support $\overline{D}$. A $\mathcal{D}$-evaluation is a function $e : {\sf For} \longmapsto \overline{D}$ s.t.:
\begin{align*}
e(\neg A) &\, := - e(A)\\
e(A \wedge B) &\, := e(A) \sqcap e(B)\\
e(A \vee B) &\, := e(A) \sqcup e(B)\\
e(A \rightarrow B) &\, := e(A) \dfcond e(B)
\end{align*}
\end{defn}

\begin{defn}
\label{A-consequence}
For every class of de Finetti algebrae $\mathfrak{D}$, every de Finetti algebra $\mathcal{D} \in \mathfrak{D}$, and every set $\{\Gamma, A\} \subseteq {\sf For}$:
\begin{itemize}
\item[] $A$ is a \emph{$\mathcal{D}$-consequence} of $\Gamma$, in symbols $\Gamma \models_{\mathcal{D}} A$ if for every $\mathcal{D}$-evaluation $e$, if for every $B \in \Gamma$,  $e(B) = {\bf 1}$ or $\nicefrac{\bf 1}{\bf 2}$, then $e(A) = {\bf 1}$ or $\nicefrac{\bf 1}{\bf 2}$. 
\item[] $A$ is a \emph{$\mathfrak{D}$-consequence} of $\Gamma$, in symbols $\Gamma \models_{\mathfrak{D}} A$, if for every $\mathcal{D} \in \mathfrak{D}$, $A$ is a $\mathcal{D}$-consequence of $\Gamma$. 
\end{itemize}
\end{defn}

Notice that, even though de Finetti algebrae include an algebraic counterpart of the \L{}ukasiewicz trivalent conditional, the latter is not used in defining an algebraic evaluation for de Finetti algebrae (and it is not going to be used to construct specific algebraic models of  {\sf DF/TT} either). The reason behind this choice is that we want to isolate the de Finetti conditional, and the respective {\sf TT}-logic, without including extraneous connectives (such as the \L{}ukasiewicz conditional). However, it would be possible to expand our definition of algebraic evaluations and algebraic consequence to include the \L{}ukasiewicz trivalent conditional, and prove the relative algebraic soundness and completeness theorems by adding suitable multi-sequent rules to {\sf DF/TTm}. 

We can finally use the structures introduced above, as well as the algebraic notions of consequence, to establish an algebraic soundness and completeness result. This is done in the next two results. 

\begin{lem}
\label{LTDlemma}
For every set $\{\Gamma, A\} \subseteq {\sf For}$:
\[\Gamma \vdash_{\sf DF/TTm} A \hspace{5pt} \mbox{ if and only if } \hspace{5pt} \Gamma \models_{\mathcal{D}(\Gamma)} A\]
\end{lem}

\begin{proof}[Proof sketch]
The left-to-right direction is proven by induction on the length of derivations in {\sf DF/TTm}. As for the right-to-left direction, suppose that $\Gamma \not\vdash_{\sf DF/TTm} A$. By Proposition \ref{P:derivation-countermodel}, the sequent $\Gamma \; | \; A \; | \; A$ has a countermodel, that is, there is a function $v : {\sf For} \longmapsto \{0, \nicefrac{1}{2}, 1\}$ s.t. for every $B \in \Gamma$, $v(B) = \nicefrac{1}{2}$ or $1$, but $v(A) = 0$.\footnote{Observe that the range of $v$ is the usual value space of de Finetti evaluations, that is, the set $\{0, \nicefrac{1}{2}, 1\}$, not the set of designated elements ${\bf 0}_{\Gamma}$, $\nicefrac{\bf 1}{\bf 2}_{\Gamma}$, and ${\bf 1}_{\Gamma}$ of $\mathcal{D}(\Gamma)$.} Let $P$ be the set of propositional variables in $\{\Gamma, A\}$. Let $\widetilde{e}$ be the partial function $\widetilde{e} : {\sf For} \rightharpoonup {\sf For}/\hspace{-4pt}\sim^{\sf df}_{\Gamma}$ defined as follows: 
\[\widetilde{e}(p_i) = \begin{cases}
\mathbf{0}_{\Gamma}, \mbox{ if } p_i \in P \mbox{ and } v(p_i) = 0\\
\nicefrac{\mathbf{1}}{\mathbf{2}}_{\Gamma}, \mbox{ if } p_i \in P \mbox{ and } v(p_i) = \nicefrac{1}{2}\\
\mathbf{1}_{\Gamma}, \mbox{ if } p_i \in P \mbox{ and } v(p_i) = 1\\
\end{cases}\]
By Zorn's Lemma, $\widetilde{e}$ can be expanded to a total function $e$ obeying the clauses of definition \ref{evaluations}, i.e. a {\sf DF}-evaluation.\footnote{More specifically, one first observes that $\widetilde{e}$ is a quasi-{\sf DF}-evaluation. It's easy to check that the collection of quasi-{\sf DF}-evaluations and {\sf DF}-evaluations forms a partially ordered set $P$, induced by the inclusion relation ($\subseteq$). Moreover, every totally ordered subset $Q$ of $P$ contains an upper bound, that is a quasi-{\sf DF}-evaluation or a {\sf DF}-evaluation which is not extended by any other element in $Q$. By Zorn's Lemma, then, $P$ itself contains a maximal element, i.e. a quasi-{\sf DF}-evaluations or a {\sf DF}-evaluation $e$ s.t. $\widetilde{e} \subseteq e$. It is then immediate to observe than $e$ is a {\sf DF}-evaluation and not a quasi-{\sf DF}-evaluation, and hence that $\widetilde{e} \subsetneqq e$, for otherwise it would be itself expanded by another element in $P$, and therefore it would not be maximal.} By Lemma \ref{LTDeFinetti}, $e$ is a $\mathcal{D}(\Gamma)$-evaluation, and by construction for every $B \in \Gamma$, $e(B) = \mathbf{1}_{\Gamma}$ or $\nicefrac{\mathbf{1}}{\mathbf{2}}_{\Gamma}$, but $e(A) = \mathbf{0}_{\Gamma}$. This shows that $\Gamma \not\models_{\mathcal{D}(\Gamma)} A$, as desired. 
\end{proof}

\begin{prop}[Algebraic soundness and completeness]
Let $\mathfrak{D}$ be the class of all the de Finetti algebrae. For every set $\{\Gamma, A\}\subseteq {\sf For}$: 
\[\Gamma \vdash_{\sf DF/TTm} A \hspace{5pt} \mbox{ if and only if } \hspace{5pt} \Gamma \models_{\mathfrak{D}} A\]\label{dfalg}
\end{prop}	

\begin{proof}[Proof sketch]
The left-to-right direction is straightforward. As for the right-to-left direction, suppose that $\Gamma \not\vdash_{\sf DF/TTm} A$. By Lemma \ref{LTDlemma}, this entails that $\Gamma \not\models_{\mathcal{D}(\Gamma)} A$, which in turn entails that $\Gamma \not\models_{\mathfrak{D}} A$, as desired. 
\end{proof}

It should be noted that the proof of algebraic completeness just given is not, strictly speaking, a genuine algebraic proof: it is parasitic on the Sch\"{u}tte-style completeness proof given in the previous Subsection \ref{S:2.2}.\footnote{See the proof of Proposition \ref{P:derivation-countermodel}.} More precisely, the Sch\"{u}tte-style proof is used to construct a countermodel based on the de Finetti algebra with just three elements, $\mathbf{0}$, $\nicefrac{\bf 1}{\bf 2}$, and $\mathbf{1}$, which is then expanded to an evaluation based on $\mathcal{D}(\Gamma)$. A typical algebraic proof would proceed by establishing a \emph{canonical model theorem}; however, {\sf DF/TT} does not seem to support this result. We can explain the specific features of the completeness theorem for {\sf DF/TT} as follows.\footnote{We thank an anonymous referee for prompting us to comment on this issue.} 

Lemma \ref{LTDlemma} and Proposition \ref{dfalg} entail that $\Gamma \models_{\mathcal{D}(\Gamma)} A$ if and only if $\Gamma \models_{\mathfrak{D}} A$. While the right-to-left direction of this biconditional is not surprising,\footnote{Lemma \ref{LTDeFinetti} establishes that $\mathcal{D}(\Gamma)$ is a de Finetti algebra, and hence $\mathcal{D}(\Gamma) \in \mathfrak{D}$, as in Proposition \ref{dfalg} we explicitly assume that $\mathfrak{D}$ is the class of all de Finetti algebrae.} its left-to-right direction yields that $A$ follows from $\Gamma$ in the specific de Finetti algebra $\mathcal{D}(\Gamma)$ just in case $A$ follows from $\Gamma$ in the class of all de Finetti algebrae. That is, the question of whether $A$ is a consequence of $\Gamma$ in all de Finetti algebrae is reduced to the question of whether $A$ is a consequence of $\Gamma$ in $\mathcal{D}(\Gamma)$. However, this appears less surprising if one considers how $\mathcal{D}(\Gamma)$ is constructed: in fact, $\mathcal{D}(\Gamma)$ is built using the very sentences that make up the target inference from $\Gamma$ to $A$.\footnote{There are immediate parallels outside of abstract algebraic logic if one looks at how Henkin models are defined in standard completeness proofs for pure classical logic, or if one considers canonical models in completeness proofs for normal modal logics.} Indeed, algebraic completeness proofs typically establish even stronger results, i.e. they establish a canonical model theorem. Such a theorem guarantees that there is \emph{exactly one} evaluation $e$ (in our case, it would be a de Finetti (algebraic) evaluation) which assigns the designated value to all the elements of $\Gamma$ such that $e$ assigns the designated value to $A$ if and only if $\Gamma \models_{\mathfrak{D}} A$. That is, rather than reducing the question of whether $A$ is a consequence of $\Gamma$ in all algebrae of a certain kind (again, in our case it would be de Finetti algebrae) to the question of whether $A$ is a consequence of $\Gamma$ in a single algebra of that kind, a canonical model theorem reduces the former question to the question of whether $A$ has a designated value in a single evaluation. And this is a much stronger result because, in order to establish that $\Gamma \models_{\mathcal{D}(\Gamma)} A$, one has to consider \emph{all} the de Finetti algebraic evaluations based on $\mathcal{D}(A)$, whereas if a canonical model were available, one would just need to consider \emph{one} such evaluation (a particularly informative and `canonical' one). This happens, for example, in proofs of algebraic completeness for classical logic (with respect to the class of all Boolean algebrae), and for several other logics as well. Many of these logics are also algebraizable (which is, typically, an even stronger result). 

One can visualize the situation as follows: 

\begin{center}
\begin{figure}[h!]
\begin{tabular}{ccccccc}
\begin{tikzpicture}[scale=0.3]
    \draw (0,0) circle (3cm);
    \clip (0,0) circle (3cm);
    \pgfmathsetseed{10294453}
    \foreach \p in {1,...,500}
    { \fill[black]  (5*rand,3*rand) circle (0.05);} 
\end{tikzpicture}
&
\begin{tikzpicture}[scale=0.3]
    \draw (0,0) circle (3cm);
    \clip (0,0) circle (3cm);
    \pgfmathsetseed{54356782}
    \foreach \p in {1,...,500}
    { \fill[black]  (5*rand,3*rand) circle (0.05);} 
\end{tikzpicture}
&
\begin{tikzpicture}[scale=0.3]
    \draw (0,0) circle (3cm);
    \clip (0,0) circle (3cm);
    \pgfmathsetseed{78493029}
    \foreach \p in {1,...,500}
    { \fill[black]  (5*rand,3*rand) circle (0.05);} 
\end{tikzpicture}
&
\begin{tikzpicture}[scale=0.3]
    \draw (0,0) circle (3cm);
    \clip (0,0) circle (3cm);
    \pgfmathsetseed{32139472}
    \foreach \p in {1,...,500}
    { \fill[black]  (5*rand,3*rand) circle (0.05);} 
\end{tikzpicture}
&
\begin{tikzpicture}[scale=0.3]
    \draw (0,0) circle (3cm);
    \clip (0,0) circle (3cm);
    \pgfmathsetseed{209384560}
    \foreach \p in {1,...,500}
    { \fill[black]  (5*rand,3*rand) circle (0.05);} 
\end{tikzpicture}
&
\begin{tikzpicture}[scale=0.3]
    \draw (0,0) circle (3cm);
    \clip (0,0) circle (3cm);
    \pgfmathsetseed{24638594}
    \foreach \p in {1,...,500}
    { \fill[black]  (5*rand,3*rand) circle (0.05);}
\end{tikzpicture}
&
$\cdots$\\
\begin{tikzpicture}[scale=0.3]
    \draw (0,0) circle (3cm);
    \clip (0,0) circle (3cm);
    \pgfmathsetseed{53222289}
    \foreach \p in {1,...,500}
    { \fill[black]  (5*rand,3*rand) circle (0.05);
    }
\end{tikzpicture}
&
\begin{tikzpicture}[scale=0.3]
    \draw (0,0) circle (3cm);
    \clip (0,0) circle (3cm);
    \pgfmathsetseed{77765120}
    \foreach \p in {1,...,500}
    { \fill[black]  (5*rand,3*rand) circle (0.05);
    }
\end{tikzpicture}
&
\begin{tikzpicture}[scale=0.3]
    \draw (0,0) circle (3cm);
    \clip (0,0) circle (3cm);
    \pgfmathsetseed{62398475}
    \foreach \p in {1,...,500}
    { \fill[black]  (5*rand,3*rand) circle (0.05);} 
\end{tikzpicture}
&
\begin{tikzpicture}[scale=0.3]
    \draw (0,0) circle (3cm);
    \clip (0,0) circle (3cm);
    \pgfmathsetseed{19234756}
    \foreach \p in {1,...,500}
    { \fill[black]  (5*rand,3*rand) circle (0.05);} 
\end{tikzpicture}
&
\begin{tikzpicture}[scale=0.3]
    \draw (0,0) circle (3cm);
    \clip (0,0) circle (3cm);
    \pgfmathsetseed{39847561}
    \foreach \p in {1,...,500}
    { \fill[black]  (5*rand,3*rand) circle (0.05);} 
\end{tikzpicture}
&
\begin{tikzpicture}[scale=0.3]
    \draw (0,0) circle (3cm);
    \clip (0,0) circle (3cm);
    \pgfmathsetseed{92837465}
    \foreach \p in {1,...,500}
    { \fill[black]  (5*rand,3*rand) circle (0.05);} 
\end{tikzpicture}
&
$\cdots$\\
$\vdots$ & $\vdots$ & $\vdots$ & $\vdots$ & $\vdots$ & $\vdots$ &
\end{tabular}
\end{figure}
\end{center}
\noindent In the above representation, a circle is an algebra, a dot in a circle is an evaluation based on that algebra, and the set of circles is a class of algebrae. Our completeness result establishes that, for every inference, there is one circle such that that inference is valid in all the circles if and only if it is valid in that circle (which in turn requires considering all the dots in that circle). A canonical model theorem establish a much stronger claim, namely that for every inference there is \emph{one single dot} in one of the circles such that that inference holds in all the circles if and only if its conclusion receives the designated value in that dot. 

Crucially, however, proving a canonical model theorem typically requires Modus Ponens. Indeed, an attempted proof of a canonical model theorem for {\sf DF/TT} breaks down exactly where Modus Ponens is required. In absence of Modus Ponens, only the weaker result that $\Gamma \models_{\mathcal{D}(\Gamma)} A$ if and only if $\Gamma \models_{\mathfrak{D}} A$ holds. However, this is enough for our purposes, that is proving algebraic completeness for {\sf DF/TT}. In the next subsection, we will see that things do not work so well for {\sf CC/TT} and the {\sf TT}-logics of Jeffrey conditionals more generally. 

We conclude this subsection with a brief discussion of the prospects for a full \emph{algebraizability} of {\sf DF/TT}.\footnote{The notion of algebraizability \citep[introduced in][]{BlokPigozzi1989} generalizes the link between a logic and its algebraic semantics, imposing stricter conditions than those required for algebraic completeness. See also \citet{Herrmann1996} and \citet[][Chapters 2 and 3]{Font2016}.} 
In this respect, however, {\sf DF/TT} reveals non-negligible limitations. Let's start with some preliminary definitions.

\begin{defn}$\,$
\begin{itemize}
\item[-] Let an \emph{equation}, in symbols $A \approx B$, be a pair of formulae of $\mathcal{L}_{\rightarrow}$. The set of all equations is therefore ${\sf For} \times {\sf For}$. A set of equations, i.e. a subset of ${\sf For} \times {\sf For}$, is denoted by $\mathbb{E}$. 
\item[-] Let $\mathcal{D}$ be a de Finetti algebra, and let $\{\mathbb{E}, C \approx D\}$ be a set of equations. $\mathcal{D}$ \emph{validates the equational inference from} $\mathbb{E}$ \emph{to} $C \approx D$, in symbols $\mathbb{E} \models_{\mathcal{D}} C \approx D$, if for every $\mathcal{D}$-evaluation $e$ if, if for every equation $(A \approx B) \in \mathbb{E}$, $e(A) = e(B)$, then $e(C) = e(D)$.\footnote{Informally, $\mathcal{D}$ validates the equational inference from $\mathbb{E}$ to $C \approx D$ if every de Finetti evaluation based on $\mathcal{D}$ that satisfies all the equations in $\mathbb{E}$, also satisfies $C \approx D$.}
\item[-] Let $\mathfrak{D}$ be a class of de Finetti algebrae. $\mathfrak{D}$ \emph{validates the equational inference from} $\mathbb{E}$ \emph{to} $C \approx D$, in symbols $\mathbb{E} \models_{\mathfrak{D}} A \approx B$, if for every $\mathcal{D} \in \mathfrak{D}$, $\mathbb{E} \models_{\mathcal{D}} A \approx B$.
\end{itemize}
\end{defn}

\noindent We can now formulate the notion of algebraizability.

\begin{defn}
A logic ${\sf L}$ is \emph{algebraizable} if there are a class of algebrae $\mathfrak{A}$ and functions $f : \mathcal{P}({\sf For}) \longmapsto \mathcal{P}({\sf For} \times {\sf For})$ and $g : \mathcal{P}({\sf For} \times {\sf For}) \longmapsto \mathcal{P}({\sf For})$ s.t. for every $\{\Gamma, A\} \subseteq {\sf For}$:\footnote{Let `$f[\Gamma]$' be a shorthand for $\{f(B) \in \mathcal{P}({\sf For} \times {\sf For}) \,|\, B \in \Gamma\}$.}
\begin{itemize}
\item[\emph{(A1)}] $\Gamma \vdash_{\sf L} A$ if and only if $f[\Gamma] \models_{\mathfrak{A}} f(A)$
\item[\emph{(A2)}] $B \approx C \hspace{2pt}\rotatebox[origin=c]{180}{$\models$} \hspace{-2.5pt}\models_{\mathfrak{A}} f(g(B \approx C))$
\end{itemize}
\label{D:algebraizability}
\end{defn}

(A1) is a generalization of algebraic completeness, where the right-hand side expresses in the object-language the requirement that $A$ has a designated value whenever all the sentences in $\Gamma$ do.\footnote{For example, consider the case of classical logic and Boolean algebrae. Letting $f$ be the function defined as $f(A) = A \approx \top$, and $\mathfrak{B}$ be the class of Boolean algebrae, $f[\Gamma] \models_{\mathfrak{B}} f(A)$ becomes $\bigcup_{B\, \in \,\Gamma} B \approx \top \models_{\mathfrak{B}} A \approx \top$,
which formalizes the idea that whenever every sentence in $\Gamma$ has value $\mathbf{1}$, so does $A$.}
(A2) ensures that the solvability of equations
is fully captured by some formula of the object-language.\footnote{Consider again classical logic and Boolean algebrae. Classically, $B$ and $C$ have the same value just in case $B \equiv_{\it m} C$ (where $\equiv_{\it m}$ denotes the material biconditional) has value $1$. Letting $g$ be s.t. $g(\mathbb{E}) = \{B \equiv_{\it m} C \in {\sf For} \,|\, B \approx C \in \mathbb{E}\}$, $B \approx C \hspace{2pt}\rotatebox[origin=c]{180}{$\models$} \hspace{-2.5pt}\models_{\mathfrak{B}} f(g(B \approx C))$ becomes $B \approx C \hspace{2pt}\rotatebox[origin=c]{180}{$\models$} \hspace{-2.5pt}\models_{\mathfrak{B}} (B \equiv_{\it m} C) \approx \top$, 
which express the idea that $B$ and $C$ have the same value whenever the corresponding biconditional holds (i.e. has value $1$) in classical logic.}

However, (A1) and (A2) do not sit well with the conditional of {\sf DF/TT}. In the context of {\sf DF/TT}, the right-hand side of (A1) expresses that whenever all the sentences in $\Gamma$ have value $\mathbf{1}$ or $\nicefrac{\mathbf{1}}{\mathbf{2}}$, so does $A$. A plausible formalization is obtained by letting $f$ be $f(A) = A \approx (A \leftrightarrow \top) \vee (A \leftrightarrow (\bot \rightarrow \top))$. 
$f[\Gamma] \models_{\mathfrak{D}} f(A)$ then becomes 
\[\bigcup_{B \,\in\, \Gamma} B \approx (B \leftrightarrow \top) \vee (B \leftrightarrow (\bot \rightarrow \top)) \models_{\mathfrak{D}} A \approx (A \leftrightarrow \top) \vee (A \leftrightarrow (\bot \rightarrow \top)),\]
Define $g$ as $g(\mathbb{E}) = \{B \leftrightarrow C \in {\sf For} \,|\, B \approx C \in \mathbb{E}\}$. (A2) becomes 
\begin{align*}
B \approx C \hspace{2pt}\rotatebox[origin=c]{180}{$\models$} \hspace{-2.5pt}\models_{\mathfrak{D}} (B \leftrightarrow C) \approx ((B \leftrightarrow C) \leftrightarrow \top) \vee ((B \leftrightarrow C) \leftrightarrow (\bot \rightarrow \top))
\label{E:second}
\end{align*}
which, however, does \emph{not} express the idea that $B$ and $C$ have the same value whenever $B \leftrightarrow C$ holds in {\sf DF/TT}. In fact, it is \emph{not} the case that, in order for $B \leftrightarrow C$ to have the same value as $((B \leftrightarrow C) \leftrightarrow \top) \vee ((B \leftrightarrow C)$, $B$ has to have the same value as $C$; an evaluation $e$ such that $e(B) = 1$ and $e(C) = \nicefrac{1}{2}$ provides a counterexample. This translates into the algebraic semantics, considering a de Finetti algebra $\mathcal{D}$ and an algebraic evaluation $e$ based on $\mathcal{D}$ s.t. $e(B) = \mathbf{1}$ and $e(C) = \nicefrac{\mathbf{1}}{\mathbf{2}}$ (for $\mathbf{1}$, $\nicefrac{\mathbf{1}}{\mathbf{2}} \in \overline{D}$). 

Of course, this observation only tells us that $B \leftrightarrow C$ does not express the fact that $B$ and $C$ have the same value---which is not surprising, given the $1$- and $\nicefrac{1}{2}$-rows of the truth table of the de Finetti conditional. However, the idea of mapping identity of semantic values to a formula that expresses `having a designated value' seems at odds with the conditional of {\sf DF/TT}, because {\sf DF/TT} does not distinguish between $1$ and $\nicefrac{1}{2}$ when it comes to designatedness, nor does its conditional. 
In a tolerant-tolerant semantics, validity (and the corresponding conditionals) does not depend on the {identity} of the semantic values that are preserved from premises to conclusion, but on their {similarity}: $1$ and $\nicefrac{1}{2}$ are not identical, but {similar enough} for ${\sf DF/TT}$ not to distinguish them. But condition (A2) can only be satisfied via a formula that captures a notion of validity based on the identity of semantic values. 

This does not show that {\sf DF/TT} is not algebraizable: here, we leave the question open. However, the above observations suggest a possible  strategy to prove non-algebraizability: if one can show that a truth-function expressing identity of truth values is not definable in the truth table semantics for {\sf DF}, this would translate into the algebraic semantics, and establish non-algebraizability.\footnote{This proof strategy seems simpler and more informative than a proof via Isomorphism Theorems, which are the standard results employed to prove non-algebraizability \citep[see][Chapter 3.5]{Font2016}. Another open question is whether {\sf DF/TT} is algebraizable over other logics. A natural choice would be \L{}ukasiewicz trivalent logic (with a {\sf TT}-notion of validity), because de Finetti algebrae are defined over $\L{}3$ algebrae. There are reasons to expect a positive result, namely that the \L{}ukasiewicz trivalent conditional can be used to express the identity of semantic values, e.g. via $\neg ((A \leftrightarrow B) \rightarrow \neg (A \leftrightarrow B))$. 
}

\subsection{An algebraic semantics for {\sf CC/TT}?}
\label{S:algebraic-sound-compl-CC}

\noindent Can we provide a proof of algebraic completeness for {\sf CC/TT} employing the Lindenbaum-Tarski method, as we did for {\sf DF/TT}? The {\sf CC/TT}-conditional appears better-behaved than the {\sf DF/TT} one---in particular because it obeys Modus Ponens---so this would appear prima facie possible. 

Let us try to apply the Lindenbaum-Tarski method to {\sf CC/TT}. First, we need an algebraic counterpart of the Cooper-Cantwell conditional. This is provided by the following definition. 

\begin{defn}
An $\L{}3$ algebra $\mathcal{A} = \langle \overline{A}, \sqcap, \sqcup, -, \rightarrowtail, \ccond, \mathbf{0}, \mathbf{1}, \nicefrac{{\bf 1}}{{\bf 2}} \rangle$ is \emph{Cooper-Cantwell} if: 
\begin{itemize}
\item[] There is a distinguished element $\nicefrac{{\bf 1}}{{\bf 2}} \in \overline{A}$ s.t. $- \nicefrac{{\bf 1}}{{\bf 2}} = \nicefrac{{\bf 1}}{{\bf 2}}$, and 
\item[] There is an operation $\ccond$ defined on $\overline{A} \times \overline{A}$ s.t. $a \ccond b = -{\sf w}(a) \sqcup ({\sf w}(a) \sqcap b)$, where ${\sf w}(a)$ is a shorthand for $- a \rightarrowtail \nicefrac{{\bf 1}}{{\bf 2}}$.
\end{itemize}
\end{defn}

We then work towards the construction of a Lindenbaum-Tarski algebra for {\sf CC/TT}. 

\begin{defn}
For every $\Gamma \subseteq {\sf For}$, let $\sim^{\sf c}_{\Gamma} \, \subseteq {\sf For} \times {\sf For}$ be the relation defined as follows: 
\[A \sim^{\sf c}_{\Gamma} B \mbox{ if and only if } \Gamma \vdash_{\sf CC/TTm} A \leftrightarrow B\]
\end{defn}

Since Modus Ponens holds in {\sf CC/TT}, the relation of {\sf CC/TTm}-provable equivalence seems better behaved than the one defined for {\sf DF/TTm}.

\begin{lem}
\label{equivalence}
$\sim^{\sf c}_{\Gamma}$ is an equivalence relation on ${\sf For} \times {\sf For}$. 
\end{lem}

\begin{proof}
Reflexivity holds since $\Gamma \vdash_{\sf CC/TTm} A \leftrightarrow A$. Symmetry also holds, since if $\Gamma \vdash_{\sf CC/TTm} A \leftrightarrow B$, then also $\Gamma \vdash_{\sf CC/TTm} B \leftrightarrow A$. Finally, transitivity holds as well, because if $\Gamma \vdash_{\sf CC/TTm} A \leftrightarrow B$, and $\Gamma \vdash_{\sf CC/TTm} B \leftrightarrow C$ then $\Gamma \vdash_{\sf CC/TTm} A \leftrightarrow C$ as well. These equivalences are quickly established semantically, i.e. considering $\models_{\sf CC/TTm}$ rather than $\vdash_{\sf CC/TTm}$, by the completeness of {\sf CC/TTm}. 
\end{proof}

We now have an equivalence relation, so we can use it to partition the set of formulae into equivalence classes. 

\begin{defn}
For every $\{\Gamma, A\} \subseteq {\sf For}$, let $[A]^{\sf c}_{\Gamma}$ denote the \emph{equivalence class of $A$ induced by $\sim^{\sf c}_{\Gamma}$}. 
The \emph{quotient} induced by $\sim^{\sf c}_{\Gamma}$ on ${\sf For}$, in symbols ${\sf For} \slash\hspace{-4pt}\sim^{\sf c}_{\Gamma}$, is the set of equivalence classes induced by $\sim^{\sf c}_{\Gamma}$. 
\end{defn}

\noindent Since we only work with Cooper-Cantwell algebrae in this subsection, we drop the superscript $^{\sf cc}$ again to improve readability, without risks of confusion. Now, in order to proceed with the proof of algebraic completeness, we would have to define a Cooper-Cantwell version of a Lindenbaum-Tarski algebra. Such a structure would look as follows: 

\[\mathcal{C}(\Gamma) = \langle {\sf For} \slash\hspace{-4pt}\sim^{\sf c}_{\Gamma}, \sqcap_{\Gamma}, \sqcup_{\Gamma}, -_{\Gamma}, \ccond_{\Gamma}, \mathbf{0}_{\Gamma}, \mathbf{1}_{\Gamma} \rangle\]
where: 
\begin{center}
\begin{minipage}{6cm}
\begin{itemize}
\item[] $[A]_{\Gamma} \sqcap_{\Gamma} [B]_{\Gamma} := [A \wedge B]_{\Gamma}$
\item[] $-_{\Gamma} [A]_{\Gamma} := [\neg A]_{\Gamma}$
\item[] $[\bot]_{\Gamma} := \mathbf{0}_{\Gamma}$
\end{itemize}
\end{minipage} 
\begin{minipage}{6cm}
\begin{itemize}
\item[] $[A]_{\Gamma} \sqcup_{\Gamma} [B]_{\Gamma} := [A \vee B]_{\Gamma}$
\item[] $ [A]_{\Gamma} \ccond_{\Gamma} [B]_{\Gamma} := [A \rightarrow B]_{\Gamma}$
\item[] $[\top]_{\Gamma} := \mathbf{1}_{\Gamma}$
\end{itemize}
\end{minipage} 
\[[\bot \rightarrow \top]_{\Gamma} := \nicefrac{{\bf 1}}{{\bf 2}}_{\Gamma}\]
\end{center}

However, the construction is blocked, because some of its defining operations turn out to be \emph{not well-defined}. In particular, the Cooper-Cantwell conditional is not substitutive with respect to negation. 

\begin{lem}\label{lem:CCsub}
There are sets $\{\Gamma, A, B\} \subseteq {\sf For}$ s.t. 
\[A \sim_{\Gamma} B \hspace{5pt} \mbox{ but it is not the case that } \hspace{5pt} \neg A \sim_{\Gamma} \neg B\] 
\end{lem}

\begin{proof}
It is sufficient to set $\Gamma = \varnothing$, $A = \top$, and $B = \bot \rightarrow \top$. 
\end{proof}

\noindent This lemma shows that the process of providing an algebraic semantics (via the standard Lindenbaum-Tarski method) for {\sf CC/TT} stops here: it does not even get off the ground. 

In fact, this negative result is more general: it applies to \emph{every} Jeffrey conditional. Recall that Jeffrey conditionals are required to obey the condition that $f_{\to}(1,0) = f_{\to}(\nicefrac{\bf 1}{\bf 2},0) = 0$. Now, the above proof employs exactly the cases in which a conditional has an antecedent with value $1$ and a consequent with value $0$, and an antecedent with value $1$ and a consequent with value $\nicefrac{\bf 1}{\bf 2}$. Therefore, no Jeffrey conditional is substitutive with respect to negation---under a {\sf TT}-notion of validity, and a Strong Kleene interpretation of conjunction and negation. In turn, this means that no `{\sf J}-Lindenbaum-Tarski algebra', where `{\sf J}' is any Jeffrey conditional, is well-defined, and therefore that no algebraic semantics (via the Lindenbaum-Tarski method) is available for any {\sf TT}-logic of a Jeffrey conditional.

\section{General Discussion}
\label{S:Discussion2}

This two-part paper has reviewed the main motivations for a trivalent semantics for indicative conditionals, interpreting them as \textit{conditional assertions}, and defining their truth conditions in analogy with the conditions that settle the winner of a conditional bet (i.e., the bet or assertion is declared void when the antecedent is false). Although the idea goes back to \citet{definetti1936logique}, and \citet{reichenbach1935wahr, reichenbach1944philosophic}, there have been few explorations of the 
logics induced by the adoption of that semantic scheme. 
Beside expounding the historical roots of trivalent semantics for conditionals, our paper has given a systematic survey of the different logics that emerge by (i) choosing a truth table for the conditional operator in agreement with the above rationale, and (ii) determining a specific notion of validity (one vs.~two designated truth values, pure vs.~mixed consequence relations). 

As reviewed in Part I, the trivalent approach yields a fully truth-functional semantics with attractive logical and inferential properties. It also provides the conceptual foundations for a probabilistic theory of assertability and reasoning with conditionals along the lines of \citet{Adams1975}. For simple conditionals, combining our semantics with defining the assertability of a sentence $A$ as the conditional probability that A is true, given that it has a classical truth value, immediately yields Adams' Thesis that $Ast(A \to C) = p(C|A)$. This property highlights the potential of the trivalent approach for guiding an account of the epistemology of conditionals, and explaining how people reason with them \citep[e.g.,][]{baratgin2013uncertainty,BaratginEtAl2018}. 
While the semantics of the trivalent conditional is \textit{factual}---that is, its truth value is a function of matters in the actual world---no such limits are imposed on the scope of the probability functions in judgments of assertability (e.g., $A$ can be practically unverifiable, but the conditional may still be highly assertable). 

With respect to the above challenges (i) and (ii), it quickly transpires that any alternative to a tolerant-to-tolerant ({\sf TT-}) notion of validity would be either too strong (in the sense of licensing undesirable inferences such as implying the converse condtional) or too weak (in the sense of violating the Identity Law $A \rightarrow A$ and not having sentential validities). Only the Cooper-Cantwell conditional, where indeterminate antecedents are exactly treated like true ones, satisfies both the full Deduction Theorem and commutation with negation. For conceptual, empirical and logical reasons (the conditional is essentially interpreted as making an assertion upon supposing the antecedent), these are eminently reasonable properties, apparently favoring \cctt\ as the best trivalent logic of the indicative conditional. 

The results of Part II nuance this judgment. For both \dftt\ and \cctt\, we can develop sound and complete calculi based on tableaux (Section \ref{S:tableaux}) and three-sided sequents (Section \ref{S:sequents}). The latter calculi have the advantage of being simpler and more direct: unlike tableau calculi, they do not establish that an inference is valid by showing that it is impossible to assign a designated value to the premises and an undesignated value to the conclusion. Moreover, many-sided sequent calculi make it easier to handle inferences with multiple conclusions, as well as inferences involving infinite sets of sentences.

As soon as we consider the algebraic semantics, however, differences between \dftt\ and \cctt\ emerge. While provable equivalence fails to be transitive and therefore induces no equivalence relation for \dftt\ (Lemma \ref{noequiv}), we can still use this relation to define a Lindenbaum-Tarski algebra and to show an algebraic soundness and completeness theorem (Proposition \ref{dfalg}). In other words, $A$ can be derived from $\Gamma$ using one of the above calculi (e.g., many-sided sequents) if and only if a consequence relation holds between $\Gamma$ and $A$ in the associated de Finetti algebrae. The failure of Modus Ponens for \dftt\, however, blocks the construction of a canonical algebraic model.

Things look bleak, by contrast, for \cctt\ and other {\sf TT-}logics based on a Jeffrey conditional. While provable equivalence induces an equivalence relation for these logics, the construction of a Lindenbaum-Tarski algebra does not get off the ground because provable equivalence fails to be substitutive under negation. More precisely, the Cooper-Cantwell biconditional $\leftrightarrow$ falls short of expressing \cctt-equivalence since $A \leftrightarrow B \not\models_{\sf CC/TT} \neg A \leftrightarrow \neg B$. Which means that there is not, and cannot be, a fruitful algebraic treatment of Jeffrey conditionals. In fact, this is grounded in a defining property of Jeffrey conditionals: to preserve Modus Ponens and to yield a full Deduction Theorem, a trivalent conditional based on the ``defective'' truth table needs to obey $f_{\to}(1,0) = f_{\to}(\half,0) = 0$. It is exactly this property which makes substitution under negation fail (Lemma \ref{lem:CCsub}), and prevents a proper algebraic semantics for Jeffrey conditionals. 

Clearly, the failure of substitution under negation is closely related to the failure of contraposition in Jeffrey conditionals---an inference that does not fail in {\sf DF/TT}. Indeed, the same evaluation provides the counterexamples employed in proving both Proposition 5.6 (Part I) and Lemma \ref{lem:CCsub} (Part II). 
So it turns out that what has been a strength of Jeffrey tolerant-tolerant logics, and \cctt\ in particular, at the level of desirable conditional principles, comes at the price of the  algebraic semantics. Importantly, the lack of an algebraic semantics is not a mere technical fact, but it has philosophical consequences as well. In particular, in every Jeffrey tolerant-tolerant logic, even if it is the case that $A \leftrightarrow B$, the same equivalence does not hold in general for logically complex sentences that result by uniform substitutions of $A$ and $B$ (see Lemma \ref{L:indipendence} for a formally precise version of this property). Therefore, Jeffrey conditionals do not provide a workable notion of equivalence.

Of course, the limitations of Jeffrey conditionals just reviewed arise from the combination of the semantics of Jeffrey conditionals, {\sf TT}-validity, and Strong Kleene conjunction, disjunction, and negation: one might therefore wonder whether they can be improved on by altering some of these parameters. However, as the results of Part I show, adopting an alternative to  {\sf TT}-validity does not seem promising. As for the semantics of the other connectives, in part I we noted that Cooper adopted alternative truth tables for conjunction and disjunction, while retaining the {\sf K3} table for negation \citep[see also][\S 7.19, 1044 and following]{humberstone2011conn}. Cooper's conjunction and disjunction, we noted, lose some classical properties: for instance, in Cooper's original system, one cannot in general infer $A \vee B$ from $A$, for when $A$ has value $\half$ and $B$ has value $0$, $A \vee B$ has value $0$ as well. With regard to the algebraic semantics, the difficulties raised in the previous section may not easily be overcome as a result.\footnote{A similar application of the Lindenbaum-Tarski method would also likely fail if one replaces strong Kleene conjunction and disjunction with Cooper's quasi-conjunction and quasi-disjunction (see Part I, Section 6), because it only depends on features of the conditional and negation (see Lemma \ref{lem:CCsub}). However, we don't have proper algebraic counterparts of Cooper's connectives, therefore we can only advance this claim as a conjecture. In particular, some features of Cooper's connectives show that they cannot receive their standard algebraic interpretation. Consider the failure of quasi-disjunction introduction. Standardly, a disjunction is interpreted as the maximum of the two disjuncts, in some order. In a bounded lattice, this means that if an element $a$ has a  position with respect to the top element, in the order induced by the lattice, then for every element $b$, $\max(a, b)$ has a position which is at least as close to the top element as the position of $a$. But this feature seemingly has to fail for an algebraic representation of Cooper's quasi-disjunction, for it corresponds to disjunction introduction. The failure of negated quasi-conjunction introduction (i.e. that the inference from $\neg A$ to $\neg (A \wedge B)$ is not {\sf TT}-valid for Cooper's quasi-conjunction) seems to have even more far-reaching consequences, for it also involves the involution defined on the corresponding lattices.}

Alternatively, one might inquire into what happens to the interaction of Jeffrey conditionals with a non-{\sf K3} negation. To be sure, the {\sf K3}-negation squares particularly well with the philosophical motivation for de Finettian conditionals: when a conditional assertion $A \to C$ is ``called off'' because $A$ is false, the same should happen for the negation of that assertion (i.e., the sentence $A \to \neg C$, thanks to the commutation scheme). Nonetheless, it might be worth investigating how Jeffrey logics (keeping a tolerant-tolerant notion of validity) fare when coupled with what \cite{chemla2018mvcons} call a ``Gentzen-regular'' negation, that is a negation obeying the Gentzen sequent calculus rules. While a Gentzen-regular negation might avoid some of the above problems, it would lose the commutation of conditional and negation, and the attached connexive principles (see Subsections 5.2 and 5.3 
of Part I). In conclusion, there seem to be structural limitations, or at least unavoidable tradeoffs, that affect Jeffrey conditionals, when it comes to their interaction with other connectives.

We therefore believe that it is not easy to justify a clear preference between the two logics \cctt\ and \dftt  that we have isolated as most promising amongst trivalent logics of indicative conditionals. Both have attractive properties, both have limitations---but they agree in essential properties such as the valuation of classical sentences, the Import-Export principle, the analysis of paradoxes of material implication, their connexive nature, and the connection to a theory of assertability. 
To solve the limitations highlighted throughout the paper, one would probably have to give up one or more of these features. 
So while there is perhaps no perfect trivalent semantics for indicative conditionals, they need to be considered carefully between two-valued logic and modal logics of conditionals.
In any event, they give rise to a promising research program, and we shall support this claim by sketching some future projects that build on our work in this paper. 

Firstly, we would like to extend the current framework to predicate logic and to investigate how the trivalent conditionals fare in that context, including how they interact with a na\"{i}ve or a compositional truth predicate. Secondly, we would like to apply trivalent semantics to McGee's famous challenge to Modus Ponens, applying our accounts of logical consequence and probabilistic assertability \citep{McGee1985,SternHartmann2018}. Thirdly, one should review the intuitions and inference schemes which fuel connexive logics (e.g., Aristotle's Thesis, Boethius' Thesis) from a trivalent perspective, and conduct a more detailed comparison. Finally, we need to develop more precise criteria as to which inferences should be validated by a trivalent \textit{logic} of conditionals, based on the concept of supposition, and which inferences can be relegated to a probabilistic theory of \textit{assertability} grounded in the truth conditions. In other words, we have to formulate a precise account of how the \textit{truth conditions} of indicative conditionals relate to \textit{reasoning} with them.


\pagebreak
\bibliographystyle{../BibFolder/mychicago-ff}

\renewcommand{\bibname}{References}
\bibliography{../BibFolder/v-c(14)}